\documentclass[12pt]{amsart}
\usepackage{amsmath,amssymb,amscd,array, mathrsfs }

\usepackage{float}

\def\url#1{\expandafter\string\csname #1\endcsname}

\def\mmat #1,#2,#3,#4,{\text{\small\arraycolsep=3pt $
\begin{pmatrix}#1&#2\\#3&#4\end{pmatrix}$}}

\usepackage{hyperref}

\usepackage{comments}
\newComments\SBe{Said}{blue}
\newComments\SBo{Sofiane}{blue}
\newComments\YM{Yoshi}{blue}
\newComments\DL{DL}{red}

\def\mmat #1,#2,#3,#4,{\text{\small\arraycolsep=3pt $
\begin{pmatrix}#1&#2\\#3&#4\end{pmatrix}$}}

\usepackage{lscape}
\usepackage{tikz-cd}

\usepackage{DLdef1}

\def\mmat #1,#2,#3,#4,{\text{\small\arraycolsep=3pt $
\begin{pmatrix}#1&#2\\#3&#4\end{pmatrix}$}}

\renewcommand {\ssbegin}[2][*]
 {\refstepcounter{subsection}%
\if#1*
\addcontentsline{toc}{subsection}{\thesubsection.\hskip 1pc #2}%
\else
\addcontentsline{toc}{subsection}{\thesubsection.\hskip 1pc #2. #1}%
\fi
 \def \secno {\gdef \secno {}{\ssecfont
\thesubsection.\hskip 2ex}%
 }%
 \begin{#2}}

\renewcommand {\sssbegin}[2][*]
 {\refstepcounter{subsubsection}
\if#1*
\addcontentsline{toc}{subsubsection}{\thesubsubsection.\hskip 1pc #2}%
\else
\addcontentsline{toc}{subsubsection}{\thesubsubsection.\hskip 1pc #2. #1}
\fi
 \def \secno {\gdef \secno {}{\ssecfont \thesubsubsection.\hskip 2ex}%
 }%
 \begin{#2}}

\renewcommand {\parbegin}[2][*]
 {\refstepcounter{paragraph}
\if#1*
\addcontentsline{toc}{paragraph}{\theparagraph.\hskip 1pc #2}%
\else
\addcontentsline{toc}{paragraph}{\theparagraph.\hskip 1pc #2. #1}
\fi
 \def \secno {\gdef \secno {}{\ssecfont \theparagraph.\hskip 2ex}%
 }%
 \begin{#2}}

\rmnameii{etr}{etr}
\rmnameii{evv}{ev}

\begin{document}
\title[Double and Lagrangian extensions for quasi-Frobenius Lie superalgebras]{Double and Lagrangian extensions for quasi-Frobenius Lie superalgebras}

\author{Sofiane Bouarroudj}

\address{Division of Science and Mathematics, New York University Abu Dhabi, P.O. Box 129188, Abu Dhabi, United Arab Emirates.}
\email{sofiane.bouarroudj@nyu.edu}

\author{Yoshiaki Maeda}
\address{Tohoku Forum for Creativity, Tohoku University, 
   2-1-1, Katahira, Aoba-ku, Sendai, Japan}
\email{ymkeiomath@gmail.com}

\thanks{SB was supported by the grant NYUAD-065.}


\date{\today}



\begin{abstract} 

A Lie superalgebra is called {\it quasi-Frobenius} if it admits a closed anti-symmetric non-degenerate bilinear form.
 
We study the notion of double extensions of quasi-Frobenius Lie superalgebra when the form is either orthosymplectic or periplectic.  We show that every quasi-Frobenius Lie superalgebra that satisfies certain conditions can be obtained as a double extension of a smaller quasi-Frobenius Lie superalgebra.

We classify all 4-dimensional quasi-Frobenius Lie superalgebras,  and show that such Lie superalgebras must be solvable. 

We study the notion of $T^*$-extensions (or Lagrangian extensions) of Lie superalgebras,  and show that they are classified by a certain cohomology space we introduce.  Several examples are provided to illustrate our construction.

\smallskip
\noindent \textbf{Keywords.} Quasi-Frobenius Lie superalgebra, orthosymplectic and periplectic forms; double extension; $T^*$-extension
\end{abstract}

\maketitle

\tableofcontents
\setcounter{tocdepth}{3}

\section {Introduction} 

\subsection{Quasi-Frobenius Lie algebras} Following, \cite{E,O} a Lie algebra $\fg$ is called a {\it quasi-Frobenius} Lie algebra if it admits a closed anti-symmetric non-degenerate bilinear from $\omega$. It is called a {\it Frobenius} Lie algebra if the bilinear $\omega$ is exact. Namely,  there exists $f\in \fg^*$ such that 
\begin{equation}\label{Kir}
\omega(x,y)=f([x,y]) \text{ for all $x,y\in \fg$}.
\end{equation}
In this case, the form $\omega$ in Eq. (\ref{Kir}) is sometimes referred to as a Kirillov form. 

According to \cite{O}, this terminology is due to Seligman. Those Lie algebras have been introduced to answer a  question raised by Jacobson: If $\fg$ is finite-dimensional,  what conditions one can put on $\fg$ in order the enveloping algebra $U(\fg)$ admits an exact simple module? The answer to this question was settled in \cite{O}.

On the other hand,  sometimes in the literature,  quasi-Frobenius Lie algebras were referred to as {\it symplectic } Lie algebras, see \cite{MR2,BC, BaBe}.   Perhaps this terminology stems from the fact that these Lie algebras are those of certain symplectic Lie groups over the field of the reals.

The superization of these notions is immediate,  see \S \ref{superFrob}. We have preferred to retain the term `quasi-Frobenius' in the super setting,  in order not to confuse them with the orthosymplectic $\mathfrak{osp}$ and the periplectic Lie superalgebras $\mathfrak{pe}$.

\subsection{Double extensions of Quasi-Frobenius Lie algebras} 
The notion of {\it symplectic double extension} of Lie algebras was introduced by Medina and Revoy \cite{MR2}. Roughly speaking,  one adds a symplectic plane to a quasi-Frobenius Lie algebra to obtain another quasi-Frobenius  Lie algebra.  This notion is the {\it anti-symmetric} analogue of the double extension of Lie algebras with symmetric non-degenerate bilinear forms introduced by Medina and Revoy  in \cite{MR1}  (such Lie algebras were also called {\it quadratic}). The notion of symplectic double extension was introduced to classify nilpotent symplectic groups by their Lie algebras.  Recall that a group admits an invariant  symplectic structure if there exists a left-invariant closed 2-form whose maximal rank is equal to the dimension of the group.  Medina and Revoy proved that the group of affine transformations admits an invariant and exact symplectic structure.  Further,  they proved that every nilpotent quasi-Frobenius Lie algebras can be obtained by a sequence of symplectic double extensions from the trivial Lie algebra.  As a by-product,  all $4$-dimensional quasi-Frobenius Lie algebras have been classified.

Regrettably,  a coboundary condition needed to carry out symplectic double extensions was missed in \cite{MR2}; this was rectified in \cite{MD}.  The converse, namely every quasi-Frobenius Lie algebra $\fg$ can be obtained as a symplectic double extension of a quasi-Frobenius Lie algebra,  is also true under certain conditions.  This double extension cannot be obtained when the center of $\fg$ is zero -- even in the case of solvable Lie algebras. Nevertheless,  in the case where $\fg$ is an abelian para-K\"ahler Lie algebra which is automatically solvable,  Benayadi and Bajo in \cite{BaBe} have generalized the notion of double extensions even if the center of $\fg$ is trivial.  Benayadi and Bajo also showed that such an algebra $\fg$ can be obtained by a sequence of generalized double extensions from the {\it affine} Lie algebra.  In fact, the case where $\fg$ is solvable but not necessarily abelian para-K\"ahler has not yet been solved.

\subsection{Double extensions of Quasi-Frobenius Lie superalgebras} In this paper,  we study the notion of {\it anti-symmetric} double extensions of Lie superalgebras,  superizing the results of Medina, Revoy and Dardi\'e  \cite{MR2,  MD}.  All our constructions are valid for Lie superalgebras defined over an arbitrary field of characteristic $p\not =2$.  This notion is the anti-symmetric analogue of the double extension of Lie superalgebras equipped with a non-degenerate symmetric bilinear form introduced and studied in \cite{ABB, ABBQ, BaBe, BeB, B}; see \cite{BeBou} for the peculiar case where the ground field is of characteristic $2$.  Our approach  is similar to that of Benayadi and Bajo \cite{BaBe} where the center might be trivial.  The converse,  namely every quasi-Frobenius Lie superalgebra can be obtained as a symplectic double extension of a quasi-Frobenius Lie superalgebra,  is also true under certain conditions,  see Theorems \ref{Rec1}, \ref{Rec2}, \ref{Rec3}, \ref{Rec4}.  To illustrate our construction,  a few examples are exhibited in Section \ref{Examples}.

We classify all 4-dimensional {\it real} quasi-Frobenius Lie superalgebras.  Our result is based on the classification of 4-dimensional real Lie superalgebras carried out by Blackhouse in \cite{Ba}.  Our list is given in Tables \ref{tab1}--\ref{tab6}, following the nomenclature of \cite{Ba}.  Surprisingly enough,  it turns out that some of the 4-dimensional Lie superalgebras can only admit  closed non-degenerate {\it non-homogeneous} bilinear forms.  This means that the purely even/odd part of the bilinear form turns into a degenerate form.  As the classification of 4-dimensional real Lie superalgebras was carried out in \cite{Ba},  we observe that every 4-dimensional real Lie superalgebra admitting a closed non-degenerate anti-symmetric form must be solvable,  superizing the result of \cite{BYC}.  

We also study the filiform Lie superalgebras $L^{n,m}$ and prove that they are quasi-Frobenius if and only if certain conditions on $n$ and $m$ are satisfied,  see Proposition \ref{filiform}.  
\subsection{$T^*$-extensions and $\Pi T^*$-extensions (also known as Lagrangian extensions)} 
As far as we know,  the notion of {\it $T^*$-extension} of Lie algebras was introduced by Bordemann \cite{Bor}.
Given a Lie algebra $\fh$ equipped with an $\fh^*$-valued 2-form $\alpha$.  A $T^*$-extension of $\fh$ is the Lie  algebra $\fg:=\fh\oplus \fh^*$ constructed by means of the form $\alpha$ that satisfies a certain cohomological requirement.   The Lie algebra $\fg$ naturally admits a non-degenerate bilinear form.   However,  in contrast to the double extension, the method
of $T^*$-extension applies not only to Lie algebras,  but to arbitrary nonassociative algebras and is a one-step rather than a multi-step extension,  see \cite{Bor}.

Bordemann showed that every finite-dimensional nilpotent Lie algebra $\fh$ with an invariant symmetric non-degenerate bilinear form   “is” a suitable $T^*$-extension.    If $\fh$ is a real finite-dimensional Lie algebra belonging to a Lie group $G$,  then the $T^*$-extension of $\fh$ is the Lie algebra of the cotangent bundle $T^*G$ of $G$,  considered as a Lie group.  This differential geometric fact justifies the terminology and the notation.

The notion of $T^*$-extension was also studied by Baues and Cort\'es in \cite{BC} in the context of Lie algebras admitting a~flat connection,  and called {\it Lagrangian extensions} since the Lie sub-algebra $\fh^*$ is a Lagrangian ideal in $\fh\oplus \fh^*$.  Baues and Cort\'es developed the notion of {\it Lagrangian cohomology} and showed that this cohomology captures all equivalence classes  of $T^*$-extensions. 

In this paper,  we {\it superize} the notion of $T^*$-extensions (or Lagrangian extensions),  generali-zing several results of \cite{BC}.  Following the same ideas of \cite{BC},  the Lie superalgebra $\fh$ is equipped with a {\it flat connection} and a {\it 2-form} that satisfies a certain cohomological requirement, see \S  \ref{TstarSec}.  It turns out that our construction remains valid over any field of characteristic $p\not =2$.  In the super setting,  we have observed that there are two ways to perform this extension: we consider either $\fg:=\fh\oplus \fh^*$ or $\fg:=\fh\oplus \Pi (\fh^*)$,  where $\Pi$ is the change of parity functor,   see \S \ref{basics}.  We call them $T^*$-extension and $\Pi T^*$-extensions of $\fh$,  respectively.  

Conversely, we prove that every {\it strongly polarized} Lie superalgebra  -- see \S \ref{TstarChap} for the definition -- is either a $T^*$-extension or $\Pi T^*$-extension of a flat Lie superalgebra,  see Theorem \ref{TstarConv}. 

We study equivalence classes of Lagrangian extensions, and show that isomorphic $T^*$-extensions (or $\Pi T^*$-extensions) give rise to the same associated quotient flat Lie superalgebra.

We show that isomorphism classes of Lagrangian extensions of a flat Lie superalgebra are classified by a suitable two-dimensional cohomology,  called the {\it Lagrangian extension cohomology},  superizing the construction of \cite{BC}.  We then show that the 2-cocycles give the same Lagrangian extensions if and and only if the 2-cocycles are cohomologous.  We further establish a map between Lagrangian cohomology and ordinary cohomology.

Several examples of four-dimensional real Lie superalgebras that can be obtained as $T^*$-extensions or $\Pi T^*$-extensions are given in \S \ref{Examples}.

\section {Background} 
Throughout this section,  $\Bbb K$ stands for an~arbitrary field of characteristic $p\not =2$. 
\subsection{Symmetric and anti-symmetric bilinear forms on superspaces}
Let $V$ and $W$ be two superspaces defined over an arbtray field ${\mathbb K}$ of characteristic $p\not =2$. 

Let $\cB\in \text{Bil}(V,W)$ be a bilinear form. The Gram matrix $B=(B_{ij})$ associated to $\cB$ is given by the formula, see \cite[Ch.1]{LSoS}, 
\be\label{martBil}
B_{ij}=(-1)^{p(B)p(v_i)}\cB(v_{i}, v_{j})\text{~~for the basis vectors $v_{i}\in V$}
\ee
in order to identify a~bilinear form $B(V, W)$ with an element of $\Hom(V, W^*)$. 
Consider the \textit{upsetting} of bilinear forms
$u\colon\Bil (V, W)\tto\Bil(W, V)$, see \cite[Ch.1]{LSoS}, given by the formula \be\label{susyB}
u(\cB)(w, v)=(-1)^{p(v)p(w)}\cB(v,w)\text{~~for any $v \in V$ and $w\in W$.}
\ee
In terms of the Gram matrix $B$ of $\cB$: the form
$\cB$ is  \textit{symmetric} if  and only if 
\be\label{BilSy}
u(B)=B,\;\text{ where $u(B)=
\mmat R^{t},(-1)^{p(B)}T^{t},(-1)^{p(B)}S^{t},-U^{t},$ for $B=\mmat R,S,T,U,$.}
\ee
Similarly, \textit{anti-symmetry} of  
$\cB$ means that $u(B)=-B$.

The phrase ``the bilinear form $\cB$ is \textit{supersymmetric} if $\cB(x,y)=(-1)^{p(x)p(y)}\cB(y,x)$'' can be often encountered in the literature,
although this property of the form $\cB$ is a~reflection of the braiding isomorphism $V\otimes W\simeq W\otimes V$ given by $(x,y)\mapsto(-1)^{p(x)p(y)}(y,x)$ and has nothing to do with (anti-)symmetry of the form.

Denote by $\Pi$ the functor that assigns to $V$ the superspace $\Pi(V)$ defined as the other copy of $V$ with the opposite parity of its homogeneous components. Namely, 
\[
(\Pi(V))_\ev = V_\od, \quad  (\Pi(V))_\od = V_\ev.
\]
The space $\Pi(V)$ consists of the linear combinations of elements that we denote by $\Pi(f)$ for every homogeneous $f\in V$.  For every superspace $V$,  denote the canonical odd homomorphism induced by the functor $\Pi$ by 
\[
\Pi: V \rightarrow \Pi(V) \quad f \mapsto \Pi(f).
\] 
We will identify $V$ and $\Pi(\Pi(V))$, naturally.

Observe that the passage from $V$ to $\Pi (V)$ turns every symmetric
form $\cB$ on $V$ into an anti-symmetric one $\cB^\Pi$ on $\Pi (V)$  and anti-symmetric $\cB$ into symmetric $\cB^\Pi$  by setting 
\be\label{ImF}
\text{$\cB^\Pi(\Pi(x), \Pi(y)):=(-1)^{p(\cB)+p(x)+p(x)p(y)}\cB(x,y)$ for any $x,y\in V$}.
\ee
{ \bf Following the custom in differential geometry, we will denote any anti-symmetric bilinear form by $\omega$.}
\sssbegin[(Symmetric and anti-symmetric forms for $p=2$)]{Remark}[Symmetric and anti-symmetric forms for $p=2$]
Over a field of charac-teristic $2$,  the definitions introduced in this section have to be amended,  see \cite{BB}.
\end{Remark}
\subsection{Lie superalgebras (Basics)}\label{basics}
A Lie superalgebra over a field of characteristic $p$, where $p\not =2,3$ is a $\mathbb Z/2$-graded vector space  $\fg=\fg_\ev\oplus \fg_\od$ endowed with a bilinear operation $[-,-]$ that satisfies  anti-commutativity and the Jacobi identities amended by the Sign Rule.  In characteristic $3$,  one has to add one more condition
\begin{equation}\label{p=3}
[f,[f,f]]=0\; \text{ for all $f \in \fg_\od$}.
\end{equation}
Superalgebras for which the bracket $[-,-]$ do not satisfy the condition (\ref{p=3}) are called {\it pre-Lie superalgebras},  see \cite{BBH}. 

A derivation of $\fg$ is a linear map on $\fg$ that satisfies:
\[
{\mathscr D}([f,g])=[ {\mathscr D}(f),g]  +(-1)^{p( {\mathscr D})p(f)} [f, {\mathscr D} (g)]) \text{ for all $f, g\in \fg$}.
\]
Denote by $\fder(\fg)$ the space of all derivations on $\fg$. 

Let $\omega \in \Bil(\fg, \fg)$ be an arbitrary non-degenerate bilinear form.  To each derivation ${\mathscr D}\in \fder(\fg)$ we can assign a unique linear map ${\mathscr D}^*$, called the {\it adjoint} of ${\mathscr D}$, that satisfies the following condition:
\[ 
\omega({\mathscr D}(f),g) = (-1)^{p(f) p({\mathscr D})} \omega(f, {\mathscr D}^*(g)).
\]
From the very definition of ${\mathscr D}^*$, it is easy to see that ${\mathscr D}^*$ is also a derivation.

\subsection{Quasi-Frobenius Lie superalgebras}\label{superFrob}
An even (resp. odd) non-degenerate bilinear form is called {\it orthosymplectic} (resp. {\it periplectic}).  

Let $\omega\in \text{Bil}(\fg, \fg)$ be an anti-symmetric bilinear form.  The form $\omega$ is called {\bf closed} if it satisfies the following condition (for all $f,g,h\in \fg$): 

\begin{equation}\label{coccond}
(-1)^{p(f)p(h)} \omega(f, [g , h]) +(-1)^{p(h)p(g)} \omega(h, [f , g] )+(-1)^{p(g)p(f)} \omega(g, [h , f] )=0.
\end{equation}


This is to say, $\omega$ is a non-degenerate 2-cocycle on $\fg$ with values in $\Kee$, see \cite{F}. Sometimes, it is more convenient to write the cocycle condition (\ref{coccond}) in the form 
\[
\omega(f, [g , h])  - (-1)^{p(f)p(g)} \omega(g, [f , h] ) -  \omega( [f , g] ,h)=0.
\]

A Lie superalgebra is called {\it quasi-Frobenius} if it is equipped with a~closed anti-symmetric non-degenerate bilinear form $\omega$. Such a Lie superalgebra will be denoted by $(\fg, \omega)$.

A quasi-Frobenius Lie superalgebra $(\fg, \omega)$ is called {\it Frobenius} if the form $\omega$ is exact.  Namely,  there exists $f\in \fg^*$ such that 
\[
\omega(x,y)=f([x,y])\text{ for all $x,y \in \fg$}.
\]
A quasi-Frobenius Lie superalgebra $(\fg, \omega)$ is called {\it orthosymplectic quasi-Frobenius} (resp.  {\it  periplectic quasi-Frobenius}) if the form $\omega$ is even (resp. odd). 

Let $(\fg, \omega)$ be a quasi-Frobenius Lie superalgebra and let $S \subseteq \fg$ be a subspace. The orthogonal of $S$ in $(\fg, \omega)$ is
\[
S^\perp = \{v \in \fg\; | \; \omega(v,w) = 0 \text{ for all } w \in S\}.
\]
The subspace $S$ is called {\it non-degenerate} if $S \cap S^\perp = \{0\}$. It is called
{\it isotropic} if ${S \subset  S^\perp}$. A maximal isotropic subspace is called {\it Lagrangian} if it satisfies $S = S^\perp$.

The following two lemmas were proved in \cite{MR2} in the non-super setting. Their superization is immediate.

\ssbegin[(Properties of orthogonal ideals)]{Lemma}[Properties of orthogonal ideals]\label{Iorth} Let $(\fg, \omega)$ quasi-Frobenius Lie superalgebra,  and let $I$ be an ideal of $\fg$.

\textup{(}i\textup{)} The space $I^\perp$ is a~Lie subsuperalgebra of $\fg$. 

\textup{(}ii\textup{)} The space $I^\perp$ is an ideal if and only if $I^\perp \subset  Z_\fg(I)$, where $Z_\fg(I)$ stands for the centralizer of $I$ in $\fg$. 

\end{Lemma}
\begin{proof}

This is an immediate consequence of the fact that $\omega$ is a  2-cocycle on the Lie superalgebra $\fg$. We also observe that if $I^\perp$ is an ideal and $\omega|_{I}$ is non degenerate, then $\fg$ is the orthogonal direct sum of the two ideals $I$ and $I^\perp$. 
\end{proof}
Consequently, any Lagrangian ideal must be abelian by Lemma \ref{Iorth}.

\ssbegin[(The center of $\fg$)]{Lemma}[The center of $\fg$]\label{center}
 Let $(\fg, \omega)$ be a Lie superalgebra equipped with a~closed anti-symmetric non-degenerate form $\omega$. Then, ${\mathfrak z}(\fg)\subseteq ([\fg,\fg])^\perp$.

\end{Lemma}
\begin{proof} 

Let $x\in \fz(\fg)$. For all $f,g\in \fg$, we have
\[
(-1)^{p(x) p(g)}\omega(x,[f,g])=-(1)^{p(g)p(f)}\omega(g, [x,f])-(-1)^{p(f)p(x)} \omega( f, [g,x])=0+0=0.\qed
\]
\noqed\end{proof}

An {\bf even} bilinear map $\nabla: \fg\times \fg \rightarrow \fg$, written as $(f,g)\mapsto \nabla_fg$, is called a connection on $\fg$. For all $f,g,h\in \fg$, the {\it torsion} $T^\nabla$ is defined by 
\[
T(f,g):=\nabla_fg-(-1)^{p(f)p(g)}\nabla_gf -[f,g]\quad  \text{ for all $f,g\in \fg$},
\] 
and the {\it curvature} $R^\nabla$ is defined by
\[
R^\nabla(f,g)h= \nabla_f \nabla_g h - (-1)^{p(f)p(g)} \nabla_g \nabla_f h -\nabla_{[f,g]}h\quad \text{ for all $f,g,h\in \fg$.}
\]
If $T=0$, then the connection $\nabla$ is called {\it torsion-free}. The curvature $R^\nabla$ vanishes if and only if the map $\sigma^\nabla:  \fg \rightarrow \End(\fg), \; f \mapsto \nabla_f$ is a representation of $\fg$ on itself. In this case, the connection $\nabla$ is called {\it flat}. 

A {\it flat Lie superalgebra} is a pair $(\fg, \nabla)$, where the connection $\nabla$ on $\fg$ is torsion-free and flat.

A subalgebra $\fh$ of $\fg$ is called a {\it totally geodesic subalgebra with respect to a connection} $\nabla$ if $\nabla_uv\in \fh$ for all $u,v\in \fh$.

Let $(\fg, \omega)$ \textup{(}resp. $(\fg, \kappa)$\textup{)} be an orthosymplectic \textup{(}resp.  periplectic\textup{)} quasi-Frobenius Lie superalgebra.  Let $I$ be an ideal of $g$ and let $I^\perp$ be its orthogonal with respect to $\omega$ \textup{(}resp. $\kappa$\textup{)} such that $[I,I^\perp]=0$. We have seen in Lemma \ref{Iorth} that the subalgebra $I^\perp$ is an ideal of $\fg$.  Let $\fh:= \fg \slash I^\perp$ be the associated quotient Lie superalgebra. If $\omega$ (resp. $\kappa$) is an~orthosymplectic (resp. periplectic) form on $\fg$, we obtain a non-degenerate bilinear pairing between $\fh$ and $I$ by declaring
\[
\omega_\fh( u, a):=\omega (\tilde u,a) \quad (\text{resp. }  \kappa_\fh(u, a):=\kappa (\tilde u,a)) \quad \text{for all $ u\in \fh$ and $a\in I$,}
\]
where $\tilde u$ is a lift of $u$ to $\fg$.

Let us show that this expression is well-defined. Indeed, suppose that $\bar u$ is another lift of $u$ to $\fg$. Let us write $\bar u- \tilde u=b\in I^\perp$. Using that fact that $I\subset I^\perp$, we get 
\begin{align*}
&\omega (\bar u ,a)=\omega (\tilde u+b,a)=\omega (\tilde u,a)+\omega (b,a)=\omega (\tilde u,a).
\end{align*}
The same argument can be used for the periplectic form $\kappa$.

\sssbegin[(Representation of $\fh:=\fg/I^\perp$)]{Proposition}[Representation of $\fh:=\fg/I^\perp$] Let $(\fg, \omega)$ \textup{(}resp. $(\fg, \kappa)$\textup{)} be an ortho-symplectic \textup{(}resp.  periplectic\textup{)} quasi-Frobenius Lie superalgebra.  Let $I$ be an ideal of $\fg$ and let $I^\perp$ be its orthogonal with respect to $\omega$ \textup{(}resp. $\kappa$\textup{)} such that $[I,I^\perp]=0$.  Let $\fh:=\fg/I^\perp$. The following statements hold:

\textup{(}i\textup{)} The map 
\[
\ad_{I}: \fh \rightarrow \End(I) \quad x \mapsto \ad_{\tilde x},
\] 
where $\tilde x$ is a lift of $x$ to $\fg$, is a representation. This representation induces a representation in the  spaces $I^*$ and $\Pi(I^*)$ as follows: \textup{(}for all $x\in \fh$ and for all $\xi \in I^*$\textup{)}
\begin{eqnarray}
\label{repadjn} \ad^*_I: \fh \rightarrow \End(I^*), &  \text{ where }  & \ad^*_I(x)(\xi):=- (-1)^{p(x)p(\xi)} \xi \circ \ad_{\tilde x}, \\[2mm]
\label{repadjpi}\Pi\ad^* _{I}: \fh \rightarrow \End(\Pi(I^*)), & \text{ where }  & \quad \Pi\ad^*_{I}(x)(\Pi(\xi)):=- (-1)^{p(x)(p(\xi)+1)} \Pi \circ \xi \circ \ad_{\tilde x}.
\end{eqnarray}

\textup{(}ii\textup{)} In the case of $(\fg, \omega)$,  the homomorphism 
\[
\fh \rightarrow I^*\quad u \mapsto \omega_\fh(u, \cdot)
\]
is an isomorphism. Moreover, it defines a one-cocycle in $Z^1_\ev(\fh, I^*)$, where the action of $\fh$ on $I^*$ is given by \eqref{repadjn}.

\textup{(}iii\textup{)} In the case of $(\fg, \kappa)$,  the homomorphism 
\[
\fh \rightarrow \Pi(I^*)\quad u \mapsto \Pi\circ  \kappa_\fh(u, \cdot)
\]
is an isomorphism. Moreover, it defines a one-cocycle in $Z^1_\ev(\fh, \Pi(I^*))$,  where the action of $\fh$ on $\Pi(I^*)$ is given by \eqref{repadjpi}.

\textup{(}iv\textup{)} The Lie superalgebra $\fh$ carries an even torsion-free flat connection $\nabla$ defined as follow:

\textup{(}a\textup{)} If the form $\omega$ is orthosymplectic on $\fg$:
\[
\omega_\fh( \nabla_u  v, a) =- (-1)^{p(\tilde v)p (\tilde u)} \omega (\tilde v, [\tilde u,a]), \quad \text{ for all $u,v \in \fh$ and for all $a\in I$}.
\]

\textup{(}b\textup{)} If the form $\kappa$ is periplectic on $\fg$:
\[
\kappa_\fh( \nabla_u  v, a) =- (-1)^{p(\tilde v)p (\tilde u)} \kappa (\tilde v, [\tilde u,a]), \quad \text{ for all $u,v \in \fh$ and for all $a\in I$}.
\]

\end{Proposition}
\begin{proof}

For Part (i), the representation $\ad_I$ is well-defined because $I$ is centralized by $I^\perp$.  Let us only show that the map (\ref{repadjpi}) is indeed a representation. First, we compute
\begin{align*}
 &\Pi\ad^*_{I}(u) \circ  \Pi\ad^*_{I}(v) (\Pi(\xi))= \Pi\ad^*_{I}(u) (- (-1)^{p(v)(p(\xi)+1)} \Pi\circ \xi \circ \ad_{\tilde v})\\[2mm]
&=(-1)^{(p(u)+p(v))(p(\xi)+1)+p(u)p(v)} \Pi \circ \xi \circ \ad_{\tilde v} \circ \ad_{\tilde u}.
\end{align*}
Now,
\begin{align*}
& \Pi\ad^*_{I}(u) \circ  \Pi\ad^*_{I}(v) (\Pi(\xi)) - (-1)^{p(u)p(v)} \Pi\ad^*_{I}(u) \circ  \Pi\ad^*_{I}(v) (\Pi(\xi)) -\Pi\ad^*_{I}([u,v]) (\Pi(\xi))\\[2mm]
& =(-1)^{(p(u)+p(v))(p(\xi)+1)} \left ( (-1)^{p(u)p(v)} \Pi \circ \xi \circ \ad_{\tilde v} \circ \ad_{\tilde u} - \Pi \circ \xi \circ \ad_{\tilde u} \circ \ad_{\tilde v}  + \Pi \circ  \xi \circ \ad_{[\tilde u, \tilde v]} \right ) =0.
\end{align*}
For Parts (ii) and (iii), let us first compute the kernel of $\omega_\fh$. It consists of all elements $u\in \fh$ such that $\omega(\tilde u, \cdot)=0$. This is equivalent to $\tilde u \in I^\perp$. Therefore,  $u=0$. Now,  since $\codim(I^\perp)=\dim(I)$, it follows that  $\dim(\fh)=\dim(I^*)$. Therefore, $\omega_\fh$ is an isomorphism. The same argument can be applied to the periplectic form $\kappa_\fh$. 

Let us check the 1-cocycle condition:
\begin{align*}
& u \cdot \Pi \circ \kappa_\fh(v, \cdot) - (-1)^{p(u)p(v)} v \cdot \Pi \circ \kappa_\fh(u, \cdot)- \Pi \circ \kappa_\fh([u,v], \cdot)  \\[2mm]
& =- (-1)^{p(u)(p( \kappa_\fh(v, \cdot))+1)} \Pi \circ \kappa_\fh(v, \cdot) \circ \ad_{\tilde u} + (-1)^{p(v)(p( \kappa_\fh(u, \cdot))+1)+p(u)p(v)} \Pi \circ \kappa_\fh(v, \cdot) \circ \ad_{\tilde u}\\[2mm]
& -\Pi \circ \kappa_\fh([\tilde u, \tilde v], \cdot) = -(-1)^{p(u)p(v)}  \Pi \circ (\kappa_\fg( \tilde v, \ad_{\tilde u} (\cdot))+  \kappa_\fg( \tilde u, \ad_{\tilde v} (\cdot) ) -\kappa_\fg([\tilde u,\tilde v], \cdot))=0,
\end{align*}
since $\kappa_\fg$ is periplectic on $\fg$.

Part (iv) follows from the fact that $\omega_\fh$ and $\kappa_\fh$ are isomorphisms together with the fact that $\omega$ and $\kappa$ are closed forms.
\end{proof}

\section{A few examples of Frobenius Lie superalgebras}
\subsection{Classification of $4$-dimensional real Frobenius Lie superalgebras}
Throughout this section the ground field $\mathbb K$ is the field of real numbers $\mathbb R$. 

This section is motivated by the following result due to Bon Yao Cho, see \cite{BYC}.

\ssbegin[(Solvability of $\fg$,  see \cite{BYC})]{Theorem}[Solvability of $\fg$, see \cite{BYC}]\label{4-nil}
Every $4$-dimensional quasi-Frobenius Lie algebra is solvable. 
\end{Theorem}

Over $\mathbb C$, it  is well-known that a Lie superalgebra $\fg$ is solvable if and only if the even part $\fg_\ev$ is solvable,  see \cite{K}.  This is not true anymore if the ground field if of characteristic 2,  see examples in \cite{BGL, BGLLS}.

We will classify all $4$-dimensional quasi-Frobenius Lie superalgebras,  and check whether an analog of Theorem \ref{4-nil}  is valid or not in the supper setting. The classification of $4$-dimensional Lie superalgebras has been carried out in \cite{Ba} (see also for \cite{H} for the nilpotent case). The following lemma will reduce our study to the case where $\sdim(\fg)=2|2$ and $\fg_\od\not =\{0\}$. 

\ssbegin[(Superdimension constraints)]{Lemma}[Superdimension constraints] \label{NH} Let $V=V_\ev\oplus V_\od$ be a finite-dimensional superspace such that $V_\od \not = \{0\}$, equipped with an~anti-symmetric non-degenerate bilinear form $\omega$.

\textup{(}i\textup{)}  If $p(\omega)=\ev$, then $\dim(V_\ev)$ is even.

\textup{(}ii\textup{)} If $p(\omega)=\od$, then $\dim(V_\ev)=\dim(V_\od)$.

\end{Lemma}
\begin{proof}

(i). The result follows from the fact that $\omega|_{\fg_\ev}$ is non-degenerate, and from the well-known result in the non-super setting stating that such forms exist only on even dimensional spaces. 

(ii). Let $0\not =x\in V_\od$. Since $\omega$ is non-degenerate and odd, there exists $y \in V_\ev$ such that $\omega(x,y)=1$. Consider now the vector superspace $S=\Span\{x,y\}$. We claim that $S^\perp \cap S=\{0\}$. Indeed, let $z\in S^\perp \cap S$ and let us write $z=\lambda_1 x+\lambda_2 y$ for some $\lambda_1, \lambda_2 \in \mathbb{K}$. It follows that
\[
0=\omega(x,z)=\omega(x, \lambda_1 x+\lambda_2 y)=\omega(x,\lambda_2 y)=\lambda_2.
\]
Similarly, $\lambda_1=0$. Let us now consider an arbitrary $w\in V$. We can write 
\[
u= (u+\omega(u,x)y-\omega(u,y)x)-(\omega(u,x)y-\omega(u,y))x \in S^\perp + S.
\]
It follows that $V=S\oplus S^\perp$, and therefore $\sdim(V)= \sdim(S^\perp)+(1|1)$. The form $\omega|_{S^\perp}$ is also non-degenerate and the same argument can be applied to $S^\perp$. The result follows because $V$ is finite-dimensional. 
\end{proof}

\ssbegin[(Existence of non-homogenous forms)]{Remark}[Existence of non-homogenous forms]
Lemma \ref{NH} considers only homo-geneous forms $\omega$. For instance, the Lie superalgebras $\fg:=D^1, D^2_{-1}, D^{11}_{pq}, D^{13}_{-1}, D^{15}, D^3_{-\frac{1}{2}, \frac{1}{2}}$ in Tables \ref{tab3},  \ref{tab4} and \ref{tab5} admit closed anti-symmetric non-degenerate bilinear forms for which $\mathrm{dim}(\fg_\ev)$ is not even and $\mathrm{dim}(\fg_\ev)\not = \mathrm{dim}(\fg_\od)$.
\end{Remark}

We have regrouped the list of Lie superalgebras admitting closed non-degenerate bilinear forms following the nomenclature given in \cite{Ba}. Tables \ref{tab1}, \ref{tab3} and \ref{tab4} contain Lie superalgebras for which $[\fg_\od, \fg_\od]=\{0\}$ -- called {\it trivial} algebras -- while Table \ref{tab2}, \ref{tab5} and \ref{tab6} contain Lie superalgebras for which $[\fg_\od, \fg_\od]\not =\{0\}$ -- called {\it non-trivial} algebras.  As these structures are not necessarily unique,  we have chosen to select,  when possible,  the bilinear form which is homogeneous (either purely even or purely odd) and satisfies the 2-cocycle condition \eqref{coccond}. The label `NH' stands for non-homogeneous non-degenerate antisymmetric bilinear forms captured by several examples in Table \ref{tab2}, \ref{tab3}, \ref{tab4} and \ref{tab5}.  The Lie superalgebra $\fg$ is spanned by $e_1, e_2,e_3$ and $e_4$.  The list of even generators and odd generators is separated by the vertical bar $|$.  

The following proposition gives a positive answer to the super analog of Theorem \ref{4-nil}. 

\ssbegin[(Solvability of $\fg$ -- the super case)]{Theorem}[Solvability of $\fg$ -- the  super case]
Every real $4$-dimensional quasi-Frobenius Lie superalgebra is solvable.

\end{Theorem}
\begin{proof}

The proof follows from the classification of real $4$-dimensional Lie superalgebras carried out in \cite{Ba},  and Tables \ref{tab1}-- \ref{tab6}.
\end{proof}

\small
\begin{table}[H]
\centering
\begin{tabular}{| c | l | m{7cm} | c| }\hline
 The LSA & Relations in the basis: $e_1, e_2 \, |\, e_3, e_4$ & The form $\omega$  & $p(\omega)$ \\ \hline
 $D^5$ & 
$
\begin{array}{lcllcl}
[e_1,e_3]&=&e_3,\\[1mm] 
[e_1, e_4]&=&e_4, \\[1mm]
[e_2, e_4]&=&e_3
\end{array}
$
& $\lambda e_1^*\wedge e_4^*+\mu (e_1^*\wedge e_3^*+ e_2^*\wedge e_4^*)$, where $\mu\not =0$ & $\od$  \\[2mm]  \hline
$D^6$ &  
$
\begin{array}{lcrlcl}
[e_1, e_3]&=&e_3, \\[1mm]
 [e_1, e_4]&=&e_4,\\[1mm]
[e_2, e_3]&=&-e_4, \\[1mm]
 [e_2, e_4]&=&e_3
\end{array}
$&$\lambda (e_1^*\wedge e_3^* + e_2^*\wedge e_4^*)- \mu (e_1^*\wedge e_4^*-e_2^*\wedge e_3^*) $, where $\mu^2+\lambda^2\not=0$ & $\od$  \\[2mm] \hline
$\begin{array}{c}
D^7_{pq},\\[1mm]
 pq\not =0,\\[1mm]
p\geq q
\end{array}
$ 
& 
$
\begin{array}{lcllcl}
[e_1,e_2]&=&e_2, \\[1mm]
 [e_1, e_3]&=&p e_3,\\[1mm]
[e_1, e_4]&=&q e_4
\end{array}$ & \text{None} &--  \\[2mm] \hline
$
\begin{array}{c}
D^7_{-1q},\\[1mm]
q\leq -1 
\end{array} 
$
&
$
\begin{array}{lcllcl}
[e_1,e_2]&=&e_2,\\[1mm]
[e_1, e_3]&=&- e_3,\\[1mm]
 [e_1, e_4]&=&q e_4
\end{array}$ & $\lambda e_1^*\wedge e_3^*+\mu e_1^*\wedge e_4^*+\nu e_2^*\wedge e_3^* $, where $\mu \nu\not =0$ & $\od$   \\[2mm] \hline
$
\begin{array}{c}
D^7_{pp},\\[1mm]
p=-1
\end{array} 
$
&
$
\begin{array}{lcllcl}
[e_1,e_2]&=&e_2,\\[1mm]
[e_1, e_3]&=&- e_3,\\[1mm]
 [e_1, e_4]&=&-e_4
\end{array}$ & $\begin{array}{l}
\lambda e_1^*\wedge e_3^*+\mu e_1^*\wedge e_4^*+\nu e_2^*\wedge e_3^* +\gamma e_2^*\wedge e_4^* \\
\text{where } \mu\nu-  \lambda \gamma\not =0
\end{array}$ & $\od$   \\[2mm] \hline
$
\begin{array}{c}
D^7_{pq},\\[1mm]
p=-q , \\[1mm]
q\not =0,-1
\end{array} 
$
&
$
\begin{array}{lcllcl}
[e_1,e_2]&=&e_2,\\[1mm]
 [e_1, e_3]&=&p e_3,\\[1mm]
 [e_1, e_4]&=&-p e_4
\end{array}$ & $ \lambda e_1^*\wedge e_2^*+\mu e_3^*\wedge e_4^* $, where $\lambda \mu\not =0$ & $\ev$   \\[2mm] \hline
$
\begin{array}{c}
D^7_{pq},\\[1mm]
p=-q , \\[1mm]
q =-1
\end{array} 
$
&
$
\begin{array}{lcllcl}
[e_1,e_2]&=&e_2,\\[1mm]
 [e_1, e_3]&=& e_3,\\[1mm]
 [e_1, e_4]&=&- e_4
\end{array}$ & $
 \lambda e_1^*\wedge e_2^*+\mu e_3^*\wedge e_4^* , \text{ where } \lambda \mu\not =0$ \newline
$\lambda e_1^*\wedge e_3^*+ \mu e_1^*\wedge e_4^*+\nu e_2^*\wedge e_4^*,
\text{where } \lambda \nu \not =0$ & $\begin{array}{l}
\ev\\[1mm]
\od
\end{array}$   \\[2mm] \hline
$
\begin{array}{c}
D^8_{p},\\
p\not=0 
\end{array} 
$
&
$
\begin{array}{lcllcl}
[e_1, e_2]&=&e_2,\\[1mm]
[e_1, e_3]&=&p e_3,\\[1mm]
 [e_1, e_4]&=&e_3+ pe_4
\end{array}$ & None & --
\\[2mm] \hline
$D^8_{-1}$ & 
$
\begin{array}{lcllcl}
[e_1, e_2]&=&e_2,\\[1mm]
 [e_1, e_3]&=&- e_3, \\[1mm]
[e_1, e_4]&=&e_3-e_4
\end{array}$ & $ \lambda e_1^*\wedge e_3^*+ \nu e_1^*\wedge e_4^*+ \mu e_2^*\wedge e_4^*,$ where $\lambda \mu\not =0$ & $\od$ \\[2mm] \hline
$
\begin{array}{c}
D^9_{pq}\\[1mm]
q>0
\end{array} 
$

& 
$
\begin{array}{lcllcl}
 [e_1, e_3]&=&pe_3-qe_4, \\[1mm]
 [e_1, e_2]&=&e_2, \\[2mm]
 [e_1, e_4]&=&q e_3+pe_4
\end{array}$ & None  &--\\[2mm] \hline
  $D^{10}_{q}$ & $
\begin{array}{lcllcl}
[e_1, e_2]&=&e_2, \\[1mm]
 [e_1, e_3]&=&(q+1)e_3,\\[1mm]
 [e_1, e_4]&=&q e_4, \\[1mm]
[e_2, e_4]&=&e_3
\end{array}$ & $ (1+q) \lambda e_1^*\wedge e_3^*+ \mu e_1^*\wedge e_4^*+\lambda e_2^*\wedge e_4^*, $ where $(1+q)\lambda \not=0$ & $\od$ \\[2mm] \hline
  $D^{10}_{-2}$ &
$
\begin{array}{lcrlcl}
[e_1, e_2]&=&e_2,\\[1mm]
[e_1, e_3]&=&-e_3,\\[1mm]
 [e_1, e_4]&=&-2 e_4, \\[1mm]
 [e_2, e_4]&=&e_3
\end{array}$ & $ \lambda (e_2^*\wedge e_4^*- e_1^*\wedge e_3^*)+\mu e_1^*\wedge e_4^*+\nu e_2^*\wedge e_3^*$, where $\mu \nu+\lambda ^2\not =0$ & $\od$ \\[2mm] \hline
\end{tabular}
\caption{Trivial algebras (i.e., $[\fg_\od, \fg_\od]=\{0\}$) with $\mathrm{sdim}=2|2$} \label{tab1}
\end{table}
\normalsize
\small
\begin{table}[H]
\centering
\begin{tabular}{| c | l | c |c|}\hline
 The LSA & Relations in the basis: $e_1, e_2 \, |\, e_3, e_4$ & The form $\omega$ & $p(\omega)$ \\ \hline
 $(D^7_{1/2\;1/2})^1$ & $
\begin{array}{lcrlcl}
[e_1, e_2]&=&e_2, & [e_1, e_3]&=&\frac{1}{2} e_3\\[1mm]
 [e_1, e_4]&=&\frac{1}{2}e_4, & [e_3,e_3]&=&e_2  \\[1mm]
[e_4, e_4]&=&e_2
\end{array}$ & $e_1^*\wedge e_2^* - \frac{1}{2} e_3^*\wedge e_3^*- \frac{1}{2} e_4^*\wedge e_4^*$ & $ \ev$\\ \hline
 $(D^7_{1/2\;1/2})^2$ & $
\begin{array}{lcrlcl}
[e_1, e_2]&=&e_2, & [e_1, e_3]&=&\frac{1}{2} e_3\\[1mm]
 [e_1, e_4]&=&\frac{1}{2}e_4, & [e_3,e_3]&=&e_2  \\[1mm]
[e_4, e_4]&=& - e_2
\end{array}$ &  $-e_1^*\wedge e_2^* + \frac{1}{2} e_3^*\wedge e_3^* -  \frac{1}{2} e_4^*\wedge e_4^*$ & $\ev$\\ \hline 
 $(D^7_{1/2\;1/2})^3$ & $
\begin{array}{lcrlcl}
[e_1, e_2]&=&e_2, & [e_1, e_3]&=&\frac{1}{2} e_3\\[1mm]
 [e_1, e_4]&=&\frac{1}{2}e_4, & [e_3,e_3]&=&e_2 
\end{array}$ & None &$-$\\ \hline 
 $\begin{array}{c}
D^7_{1-p,p}\\[1mm]
p\leq \frac{1}{2}
\end{array}$ & $
\begin{array}{lcllcl}
[e_1, e_2]&=&e_2, &
 [e_1, e_3]&=& p e_3\\[1mm]
 [e_1, e_4]&=&(1-p) e_4, \\[1mm]
 [e_3,e_4]&=&e_2 
\end{array}$ & $e_1^*\wedge e_2^* -  e_3^*\wedge e_4^*$ & $\ev$\\ \hline 
 $D^8_{1/2}$ & $
\begin{array}{lcllcl}
[e_1, e_2]&=&e_2, &[e_1, e_3]&=& \frac{1}{2} e_3\\[1mm]
 [e_1, e_4]&=&e_3+\frac{1}{2} e_4, \\[1mm]
 [e_4,e_4]&=&e_2 
\end{array}$ & None & $-$\\ \hline 
 $
\begin{array}{c}
D^9_{1/2,p}\\[1mm]
p>0
\end{array}$ & $
\begin{array}{lcllcl}
[e_1, e_2]&=&e_2,  & [e_3,e_3]&=&e_2 \\[1mm]
 [e_1, e_4]&=&p e_3+\frac{1}{2} e_4, & [e_4,e_4]&=&e_2 \\[1mm]
[e_1, e_3]&=& \frac{1}{2} e_3 -pe_4
\end{array}$ & $e_1^*\wedge e_2^* - \frac{1}{2}e_3^* \wedge e_3^* - \frac{1}{2}e_4^*\wedge e_4^*$ & $\ev$\\ \hline
 $(D^{10}_{0})^1$ & $
\begin{array}{lcllcl}
[e_1, e_2]&=&e_2,  & [e_1,e_3]&=&e_3 \\[1mm]
 [e_2, e_4]&=& e_3, & [e_4,e_4]&=&e_1 \\[1mm]
[e_3, e_4]&=&- \frac{1}{2} e_2
\end{array}$ & $\begin{array}{l}
-2 \lambda e_1^*\wedge e_2^*+\mu e_1^*\wedge e_3^* +\mu e_1^*\wedge e_4^*\\[1mm]
- \mu e_3^*\wedge e_4^* - \frac{\nu}{2} e_4^* \wedge e_4^* \\[1mm]
 \text{where } \mu^4-4 \lambda^4\not =0
\end{array}$ & NH \\ \hline
 $(D^{10}_{0})^2$ & $
\begin{array}{lcllcl}
[e_1, e_2]&=&e_2,  & [e_1,e_3]&=&e_3 \\[1mm]
 [e_2, e_4]&=& e_3, & [e_4,e_4]&=&-e_1 \\[1mm]
[e_3, e_4]&=& \frac{1}{2} e_2
\end{array}$ & $\begin{array}{l}
2 \lambda e_1^*\wedge e_2^*+\mu e_1^*\wedge e_3^* +\mu e_1^*\wedge e_4^*\\[1mm]
- \mu e_3^*\wedge e_4^* - \frac{\nu}{2} e_4^* \wedge e_4^* \\[1mm]
 \text{where } \mu^4-4 \lambda^4\not =0
\end{array}$ & NH\\\hline
 $(2A_{1,1}+2A)^1$ & $
\begin{array}{lcllcl}
[e_3, e_3]&=&e_1,  & [e_4,e_4]&=&e_2 
\end{array}$ & None & $-$\\ \hline
 $(2A_{1,1}+2A)^2$ &   $
\begin{array}{lcllcl}
[e_3, e_3]&=&e_1,  & [e_4,e_4]&=&e_2 \\[1mm]
[e_3, e_4]&=&e_1\\
\end{array}$ & None & $-$\\ \hline
 $(2A_{1,1}+2A)^3_p$ & $
\begin{array}{lcllcl}
[e_3, e_3]&=&e_1,  & [e_4,e_4]&=&e_2 \\[1mm]
[e_3, e_4]&=&p(e_1+e_2)\\
\end{array}$  & $\begin{array}{l}
\text{For $p\not =\frac{1}{2}$: None}\\[1mm]
\text{For $p=\frac{1}{2}:\, e_2^*\wedge e_3^* - e_1^*\wedge e_4^*$}
\end{array}$ &
$\begin{array}{l}  
-\\
\od
\end{array}
$\\ \hline
 $(2A_{1,1}+2A)^4_p$ & $
\begin{array}{lcllcl}
[e_3, e_3]&=&e_1,  & [e_4,e_4]&=&e_2 \\[1mm]
[e_3, e_4]&=&p(e_1-e_2)\\
\end{array}$ & None &
$-$\\ \hline
 $C^1_1+A$ & $
\begin{array}{lcllcl}
[e_1, e_2]&=&e_2,  & [e_1,e_3]&=&e_3 \\[1mm]
[e_3, e_4]&=&e_2\\
\end{array}$ & $e_1^*\wedge e_2^* - e_3^*\wedge e_4^* - \frac{1}{2} e_4^*\wedge e_4^*$ & $\ev$\\ \hline
 $C^1_{1/2}+A$ & $
\begin{array}{lcllcl}
[e_1, e_2]&=&e_2,  & [e_1,e_3]&=&\frac{1}{2}e_3 \\[1mm]
[e_3, e_3]&=&e_2\\
\end{array}$ & $\begin{array}{l}
\lambda e_1^*\wedge e_2^* - \frac{\lambda}{2} e_3^*\wedge e_3^* - \frac{\mu}{2}e_4^*\wedge e_4^*, \\
\text{where }\lambda\mu\not =0
\end{array}$ & $\ev$\\ \hline
 $C^2_{-1}+A$ & $
\begin{array}{lcllcl}
[e_1, e_3]&=&e_3,  & [e_1,e_4]&=& - e_4 \\[1mm]
[e_3, e_4]&=&e_2\\
\end{array}$ & None & $-$\\ \hline
 $C^3+A$ & $
\begin{array}{lcllcl}
[e_1, e_4]&=&e_3,  & [e_4,e_4]&=& e_2 
\end{array}$ & $2 e_1^*\wedge e_2^* - e_3^*\wedge e_4^*$ & $\ev$\\ \hline
 $C^5_0+A$ & $
\begin{array}{lcllcl}
[e_1, e_3]&=& - e_4,  & [e_1,e_4]&=& e_3 \\[1mm]
[e_3, e_3]&=&e_2 & [e_4, e_4]&=&e_2\\
\end{array}$ & None & $-$\\
\hline
\end{tabular}
\caption{Non-trivial algebras (i.e., $[\fg_\od, \fg_\od]\not=\{0\}$) with $\mathrm{sdim}=2|2$}\label{tab2}
\end{table}
\small
\begin{table}[H]
\centering
\begin{tabular}{| c | l | c |c|}\hline
 The LSA & Relations in the basis: $e_1, e_2, e_3 \, |\, e_4$ & The form $\omega$ & $p(\omega)$ \\ \hline
 $D^1$ & $
\begin{array}{lcrlcl}
[e_2, e_3]&=&e_1, & [e_2, e_4]&=&e_4
\end{array}$ & $\begin{array}{l}
\lambda e_1^*\wedge e_3^*+ \mu e_2\wedge e_4^*\\[1mm]
+\nu e_1^*\wedge e_2^*+ \gamma e_2\wedge e_3^*\\[1mm]
\text{where } \lambda \mu \not =0
\end{array}$ & NH\\ \hline
 $\begin{array}{l}
D^2_q\\[1mm]
q \not =-1,0
\end{array}$ & $
\begin{array}{lcrlcl}
[e_1, e_2]&=& e_1, & [e_2, e_3]&=& e_1+ e_2,\\[1mm]
 [e_3, e_4]&=& q e_4
\end{array}$ & None & $-$\\ \hline
 $\begin{array}{l}
D^2_{-1}
\end{array}$ & $
\begin{array}{lcrlcl}
[e_1, e_2]&=& e_1, & [e_2, e_3]&=& e_1+ e_2,\\[1mm]
 [e_3, e_4]&=& - e_4
\end{array}$ & $
\begin{array}{l}
\lambda e_1^*\wedge e_3^*+ \mu e_2^*\wedge e_3^*+\\[2mm]
 \nu e_2^*\wedge e_4^*+ \gamma e_3^*\wedge e_4^*\\[2mm]
\text{where } \lambda \nu \not =0
\end{array}$ & NH\\ \hline
 $\begin{array}{l}
D^3_{pq}\\[1mm]
pq\not =0
\end{array}$ & $
\begin{array}{lcrlcl}
[e_1, e_3]&=& pe_1-e_2, & [e_2, e_3]&=& e_1+ pe_2,\\[1mm]
 [e_3, e_4]&=&q e_4
\end{array}$ & None & $-$\\ \hline
\end{tabular}
\caption{Trivial algebras (i.e., $[\fg_\od, \fg_\od]=\{0\}$) with $\mathrm{sdim}=3|1$}\label{tab3}
\end{table}
\begin{table}[H]
\centering
\begin{tabular}{| c | l | c |c|}\hline
 The LSA & Relations  in the basis: $e_1\,|\,e_2,e_3,e_4$ & The form $\omega$ & $p(\omega)$ \\ \hline
 $\begin{array}{c}
D^{11}_{pq}\\
0< |p|\leq |q|\leq 1
\end{array}$ & $
\begin{array}{lcrlcl}
[e_1, e_2]&=&e_2, & [e_1, e_3]&=&p e_3,\\[1mm]
 [e_1, e_4]&=&q e_4
\end{array}$ & None& $-$\\ \hline
$\begin{array}{c}
D^{11}_{pq}\\
(p,q)=(-1,-1)
\end{array}$ & $
\begin{array}{lcrlcl}
[e_1, e_2]&=&e_2, & [e_1, e_3]&=&p e_3,\\[1mm]
 [e_1, e_4]&=&q e_4
\end{array}$ & $
\begin{array}{l}\lambda e_1^*\wedge e_2^*+ \mu e_1^*\wedge e_3^*+\nu e_1^*\wedge e_4^*\\[1mm]
+\gamma e_2^*\wedge e_3^*+\delta e_2^*\wedge e_4^*,\\[1mm]
\text{where } \nu \gamma-\mu \delta\not =0
\end{array}$ & NH\\ \hline
$\begin{array}{c}
D^{11}_{pq}\\
0< |p|\leq 1\\
p=-q
\end{array}$ & $
\begin{array}{lcrlcl}
[e_1, e_2]&=&e_2, & [e_1, e_3]&=&p e_3,\\[1mm]
 [e_1, e_4]&=&q e_4
\end{array}$ & $
\begin{array}{l}\lambda e_1^*\wedge e_2^*+ \mu e_1^*\wedge e_3^*+\nu e_1^*\wedge e_4^*\\[1mm]
+\gamma e_3^*\wedge e_4^*,\text{ where } \lambda \gamma\not =0
\end{array}$ & NH\\ \hline
$\begin{array}{c}
D^{11}_{pq}\\
0< |p|< 1\\
q=-1
\end{array}$ & $
\begin{array}{lcrlcl}
[e_1, e_2]&=&e_2, & [e_1, e_3]&=&p e_3,\\[1mm]
 [e_1, e_4]&=&q e_4
\end{array}$ & $
\begin{array}{l}\lambda e_1^*\wedge e_2^*+ \mu e_1^*\wedge e_3^*+\nu e_1^*\wedge e_4^*\\[1mm]
+\gamma e_2^*\wedge e_4^*, \text{ where } \mu \gamma\not =0
\end{array}$
 & NH\\ \hline
 $D^{12}$ & $
\begin{array}{lcrlcl}
[e_1, e_2]&=&e_2, & [e_1, e_4]&=&e_3,
\end{array}$ & None & $-$\\ \hline
 $\begin{array}{c}
D^{13}_p\\[1mm]
p\not =0 \text{ (generic)}
\end{array}$ & $
\begin{array}{lcllcl}
[e_1, e_2]&=&p e_2, & [e_1, e_3]&=&e_3,\\[1mm]
 [e_1, e_4]&=&e_3+e_4
\end{array}$ & None & $-$\\ \hline
$\begin{array}{c}
D^{13}_{-1}
\end{array}$ & $
\begin{array}{lcllcl}
[e_1, e_2]&=& - e_2, & [e_1, e_3]&=&e_3,\\[1mm]
 [e_1, e_4]&=&e_3+e_4
\end{array}$ & $
\begin{array}{l}
\lambda e_1^*\wedge e_2^*+ \mu e_1^*\wedge e_3^*+ \gamma e_1^*\wedge e_4^*\\[1mm]
+\nu e_2^*\wedge e_4^*, \text{ where } \mu \nu\not =0
\end{array}$ & NH\\ \hline
$\begin{array}{c}
D^{14}_{pq}\\[1mm]
p\not =0, q\geq 0
\end{array}$ & $
\begin{array}{lcllcl}
[e_1, e_2]&=&p e_2, \\[1mm]
 [e_1, e_3]&=&q e_3-e_4,\\[1mm]
 [e_1, e_4]&=&q e_4+e_3
\end{array}$ & None & $-$\\ \hline
$D^{15}$ & $
\begin{array}{lcrlcl}
[e_1, e_3]&=&e_2 ,& [e_1, e_4]&=&e_3,
\end{array}$ & 
$\begin{array}{l}
\lambda e_1^*\wedge e_2^*+\mu e_1^*\wedge e_3^*+ \nu e_1^*\wedge e_1^* \\[1mm]
+ \delta e_2^*\wedge e_4^* -\frac{1}{2} \delta e_3^*\wedge e_3^*,\\[1mm]
\text{where } \delta^2(\mu^2-2 \lambda  \nu)\not =0.
\end{array}$ & NH\\ \hline
$D^{16}$ & $
\begin{array}{lcl}
[e_1, e_2]&=&e_2 ,\\[1mm]
 [e_1, e_3]&=&e_2+e_3,\\[1mm]
 [e_1, e_4]&=&e_3+e_4
\end{array}$ & None & $-$\\ \hline
\end{tabular}
\caption{Trivial algebras (i.e., $[\fg_\od, \fg_\od]=\{0\}$) with $\mathrm{sdim}=1|3$}\label{tab4}
\end{table}
\small
\begin{table}[H]
\centering
\begin{tabular}{| c | l | c |c|}\hline
 The LSA & Relations in the basis: $e_1, e_2, e_3 \, |\, e_4$& The form $\omega$ & $p(\omega)$ \\ \hline
 $A_{3,1}+A$ & $
\begin{array}{lcrlcl}
[e_2, e_3]&=&e_1, & [e_4, e_4]&=&e_1
\end{array}$ & None & $-$\\ \hline
$\begin{array}{l}
D^3_{p, -1/2}\\[1mm]
p\not =0
\end{array}$ & $
\begin{array}{lcrlcl}
[e_1, e_2]&=&e_2, & [e_1, e_3]&=&p e_3 , \\[1mm]
[e_1, e_4]&=&\frac{1}{2}e_4, & [e_4, e_4]&=& e_2\\[1mm]
\end{array}$ & None & $-$\\ \hline
$D^3_{-1/2, -1/2}$ & $
\begin{array}{lcrlcl}
[e_1, e_2]&=&e_2, & [e_1, e_3]&=&- \frac{1}{2} e_3 , \\[1mm]
[e_1, e_4]&=&\frac{1}{2}e_4, & [e_4, e_4]&=& e_2\\[1mm]
\end{array}$ & $\begin{array}{l}
\frac{1}{2} \lambda e_1^*\wedge e_2^*+ \mu e_1^*\wedge e_3^*\\[1mm]
+ \nu e_1^*\wedge e_4^*+\gamma e_3^*\wedge e_4^*\\[1mm]
+\frac{1}{2}\lambda e_4^*\wedge e_4^*, \text{ where }\gamma \lambda\not =0
\end{array}$ & NH\\ \hline
 $(D^2_{-1/2})^1$ & $
\begin{array}{lcrlcl}
[e_1, e_2]&=&e_2, & [e_1, e_3]&=&e_2+e_3 , \\[1mm]
[e_1, e_4]&=&\frac{1}{2}e_4, & [e_4, e_4]&=&e_2
\end{array}$ & None & $-$\\ \hline
 $(D^2_{-1/2})^2$ & $
\begin{array}{lcrlcl}
[e_1, e_2]&=&e_2, & [e_1, e_3]&=&-e_2+e_3 , \\[1mm]
[e_1, e_4]&=&\frac{1}{2}e_4, & [e_4, e_4]&=&e_2
\end{array}$ & None & $-$\\ \hline
\end{tabular}
\caption{Trivial algebras (i.e. , $[\fg_\od, \fg_\od]\not =\{0\}$) with $\mathrm{sdim}=3|1$}\label{tab5}
\end{table}
\small
\small
\begin{table}[H]
\centering
\begin{tabular}{| c | l | c |c|}\hline
 The LSA & Relations in the basis: $e_1\,|\,e_2,e_3,e_4$ & The form $\omega$ & $p(\omega)$ \\ \hline
 $(A_{1,1}+3 A)^1$ & $
\begin{array}{lcrlcl}
[e_2, e_2]&=&e_1, & [e_3, e_3]&=&e_1,\\[1mm]
 [e_4, e_4]&=&e_1
\end{array}$ & None & $-$\\ \hline
$(A_{1,1}+3 A)^2$ & $
\begin{array}{lcrlcl}
[e_2, e_2]&=&e_1, & [e_3, e_3]&=&e_1,\\[1mm]
 [e_4, e_4]&=&-e_1
\end{array}$ & None & $-$\\ \hline
\end{tabular}
\caption{Trivial algebras (i.e., $[\fg_\od, \fg_\od]\not =\{0\}$) with $\mathrm{sdim}=1|3$}\label{tab6}
\end{table}

\newpage
\subsection{Filiform Lie superalgebras}
Following \cite{V}, a Lie algebra $\fg$ is called filiform if $\fg$ is nilpotent whose nil index is maximal, that is equal to $\dim(\fg)-1$. The  simplest filiform Lie algebra, $L^{n}$, is defined by the only non-zero
bracket products:
\[
\begin{array}{lcl}
[X_1,X_i] & = & X_{i+1};\quad  2\leq i \leq n-1.
\end{array}
\]
All other filiform Lie algebras can be obtained from it by deformations, see \cite{V}. 

The simplest filiform Lie superalgebra, $L^{n,m}$, is defined by the only non-zero
bracket products:
\[
\begin{array}{lcl}
[X_1,X_i] & = & X_{i+1};\quad  2\leq i \leq n-1\\[2mm]
[X_1, Y_j ] & = & Y_{j+1}; \quad 1\leq j \leq m-1
\end{array}
\]
where $p(X_i)=\ev$, for $j=1,\ldots,n$, and $p(Y_j)=\od$, for $j=1,\ldots, m$. This filiform Lie superalgebra is in complete analogy with filiform Lie algebras, since all the other filiform Lie superalgebras can be obtained from it by deformations, see \cite{BGKN}.

\ssbegin[(Orthosymplectic and periplectic forms on $L^{n,m}$)]{Proposition} [Orthosymplectic and periplectic forms on $L^{n,m}$]\label{filiform}\textup{(}i\textup{)} For $n=m$, the Lie superalgebra $L^{n,m}$ admits a~closed anti-symmetric {\bf periplectic} from given by \textup{(}where $\lambda\mu\not=0$\textup{)} 
\[
\omega=\lambda X_1^*\wedge Y_{n}^*+ \mu \sum_{i=2}^{n} (-1)^{n-i} X_i^* \wedge Y_{n+1-i}^*.
\]
\textup{(}ii\textup{)} For $n$ even and $m$ odd, the Lie superalgebra $L^{n,m}$ admits a~closed anti-symmetric  {\bf ortho-symplectic} form given by \textup{(}$\lambda\mu\nu\not =0$\textup{)}
\[
\omega= \lambda X_1^*\wedge X_n^*+ \mu \sum_{i=2}^{\frac{n}{2}} (-1)^{\frac{n}{2}-i} X_i^* \wedge X_{n+1-i}^*+(-1)^{\frac{m+1}{2}} \nu \sum_{i=1}^{\frac{m-1}{2}} (-1)^iY_i^*\wedge Y_{m+1-i}^*+ \frac{\nu}{2} \, Y_\frac{m+1}{2}^*\wedge Y_\frac{m+1}{2}^*.
\]
\textup{(}iii\textup{)} For $n$ even, the Lie superalgebra $L^{n,n-2}$ admits a~closed anti-symmetric non-degenerate {\bf non-homogeneous} form given by \textup{(}where $\lambda\mu\not=0$\textup{)}
\[
\omega= \lambda X_1^*\wedge X_n^*- \mu  \sum_{i=2}^{n-1} (-1)^i X_i^* \wedge Y_{n-i}^*.
\]

\end{Proposition}
\begin{proof} 

We will only prove Part (i). The proof of Parts (ii) and (iii) are almost identical. Let $a=a_\ev+a_\od,b=b_\ev+b_\od,c=c_\ev+c_\od\in L^{n,m}$, and let us express them as (where $\lambda_i, \tilde \lambda_i, \mu_i, \tilde \mu_i, \gamma_i, \tilde \gamma_i\in \Bbb K$ for $i=1,\ldots,n$):
\[
a=\sum_{i=1}^n\lambda_i X_i + \sum_{i=1}^n\tilde \lambda_i Y_i,\quad b=\sum_{i=1}^n\mu_i X_i + \sum_{i=1}^n\tilde \mu_i Y_i, \quad c=\sum_{i=1}^n\gamma_i X_i + \sum_{i=1}^n\tilde \gamma_i Y_i,
\]
A direct computation shows that 
\[
\begin{array}{lcl}
[a,b]& = & \displaystyle  \sum_{i=2}^{n-1} (\lambda_1 \mu_i-\lambda_i\mu_1)X_{i+1}+ \sum_{i=1}^{n-1} (\lambda_1 \tilde \mu_i-\tilde \lambda_i\mu_1)Y_{i+1}.
\end{array}
\]
Obviously, the expressions $[b,c]$ and $[c,a]$ can be obtained from that of $[a,b]$ by changing the roles of $a$ and $b$. On the other hand, using the fact that $[L^{n,m}_\od, L^{n,m}_\od]=\{0\}$ and $\omega$ is odd, we have
\[
\begin{array}{lcl}
(-1)^{p(a)p(c)}\omega(a,[b,c]) & = & \omega(a_\ev, [b_\ev,c_\od])+ \omega (a_\ev, [b_\od, c_\ev])+ \omega (a_\od, [b_\ev, c_\ev])\\[2mm]
& =& \displaystyle \lambda \lambda_1 (\mu_1 \tilde \gamma_{n-1}- \tilde \mu_{n-1} \gamma_1)+ \mu \sum_{i=2}^{n-1} (-1)^{n-i} \lambda_i (\mu_1 \tilde \gamma_{n-i}- \tilde \mu_{n-i} \gamma_1)
\\[2mm]
&& \displaystyle - \mu \sum_{i=3}^n (-1)^{n-i} \tilde \lambda_{n+1-i} (\mu_1 \gamma_{i-1} - \mu_{i-1} \gamma_1).
\end{array}
\]
Now, it is a matter of direct computation to show that $(-1)^{p(a)p(c)}\omega(a,[b,c])+ \circlearrowleft (a,b,c)=0$.

The matrix associated with the form $\omega$ retains the shape 
\[
\mathrm{antidiag}(\lambda, (-1)^{n}\mu, (-1)^{n-1}\mu,\ldots, \mu, \ldots, -(-1)^{n-1}\mu, -(-1)^{n}\mu, -\lambda).
\]
This matrix is invertible if and only if $\lambda\mu\not=0.$
\end{proof}

\ssbegin[(Non-unicity)]{Remark}[Non-unicity]
The closed anti-symmetric non-degenerate forms given in Proposition \ref{filiform}  are not unique.
\end{Remark}

\ssbegin[(Checked for $n=4$ and $m=6,8$)]{Conjecture}[Checked for $n=4$ and $m=6,8$] \textup{(}i\textup{)} If $n$ and $m$ are both even and $m\not =n-2$, then the Lie superalgebra $L^{n,m}$ does not admit any  anti-symmetric non-degenerate form satisfyning the $2$-cocycle condition \eqref{coccond}.

\textup{(}ii\textup{)} The solvable Lie superalgebra $SL^{n,m}$ with model the filiform Lie superalgebra $L^{n,m}$, see \cite{BN}, does not admit any non-degenerate and anti-symmetric form satisfying the 2-cocycle condition \textup{(}\ref{coccond}\textup{)}.
\end{Conjecture}
\section{Orthosymplectic double extensions}
 \subsection{${\mathscr D}_\ev$-extensions}
Let $\fa$ be an~ orthosymplectic quasi-Frobenius Lie superalgebra. Let ${\mathscr D}\in \fder_\ev(\fa)$ be a derivation such that the following 2-cocycle
\begin{equation} \label{omegatrivial}
(a,b)\mapsto \Omega(a,b):= \omega_\fa( {\mathscr D} \circ {\mathscr D}(a) + 2 {\mathscr D}^* \circ {\mathscr D}(a)+ {\mathscr D}^* \circ {\mathscr D}^*(a) +\lambda (  {\mathscr D} +  {\mathscr D} ^*  )(a)  , b)  
\end{equation}
is in $B^2(\fg;\Kee)$ \textup{(}i.e., a 2-coboundary\textup{)}; namely, there exists $\chi \in C^1(\fg,\Kee)$ such that $\Omega(a,b)=\chi([a,b]_\fa)$. Since $\omega_\fa$ is non-degenerate, there exists $Z_\Omega\in \fa$ such that 
\begin{equation}\label{zomega}
\Omega(a,b)=\omega_\fa(Z_\Omega,[a,b]_\fa) \quad \text{ for all $a,b\in \fa$}.
\end{equation}
\sssbegin[(${\mathscr D}_\ev$-extension -- the case where $\omega$ is orthosymplectic)]{Theorem}[${\mathscr D}_\ev$-extension -- the case where $\omega$ is orthosymplectic]\label{MainTh} Let $\fa$ be an~ orthosymplectic quasi-Frobenius Lie superalgebra. Let ${\mathscr D}\in \fder_\ev(\fa)$ be a derivation such that condition \textup{(}\ref{omegatrivial}\textup{)} is satisfied. Then, there exists a~Lie superalgebra structure on $\fg:=\mathscr{K} \oplus \fa \oplus \mathscr{K} ^*$, where $\mathscr{K} :=\Span\{x\}$ for $x$ even, defined as follows \textup{(}for any $a, b\in \fa$\textup{)}:
\[
\begin{array}{l}
[x ,x^*]_\fg=\lambda x, \; [a,b]_\fg  :=  [a,b]_\fa +  \omega_\fa( {\mathscr D}(a)+{\mathscr D}^*(a), b)x, \; [x^*,a]_\fg :=  {\mathscr D}(a)+ \omega_\fa(Z_\Omega,a)x,
\end{array}
\]
where $\lambda\in \Kee$ and $Z_\Omega\in \fa$ is as in Eq. \eqref{zomega}. There exists a~closed anti-symmetric orthosymplectic form $\omega_\fg$ on $\fg$ defined as follows: 
\begin{eqnarray*}
{\omega_\fg}\vert_{\fa \times \fa}:= \omega_\fa, \quad \omega_\fg(\fa,\mathscr{K} ):=\omega_\fg(\fa,\mathscr{K} ^*):=0, \quad \omega_\fg(x^*,x):=1, \quad 
\omega_\fg(x,x):=\omega_\fg(x^*,x^*) :=0.
\end{eqnarray*}
\end{Theorem}

The orthosymplectic  quasi-Frobenius Lie superalgebra $(\fg, \omega_\fg)$ will be called the {\it double extension} by means of a 1-dimensional space of the  orthosymplectic  quasi-Frobenius Lie superalgebra $(\fa, \omega_\fa)$. In the particular case where $x$ is central in $\fg$, namely $\lambda=0$, then the double extension is called {\it classical} (see \cite{BaBe}).

\begin{proof} 

For simplicity and clear presentation, let us introduce the following notations 
\[
\begin{array}{lcl}
\phi(a,b) & := & \omega_\fa( {\mathscr D}(a), b) + \omega_\fa(a, {\mathscr D}(b)) \quad \text{for any $a,b\in \fa$}\\[2mm]
J_\fg(f,g,h)& := & (-1)^{p(f)p(h)}[f,[g,h]_\fg]_\fg+\circlearrowleft (f,g,h) \quad \text{for any $f,g,h\in \fg$}.
\end{array}
\]
To check the Jacobi identity, we proceed as follows. 

The identity $J_\fg(x,f,g)=0$ is certainly satisfied since $x$ is central in $\fa\oplus {\mathscr K}$ and $x^*\not \in [\fg,\fg]_\fg$. 

The identity $J_\fg(x^*,f,g)=0$ is also satisfied for the following reasons. If either $f$ or $g$ is $x$, then we are done. Now if $f=x^*$ (or the other way round, $g=x^*$), then
\[
J_\fg(x^*,x^*,g)=[x^*,[x^*,g]_\fg]_\fg-[x^*,[x^*,g]_\fg]_\fg=0.
\]
Let us assume now that $f, g \in\fa.$ We have
\[
\begin{array}{lcl}
J_\fg(x^*,f,g)&=&[x^*, \phi(a,b)x+[f,g]_\fa]_\fg+ (-1)^{p(f)p(g)} [g,{\mathscr D}(f)+\omega_\fa(Z_\Omega,f) x]_\fg\\[2mm]
&&-[f,{\mathscr D}(g)+\omega_\fa(Z_\Omega,g), x]_\fg\\[2mm]
&=&- \phi(a,b)\lambda x+{\mathscr D}([f,g]_\fa)+\omega_\fa(Z_\Omega,[f,g]_\fa)x\\[2mm]
&&+(-1)^{p(f)p(g)}\left (  [g,{\mathscr D}(f)]_\fa + \phi(g, {\mathscr D}(f))x \right )-[f,{\mathscr D}(g)]_\fa - \phi(f, {\mathscr D}(g))x=0
\end{array}
\]
since ${\mathscr D}\in\fder(\fa)$ and the thanks to Lemma \ref{lemma0}. 

\sssbegin[(Technical)]{Lemma}[Technical] \label{lemma0} For all $f,g\in \fa$, we have 
\[
\omega_\fa(Z_\Omega,[f,g]_\fa)+(-1)^{p(f)p(g)}\phi(g, {\mathscr D}(f))-\phi(f, {\mathscr D}(g)) - \lambda  \phi(a,b) =0.
\]
\end{Lemma}

\begin{proof}
Indeed, 
\begin{align*}
&-(-1)^{p(f)p(g)} \phi(g, {\mathscr D}(f))+\phi(f, {\mathscr D}(g)) + \lambda  \phi(a,b)\\[2mm]
 &=  \phi({\mathscr D}(f), g)+\phi(f, {\mathscr D}(g)) + \lambda \phi(a,b) \\[2mm]
&= \omega_\fa( {\mathscr D}( {\mathscr D}(f)),  g) + \omega_\fa({\mathscr D}(f), {\mathscr D}(g)) +  \omega_\fa( {\mathscr D}(f), {\mathscr D}(g)) + \omega_\fa(f, {\mathscr D}({\mathscr D}(b)))\\[2mm]
& + \lambda (\omega_\fa({\mathscr D}(a),b)+ \omega_\fa(a, {\mathscr D}(b)))\\[2mm]
&= \omega_\fa ( ({\mathscr D} \circ  {\mathscr D} + 2  {\mathscr D}^* \circ  {\mathscr D} +{\mathscr D}^*\circ  {\mathscr D}^* + \lambda {\mathscr D} +\lambda {\mathscr D}^*)(f),  g)= \omega_\fa(Z_\Omega, [f,g]_\fa).\qed
\end{align*}
\noqed\end{proof}

From now on we assume that $f,g, h \in \fa$. We have
\[
\begin{array}{lcl}J_\fg(f,g,h)&=&(-1)^{p(f)p(h)}[f, \phi(g,h)x+[g,h]_\fa]_\fg+\circlearrowleft (f,g,h)\\[2mm]
&=&(-1)^{p(f)p(h)}( [f,[g,h]_\fa]_\fa+\phi(f,[g,h]_\fa)x) +\circlearrowleft (f,g,h)\\[2mm]
&=&0,
\end{array}
\]
because the JI holds on $\fa$ and thanks to Lemma \ref{lemma1}.

\sssbegin[(Technical)]{Lemma}[Technical]
\label{lemma1}
Let $(\fa,\omega_\fa)$ be an~orthosymplectic quasi-Frobenius Lie superalgebra. Let ${\mathscr D}$ be in $\fder(\fa)$. Then, 
\[
(-1)^{p(f)p(h)}\phi(f,[g,h]_\fa)+ \circlearrowleft (f,g,h)=0 \text{ for any $f,g,h\in\fa$}.
\]

\end{Lemma}
\begin{proof}

The proof follows from the fact that $\omega_\fa$ is orthosymplectic and ${\mathscr D}$ is a derivation. Indeed, 
\[
\begin{array}{lcl}
(-1)^{p(f)p(h)}\phi(f,[g,h]_\fa)+\ \circlearrowleft (f,g,h) &=&(-1)^{p(f)p(h)} (\omega_\fa( {\mathscr D}(f), [g,h]_\fa) + \omega_\fa (f, {\mathscr D}([g,h]_\fa)) ) \\[2mm]
& &+ (-1)^{p(h)p(g)} (\omega_\fa( {\mathscr D}(h), [f,g]_\fa) + \omega_\fa (h, {\mathscr D}([f,g]_\fa)) ) \\[2mm]
& &+ (-1)^{p(g)p(f)} (\omega_\fa( {\mathscr D}(g), [h,f]_\fa) + \omega_\fa (g, {\mathscr D}([h,f]_\fa))).
\end{array}
\]
We will collect terms involving ${\mathscr D}(f)$ and shows that they vanish. The same applies to terms involving ${\mathscr D}(g)$ and ${\mathscr D}(h)$. Indeed, 
\[
\begin{array}{lcl}
(-1)^{p(f)p(h)} \omega_\fa( {\mathscr D}(f), [g,h]_\fa) + (-1)^{p(h)p(g)} \omega_\fa (h, [{\mathscr D}(f),g]_\fa)  +(-1)^{p(g)p(f)} \omega_\fa (g, [h, {\mathscr D}(f)]_\fa)=0,\\[2mm]
\end{array}
\]
since $\omega_\fa$ is closed. \end{proof}
Let us now show that $\omega_\fg$ is closed on $\fg$. We should check that
\begin{equation}\label{omegag}
(-1)^{p(f)p(h)}\omega_\fa(f,[g,h]_\fg)+\circlearrowleft (f,g,h)=0 \quad \text{for any $f,g$ and $h$  $\in \fg$.}
\end{equation}
If $f=x^*$ and $g=x$, but $h\in \fa$, we have
\[
\begin{array}{lcl}
(-1)^{p(f)p(h)}\omega_\fa(f,[g,h]_\fg)+\circlearrowleft (f,g,h)& = & \omega_\fg(x^*, [x,h]_\fg)+ \omega_\fg(h, [x^*,x]_\fg)+ \omega_\fg(x, [h,x^*]_\fg)\\[2mm]
& = &  \omega_\fg(h, - \lambda x)+ \omega_\fg(x, - {\mathscr D}(h)-\omega_\fa(Z_\Omega, h)x)=0. 
\end{array}
\]
If $f=g=x^*$, but $ h \in \fa$, we have
\[
\begin{array}{lcl}
(-1)^{p(f)p(h)}\omega_\fg(f,[g,h]_\fg)+\circlearrowleft (f,g,h)&=&\omega_\fg(x^*,[x^*,h]_\fg) +  \omega_\fa(h,[x^*,x^*]_\fg) + \omega_\fg(x^*,[h,x^*]_\fg)  \\[2mm]
&=&\omega_\fg(x^*,[x^*,h]_\fg) - \omega_\fg (x^*,[x^*,h]_\fg)=0
\end{array}
\]
If $f=x^*$, but $g, h \in \fa$, we have
\begin{align*}
&(-1)^{p(f)p(h)}\omega_\fg(f,[g,h]_\fg)+\circlearrowleft (f,g,h)\\[2mm]
& =\omega_\fg(x^*,[g,h]_\fg) +(-1)^{p(h)p(g)} \omega_\fg(h,[x^*,g]_\fg) + \omega_\fg(g,[h,x^*]_\fg)  \\[2mm]
&=\omega_\fg(x^*,[g,h]_\fa+ \phi(g,h)x) +(-1)^{p(h)p(g)} \omega_\fg(h, {\mathscr D}(g)+\omega_\fg(Z_\Omega,g)x)- \omega_\fg(g, {\mathscr D}(h)+\omega_\fa(Z_\Omega,h)x)  \\[2mm]
&=\phi(g,h) - \omega_\fa({\mathscr D}(g),h)- \omega_\fa(g, {\mathscr D}(h))=0.
\end{align*}
If $f=x$ but $g, h \in \fa$, Eq. (\ref{omegag}) is trivially satisfied because of the very definition of $\omega_\fg$.

If $f,g,h\in \fa$, then we have
\[
\begin{array}{lcl}
(-1)^{p(f)p(h)}\omega_\fa(f,[g,h]_\fg)+\circlearrowleft (f,g,h)&=&(-1)^{p(f)p(h)}\omega_\fg(f,[g,h]_\fa+\phi(f,g)x)+\circlearrowleft (f,g,h)  \\[2mm]
& = & (-1)^{p(f)p(h)}\omega_\fg(f,[g,h]_\fa)+\circlearrowleft (f,g,h)  \\[2mm]
& = & (-1)^{p(f)p(h)}\omega_\fa(f,[g,h]_\fa)+\circlearrowleft (f,g,h)=0,
\end{array}
\]
since $\omega_\fa$ is closed on $\fa$.\end{proof}

\sssbegin[(Converse of Theorem \ref{MainTh})]{Theorem}[Converse of Theorem \ref{MainTh}]
\label{Rec1}
Let $(\fg,\omega_\fg)$ be an~orthosymplectic quasi-Frobenius Lie superalgebra. Suppose there exists an even non-zero $x\in ([\fg, \fg]_\fg)^\perp$ such that $\mathscr{K}:=\Span\{x\}$ is an ideal. Then, $(\fg,\omega_\fg)$ is obtained as an ${\mathscr D}_\ev$-extension by a 1-dimensional space from an~orthosymplectic quasi-Frobenius Lie superalgebra $(\fa,\omega_\fa)$. Moreover, if $\fz_\ev(\fg)\not =0$, then we can choose $x \in \fz_\ev(\fg)$, so  the double extension is classical.

\end{Theorem}
\begin{proof}

The space $\mathscr{K}^\perp$ is an ideal in $(\fg,\omega_\fg)$. Indeed, let us show first that $(\mathscr{K}^\perp)_\od=\fg_\od$. This is true because $x$ is orthogonal to any odd element since $\omega_\fg$ is even. Let $f\in \mathscr{K}^\perp$ and $g\in \fg$, then $[f,g]_\fg \in \mathscr{K}^\perp$ because $ \omega_\fg(x,[f,g]_\fg)=0$ since $x\in [\fg, \fg]_\fg^\perp$.

Since $\mathscr{K}$ is 1-dimensional and $\omega_\fg(x,x)=0$, it follows that $\mathscr{K}\subset \mathscr{K}^\perp$ and $\dim(\mathscr{K}^\perp)=\dim(\fg)-1$. Therefore, there exists $x^* \in \fg_\ev$ (since $(\mathscr{K}^\perp)_\od=\fg_\od$) such that
\[
\fg=\mathscr{K}^\perp\oplus \mathscr{K}^*, \quad \text{ where  $\mathscr{K}^*:=\Span\{x^*\}$}.
\]
We can normalize $x^*$ so that $\omega_\fg(x^*,x)=1$. 

Let us define $\fa:=(\mathscr{K} +\mathscr{K}^*)^\perp$. We then have a decomposition $\fg=\mathscr{K} \oplus \fa \oplus \mathscr{K}^*$. 

Let us define an~orthosymplectic form on $\fa$ by setting: 
\[
\omega_\fa={\omega_\fg}\vert_{\fa\times \fa}.
\]  
The form $\omega_\fa$ is non-degenerate on $\fa$. Indeed, suppose there exists an $a\in\fa$ such that 
\[
\omega_\fa(a,f)=0\quad \text{ for any $f\in \fa$}.
\] 
But $a$ is also orthogonal to $x$ as well as to $x^*$. It follows that $\omega_\fg(a,f)=0$ for any $f\in\fg$. Hence, $a=0$, since $\omega_\fg$ is nondegenerate on $\fg$. 

Let us now show that the vector space $\fa$ can be endowed with a Lie superalgebra structure, and there exists an~orthosymplectic structure on $\fa$ for which $\fg$ is its double extension. Let $a,b \in \fa$. The bracket $[a,b]_\fg$ belongs to $\mathscr{K} \oplus \fa$ because $\fa \subset \mathscr{K}^\perp=\mathscr{K}\oplus \fa$ and the latter is an ideal. It follows that 
\[
[x,a]_\fg=\lambda(a)x, \quad [x, x^*]_\fg=\lambda\, x, \quad [a,b]_\fg=\psi(a,b)+\phi(a,b)x, \quad [x^*,a]_\fg={\mathscr D}(a)+\chi(a)x.
\]
where $\psi(a,b), {\mathscr D}(a) \in \fa$, and $\lambda, \lambda(a), \phi(a,b), \chi (a) \in \Bbb K$.
 Let us show that 
\[
\phi(a,b)=  \omega_\fa( {\mathscr D}(a), b) + \omega_\fa(a, {\mathscr D}(b)).
\] 
Indeed, since 
\[
\omega_\fg(x^*,[a,b]_\fg)= - (-1)^{p(a)p(b)}\omega_\fg(b,[x^*,a]_\fg)  - \omega_\fg(a,[b,x^*]_\fg),
\] 
it follows that 
\begin{equation}
\label{phi}
\begin{array}{lcl}
\omega_\fg(x^*,\phi(a,b)x+\psi(a,b)) & = & - (-1)^{p(a)p(b)}\omega_\fg(b, \chi(a)x+{\mathscr D}(a)) \\[2mm]
&& + \omega_\fg(a, \chi(b)x+{\mathscr D}(b)).
\end{array}
\end{equation}
Since $\psi(a,b)$ is orthogonal to $x^*$, and $a$ is orthogonal to $x$, and $\omega_\fg(x^*,x)=1$, it follows that the LHS of Eq. (\ref{phi}) gives $\phi(a,b)$ while the RHS gives
\[
 - (-1)^{p(a)p(b)}\omega_\fa(b, {\mathscr D}(a)) + \omega_\fa(a, {\mathscr D}(b)) = \omega_\fa( {\mathscr D}(a), b) + \omega_\fa(a, {\mathscr D}(b)).
\]
Let us show that $\lambda(a)=0$ for all $a\in 
\fa$. Indeed, since \[
\omega_\fg(x^*,[x,a]_\fg)= -\omega_\fg(a,[x^*,x]_\fg)  - \omega_\fg(x,[a,x^*]_\fg).
\]
The LHS reads
\[
\omega_\fg(x^*,[x,a]_\fg)= \omega_\fg(x^*, \lambda(a) x) =  \lambda(a),
\]
while the RHS reads 
\[
-\omega_\fg(a,[x^*,x]_\fg)  - \omega_\fg(x,[a,x^*]_\fg)=-\omega_\fg(a, \lambda x)  + \omega_\fg(x,{\mathscr D}(a)+\chi(a)x)=0+0=0.
\]
The map $\psi$ is a bilinear on $\fa$ because $[\cdot , \cdot]_\fg$ is bilinear, so for convenience let us re-denote $\psi$ by $[\cdot , \cdot]_\fa$. We still have to show that $[\cdot , \cdot]_\fa$ satisfies the Jacobi identity. 

For any $a,b \in\fa$, the Jacobi identity on $\fg$  reads
\[
\begin{array}{lcl}
[x^*, [a,b]_\fg]_\fg & = & -(-1)^{p(a)p(b)}[b, [x^*,a]_\fg]_\fg - [a, [b, x^*]_\fg]_\fg=  [[x^*,a]_\fg, b]_\fg + [a, [x^*,b]_\fg]_\fg.
\end{array}
\]
The LHS reads
\[
\begin{array}{lcl}
[x^*, [a,b]_\fg]_\fg & =& [x^*, [a,b]_\fa+ \phi(a,b)x]_\fg = {\mathscr D}([a,b]_\fa)+ \chi([a,b]_\fa)x - \phi(a,b) \lambda x.
\end{array}
\]
The RHS reads
\[
\begin{array}{lcl}
 [[x^*,a]_\fg, b]_\fg + [a, [x^*,b]_\fg]_\fg & = & [ {\mathscr D}(a)+\chi(a)x, b]_\fg + [a, {\mathscr D}(b)+\chi(b)x]_\fg\\[2mm]
& = &[{\mathscr D}(a), b]_\fg + [a, {\mathscr D}(b)]_\fg\\[2mm]
& = &[{\mathscr D}(a), b]_\fa +\phi({\mathscr D}(a),b)x+ [a, {\mathscr D}(b)]_\fa+ \phi(a, {\mathscr D}(b))x.\\[2mm]
\end{array}
\]
By comparaison, it follows that ${\mathscr D}([a,b]_\fa)=[{\mathscr D}(a),b]_\fa+[a,{\mathscr D}(b)]_\fa$; this would imply that ${\mathscr D}$ is a derivation once we show that $[\cdot,\cdot]_\fa$ is a Lie bracket on $\fa$. Moreover,
\[
\begin{array}{lcl}
\chi([a,b]_\fa) & = & \phi({\mathscr D}(a),b)+  \phi(a, {\mathscr D}(b)) + \lambda \phi(a,b)\\[2mm]
& = &  \omega_\fa( {\mathscr D}({\mathscr D}(a)), b) + \omega_\fa({\mathscr D}(a), {\mathscr D}(b))+ \omega_\fa( {\mathscr D}(a), {\mathscr D}(b)) + \omega_\fa(a, {\mathscr D}({\mathscr D}(b))) \\[2mm]
&&- \lambda \phi(a,b)\\[2mm]
&=&   \omega_\fa( {\mathscr D} \circ {\mathscr D}(a) + 2 {\mathscr D}^* \circ {\mathscr D}(a)+ {\mathscr D}^* \circ {\mathscr D}^*(a) + \lambda ({\mathscr D}+{\mathscr D}^*)(a)  , b) .
\end{array}
\]
Therefore, the map 
\[
(a,b)\mapsto   \omega_\fa( {\mathscr D} \circ {\mathscr D}(a) + 2 {\mathscr D}^* \circ {\mathscr D}(a)+ {\mathscr D}^* \circ {\mathscr D}^*(a) + \lambda ({\mathscr D}+{\mathscr D}^*)(a) , b) 
\]
is a 2-coboundary (again, this will be justified when we show that $(\fa, [\cdot, \cdot]_\fa)$ is a Lie superalgebra).

For any $a,b$ and $c$ in $\fa$, we have
\[
\begin{array}{lcl}
0&=&(-1)^{p(a)p(c)}[a,[b,c]_\fg]_\fg+ \circlearrowleft (a,b,c)\\[2mm]
&=&(-1)^{p(a)p(c)}[a,[b,c]_\fa + \phi (b,c)x]_\fg + \circlearrowleft (a,b,c)\\[2mm]
&=&(-1)^{p(a)p(c)} [a,[b,c]_\fa]_\fa+ (-1)^{p(a)p(c)} \phi(a, [b,c]_\fa)x+ \circlearrowleft (a,b,c) \end{array}
\]
This implies that $[.,.]_\fa$ satisfies the Jacobi identity and 
\begin{equation} \label{DsymplO}
(-1)^{p(a)p(c)} \phi(a, [b,c]_\fa) + (-1)^{p(c)p(b)} \phi(c, [a,b]_\fa)+ (-1)^{p(b)p(a)} \phi(b, [c,a]_\fa)=0.
\end{equation}
Before checking Eq. (\ref{DsymplO}) we will prove first that $\omega_\fa:={\omega_\fg}\vert_{\fa \times \fa}$ is closed. Indeed, for any $a,b,c\in\fa$, we have
\[
\begin{array}{lcl}
(-1)^{p(a)p(c)}\omega_\fa(a,[b,c]_\fa)+\circlearrowleft (a,b,c)&=& (-1)^{p(a)p(c)}\omega_\fg(a,[b,c]_\fa)+\circlearrowleft (a,b,c) \\[2mm]
&=& (-1)^{p(a)p(c)}\omega_\fg(a,[b,c]_\fg- \phi(a,b)x)+\circlearrowleft (a,b,c) \\[2mm]
&=& (-1)^{p(a)p(c)}\omega_\fg(a,[b,c]_\fg)+\circlearrowleft (a,b,c) \\[2mm]
& = & 0,
\end{array}
\]
since $\omega_\fg$ is closed. Now, Let us rewrite Eq. (\ref{DsymplO}); explicitly:
\[
\begin{array}{l}
(-1)^{p(a)p(c)} \left ( \omega_\fa( {\mathscr D}(a), [b,c]_\fa) + \omega_\fa(a, [{\mathscr D}(b),c]_\fa)+\omega_\fa(a, [b,{\mathscr D}(c)]_\fa) \right )\\[2mm]
+(-1)^{p(c)p(b)} \left ( \omega_\fa( {\mathscr D}(c), [a,b]_\fa) + \omega_\fa(c, [{\mathscr D}(a),b]_\fa)+\omega_\fa(c, [a,{\mathscr D}(b)]_\fa) \right )\\[2mm]
+(-1)^{p(b)p(a)} \left ( \omega_\fa( {\mathscr D}(b), [c,a]_\fa) + \omega_\fa(b, [{\mathscr D}(c),a]_\fa)+\omega_\fa(b, [c,{\mathscr D}(a)]_\fa) \right )=0.
\end{array}
\]
We see that all these terms will cancel because $\omega_\fa$ is closed. The proof now is complete.
\end{proof}
 \subsection{${\mathscr D}_\od$-extensions}
Let $(\fa, \omega_\fa)$ be an~orthosymplectic quasi-Frobenius  Lie superalgebra. Let ${\mathscr D}\in \fder_\od(\fa)$ be a derivation such that
\begin{equation}\label{omegaeo}
(a,b)\mapsto \Omega(a,b):= \omega_\fa( {\mathscr D} \circ {\mathscr D}(a) -  {\mathscr D}^* \circ {\mathscr D}^*(a)  , b) 
\end{equation}
is in $B^2(\fg;\Kee)$ \textup{(}i.e., a 2-coboundary\textup{)}. Let us write $\Omega= \delta (\chi)$ for some $\chi\in C^1(\fa; \Kee)$; namely, $\Omega(a,b)=\chi([a,b]_\fa)$. Since $\omega_\fa$ is non-degenerate, there exists $a\in \fa_\ev$ such that $\Omega(a,b)=\omega_\fa(a_0,[a,b]_\fa)$.

\sssbegin[(${\mathscr D}_\od$-extension  -- the case where $\omega$ is orthosymplectic)]{Theorem}[${\mathscr D}_\od$-extension  -- the case where $\omega$ is orthosymplectic] \label{MainThO} 
Let $(\fa, \omega_\fa)$ be an~ortho-symplectic quasi-Frobenius Lie superalgebra. Let ${\mathscr D}\in \fder_\od(\fa)$ be a~derivation such that condition \eqref{omegaeo} is satisfied together with
\[
\begin{array}{l}
{\mathscr D}^2=\ad_{a_0} \quad \text{ and }\quad  {\mathscr D}(a_0)=0.
\end{array}
\]
Then, there exists a Lie superalgebra structure on $\fg:=\mathscr{K} \oplus \fa \oplus \mathscr{K} ^*$, where $\mathscr{K} :=\Span\{x\}$ for $x$ odd, defined as follows: \textup{(}for any $a, b\in \fa$\textup{)}
\[
\begin{array}{ll}
[\mathscr{K} ,\fg]_\fg=0, & [a,b]_\fg  :=  [a,b]_\fa +  \left (\omega_\fa( {\mathscr D}(a), b) +(-1)^{p(a)} \omega_\fa(a, {\mathscr D}(b))\right )x, \\[2mm]
[x^*, x^*]_\fg=2 a_0, & [x^*,a]_\fg :=  {\mathscr D}(a) - \omega_\fa(a,a_0)x.
\end{array}
\]
There exists a~closed anti-symmetric orthosymplectic form $\omega_\fg$ on $\fg$ defined as follows: 
\begin{eqnarray*}
{\omega_\fg}\vert_{\fa \times \fa}:= \omega_\fa, \quad \omega_\fg(\fa,\mathscr{K} ):=\omega_\fg(\fa,\mathscr{K} ^*):=0, \quad \omega_\fg(x^*,x):=1, \quad 
\omega_\fg(x,x):=\omega_\fg(x^*,x^*) :=0.
\end{eqnarray*}

\end{Theorem}
\begin{proof}

As we did in proving Theorem \ref{MainTh}, let us introduce the following notations:
\[
\begin{array}{lcl}
\phi(a,b) & := & \omega_\fa( {\mathscr D}(a), b) +(-1)^{p(a)} \omega_\fa(a, {\mathscr D}(b)) \quad \text{for any $a,b\in \fa$}\\[2mm]
J_\fg(f,g,h)& := & (-1)^{p(f)p(h)}[f,[g,h]_\fg]_\fg+\circlearrowleft (f,g,h) \quad \text{for any $f,g,h\in \fg$}.
\end{array}
\]
To check the Jacobi identity, we proceed as follows. 

If $h=x$, the identity $J_\fg(x,f,g)=0$ is certainly satisfied since $x$ is central in $\fg$. 

If $h=x^*$, the identity $J_\fg(x^*,f,g)=0$ is also satisfied for the following reasons. If either $f$ or $g$ is $x$, then we are done since $x$ is central in $\fg$. Now if $f=x^*$ (or the other way round, $g=x^*$), then
\[
\begin{array}{lcl}
J_\fg(x^*,x^*,g) & = & (-1)^{p(g)}[x^*,[x^*,g]_\fg]_\fg + (-1)^{p(g)}[g,[x^*,x^*]_\fg]_\fg -  [x^*,[g,x^*]_\fg]_\fg=0\\[2mm]
& = &  (-1)^{p(g)}[x^*,{\mathscr D}(g)-\omega_\fa(a,a_0)x]_\fg +(-1)^{p(g)} [g, 2 a_0]_\fg \\[2mm]
&&+  (-1)^{p(g)}[x^*,{\mathscr D}(g)-\omega_\fa(a,a_0)x]_\fg\\[2mm]
& = &  (-1)^{p(g)} {\mathscr D}({\mathscr D}(g))-\omega_\fa({\mathscr D}(g),a_0)x + (-1)^{p(g)}([g, 2 a_0]_\fa+\phi(g,2a_0)x) \\[2mm]
&&+  (-1)^{p(g)}  \left ( {\mathscr D}({\mathscr D}(g))-\omega_\fa({\mathscr D}(g),a_0)x \right )=0\\[2mm]
\end{array}
\]
since ${\mathscr D}^2=\ad_{a_0}$ and ${\mathscr D}(a_0)=0$.

Let us assume now that $f, g \in\fa.$ We have
\[
\begin{array}{lcl}
J_\fg(x^*,f,g)&=&(-1)^{p(g)}[x^*, \phi(f,g)x+[f,g]_\fa]_\fg+ (-1)^{p(f)p(g)} [g,{\mathscr D}(f)-\omega(f,a_0) x]_\fg\\[2mm]
&&- (-1)^{p(f)+p(g)}[f,{\mathscr D}(g)-\omega(g,a_0), x]_\fg \\[2mm]
&=&(-1)^{p(g)}({\mathscr D}([f,g]_\fa) - \omega_\fa([f,g]_\fa,a_0)x+(-1)^{p(f)p(g)}\left (  [g,{\mathscr D}(f)]_\fa + \phi(g, {\mathscr D}(f))x \right )\\[2mm]
&&- (-1)^{p(f)+p(g)}([f,{\mathscr D}(g)]_\fa + \phi(f, {\mathscr D}(g))x)=0
\end{array}
\]
since ${\mathscr D}\in\fder(\fa)$ and thanks to Lemma \ref{lemma01}.

\sssbegin[(Technical)]{Lemma}[Technical] \label{lemma3} \label{lemma01} For all $f,g\in \fa$, we have 
\[
\omega_\fa([f,g]_\fa,a_0)-(-1)^{(p(f)+1)p(g)}\phi(g, {\mathscr D}(f)) +  (-1)^{p(f)} \phi(f, {\mathscr D}(g))=0.
\]

\end{Lemma}
\begin{proof}

Indeed, 
\begin{align*}
& -(-1)^{(p(f)+1)p(g)} \phi(g, {\mathscr D}(f)) +  (-1)^{p(f)} \phi(f, {\mathscr D}(g))  =   \phi({\mathscr D}(f), g) + (-1)^{p(f)}\phi(f, {\mathscr D}(g))\\[2mm]
 & =  \omega_\fa( {\mathscr D}( {\mathscr D}(f)),  g) +(-1)^{p(f)+1} \omega_\fa({\mathscr D}(f), {\mathscr D}(g))\\[2mm]
& +(-1)^{p(f)}(  \omega_\fa( {\mathscr D}(f), {\mathscr D}(g))+(-1)^{p(f)} \omega_\fa(f, {\mathscr D}({\mathscr D}(g))))\\[2mm]
&= \omega_\fa ( ({\mathscr D} \circ  {\mathscr D} - {\mathscr D}^*\circ  {\mathscr D}^*)(f),  g)= \omega_\fa(a_0, [f,g]_\fa).\qed
\end{align*}
\noqed\end{proof}

From now on we assume that $f,g, h \in \fa$. We have
\[
\begin{array}{lcl}J_\fg(f,g,h)&=&(-1)^{p(f)p(h)}[f, \phi(g,h)x+[g,h]_\fa]_\fg+\circlearrowleft (f,g,h)\\[2mm]
&=&(-1)^{p(f)p(h)}( [f,[g,h]_\fa]_\fa+\phi(f,[g,h]_\fa)x) +\circlearrowleft (f,g,h),\\[2mm]
&=&0,
\end{array}
\]
because the JI holds on $\fa$ and thanks to Lemma \ref{lemma4}.

\sssbegin[(Technical)]{Lemma}[Technical]
\label{lemma4}
Let $(\fa,\omega_\fa)$ be an~orthosymplectic quasi-Frobenius Lie superalgebra. Let ${\mathscr D}$ be in $\fder(\fa)$. Then 
\[
(-1)^{p(f)p(h)}\phi(f,[g,h]_\fa)+ \circlearrowleft (f,g,h)=0 \text{ for any $f,g,h\in\fa$}.
\]

\end{Lemma}
\begin{proof} 

Similar to that of Lemma \ref{lemma1} taking into consideration that ${\mathscr D}$ is an odd derivation. 
\end{proof}
Let us now show that $\omega_\fg$ is closed.  We should check that
\[
(-1)^{p(f)p(h)}\omega_\fa(f,[g,h]_\fg)+\circlearrowleft (f,g,h)=0 \quad \text{for any $f,g$ and $h$  $\in \fg$.}
\]
This is true if $f=x$ because $x$ is central and by the very definition of $\omega_\fg$.  If $f=x^*$, but $g, h \in \fa$, we have
\begin{align*}
& (-1)^{p(f)p(h)}\omega_\fg(f,[g,h]_\fg)+\circlearrowleft (f,g,h)\\[2mm]
&=(-1)^{p(h)}\omega_\fg(x^*,[g,h]_\fg) +(-1)^{p(h)p(g)} \omega_\fg(h,[x^*,g]_\fg) +(-1)^{p(g)} \omega_\fg(g,[h,x^*]_\fg)  \\[2mm]
&=(-1)^{p(h)}\omega_\fg(x^*,[g,h]_\fa+ \phi(g,h)x) +(-1)^{p(h)p(g)} \omega_\fg(h, {\mathscr D}(g) - \omega_\fg(g,a_0)x)\\[2mm]
&- (-1)^{p(g)+p(h)} \omega_\fg(g, {\mathscr D}(h)- \omega_\fa(h,a_0)x)  \\[2mm]
&= (-1)^{p(h)}\phi(g,h) - (-1)^{p(h)}\omega_\fa({\mathscr D}(g),h) - (-1)^{p(g)+p(h)} \omega_\fa(g, {\mathscr D}(h))=0.
\end{align*}
 If $f=g=x^*$, but $ h \in \fa$, we have
\[
\begin{array}{lcl}
(-1)^{p(f)p(h)}\omega_\fg(f,[g,h]_\fg)+\circlearrowleft (f,g,h)&=&(-1)^{p(h)}\omega_\fg(x^*,[x^*,h]_\fg) +(-1)^{p(h)}  \omega_\fa(h,[x^*,x^*]_\fg)\\[2mm]
 && - \omega_\fg(x^*,[h,x^*]_\fg)  \\[2mm]
&=&-2 (-1)^{p(h)} \omega_\fa(h,a_0)+(-1)^{p(h)} \omega_\fa(h, 2 a_0)=0
\end{array}
\]
If $f,g,h\in \fa$, then we have
\[
\begin{array}{lcl}
(-1)^{p(f)p(h)}\omega_\fa(f,[g,h]_\fg)+\circlearrowleft (f,g,h)&=&(-1)^{p(f)p(h)}\omega_\fg(f,[g,h]_\fa+\phi(f,g)x)+\circlearrowleft (f,g,h)  \\[2mm]
& = & (-1)^{p(f)p(h)}\omega_\fg(f,[g,h]_\fa)+\circlearrowleft (f,g,h)  \\[2mm]
& = & (-1)^{p(f)p(h)}\omega_\fa(f,[g,h]_\fa)+\circlearrowleft (f,g,h)=0,
\end{array}
\]
since $\omega_\fa$ is closed on $\fa$.
\end{proof}

\sssbegin[(Converse of Theorem \ref{MainThO})]{Theorem}[Converse of Theorem \ref{MainThO}]
\label{Rec2}
Let $(\fg,\omega_\fg)$ be an irreducible orthosymplectic quasi-Frobenius Lie superalgebra. Suppose there exists a non-zero $x \in \fz(\fg)_\od$ such that $\omega(x,x)=~0$. Then, $(\fg,\omega_\fg)$ is a ${\mathscr D}_\od$-extension of an~orthosymplectic quasi-Frobenius Lie superalgebra $(\fa,\omega_\fa)$. 

\end{Theorem}

\ssbegin[(The condition $\omega(x,x)=0$)]{Remark}[The condition $\omega(x,x)=0$] The condition $\omega(x,x)=0$ is actually necessary. A counterexample for the Lie superalgebra $C^1_{1/2}+A$ -- see Table \ref{tab2} -- will be provided in Proposition \ref{centdame}.

\end{Remark}

\begin{proof}[Proof of Theorem \ref{Rec2} ]

Let $x$ be a non-zero element in $\fz(\fg)_\od$. The subspace $\mathscr{K}:=\Span\{x\}$ is an ideal in $\fg$  because $x$ is central in $\fg$. Moreover, $\mathscr{K}^\perp$ is also an ideal in $(\fg,\omega_\fg)$. Indeed, let us show first that $(\mathscr{K}^\perp)_\ev=\fg_\ev$. This is true because $x$ is orthogonal to any even  element since $\omega_\fg$ is even. Let $f\in \mathscr{K}^\perp$ and $g\in \fg$, then $[f,g]_\fg \in \mathscr{K}^\perp$ because
\[
(-1)^{p(g)}\omega_\fg(x,[f,g]_\fg)=-(-1)^{p(g)p(f)} \omega_\fg(g, [x,f]_\fg) - (-1)^{p(f)}\omega_\fg(f, [g,x]_\fg)=0-0=0.
\]
Since $\mathscr{K}$ is 1-dimensional and $\omega(x,x)=0$, it follows that $\mathscr{K}\subset \mathscr{K}^\perp$ and $\dim(\mathscr{K}^\perp)=\dim(\fg)-1$. Therefore, there exists $x^* \in \fg_\od$ such that
\[
\fg=\mathscr{K}^\perp\oplus \mathscr{K}^*, \quad \text{ where  $\mathscr{K}^*:=\Span\{x^*\}$}.
\]
This $x^*$ can be normalized to have $\omega_\fg(x^*,x)=1$. 

Let us define $\fa:=(\mathscr{K} +\mathscr{K}^*)^\perp$. We then have a decomposition $\fg=\mathscr{K} \oplus \fa \oplus \mathscr{K}^*$. 

Let us define an~orthosymplectic form on $\fa$ by setting: 
\[
\omega_\fa={\omega_\fg}\vert_{\fa\times \fa}.
\]  
The form $\omega_\fa$ is non-degenerate on $\fa$. Indeed, suppose there exists an $a\in\fa$ such that 
\[
\omega_\fa(a,f)=0\quad \text{ for any $f\in \fa$}.
\] 
But $a$ is also orthogonal to $x$ and $x^*$. It follows that $\omega_\fg(a,f)=0$ for any $f\in\fg$. Hence, $a=0$, since $\omega_\fg$ is nondegenerate. 

Let us show now that there exists an~ orthosymplectic structure on the vector space $\fa$ and $\fg$ is its double extension. Let $a,b \in \fa$. The bracket $[a,b]_\fg$ belongs to $\mathscr{K} \oplus \fa$ because $\fa \subset \mathscr{K}^\perp=\mathscr{K}\oplus \fa$ and the latter is an ideal. It follows that 
\[
[a,b]_\fg=\phi(a,b)x+\psi(a,b), \quad [x^*,a]_\fg=\chi(a)x+{\mathscr D}(a),\quad [x^*,x^*]_\fg=\nu x+2 a_0  .
\]
where $a_0, \psi(a,b), {\mathscr D}(a) \in \fa$, and $\phi(a,b), \chi (a), \nu \in \Bbb K$. Let us show that 
\begin{equation}
\label{phioBB}
\phi(a,b)=  \omega_\fa( {\mathscr D}(a), b) + (-1)^{p(a)}\omega_\fa(a, {\mathscr D}(b)).
\end{equation}
Indeed, since 
\[
(-1)^{p(b)}\omega_\fg(x^*,[a,b]_\fg)= - (-1)^{p(a)p(b)}\omega_\fg(b,[x^*,a]_\fg)  -(-1)^{p(a)} \omega_\fg(a,[b,x^*]_\fg),
\] 
it follows that 
\begin{equation}
\label{phioM}
\begin{array}{lcl}
(-1)^{p(b)}\omega_\fg(x^*,\phi(a,b)x+\psi(a,b)) & = & - (-1)^{p(a)p(b)}\omega_\fg(b, \chi(a)x+{\mathscr D}(a)) \\[2mm]
&&  +(-1)^{p(a)+p(b)} \omega_\fg(a, \chi(b)x+{\mathscr D}(b))
\end{array}
\end{equation}
Since $x^*$ is perpendicular to $\psi(a,b)$, and $a$ is orthogonal to $x$, and $\omega_\fg(x^*,x)=1$, it follows that Eq. (\ref{phioBB}) holds.

Similarly, since 
\[
(-1)^{p(a)}\omega_\fg(x^*,[x^*,a]_\fg)= - (-1)^{p(a)}\omega_\fg(a,[x^*,x^*]_\fg)  + \omega_\fg(x^*,[a,x^*]_\fg),
\] 
it follows that 
\begin{equation}
\label{phiod}
\begin{array}{lcl}
(-1)^{p(a)}\omega_\fg(x^*,\chi(a)x+{\mathscr D}(a)) & = & - (-1)^{p(a)}\omega_\fg(a,2 a_0 + \nu  x)  - (-1)^{p(a)}\omega_\fg(x^*, \chi(a)x+{\mathscr D}(a))
\end{array}
\end{equation}
Since ${\mathscr D}(a)$ is orthogonal to $x^*$, and $a$ is orthogonal to $x$, and $\omega_\fg(x^*,x)=1$, it follows that the LHS of Eq. (\ref{phiod}) gives $(-1)^{p(a)}\chi(a)$ while the RHS gives
\[
 - (-1)^{p(a)}\omega_\fa(a, 2a_0) - (-1)^{p(a)} \chi(a).
\]
It follows that
\begin{equation}\label{chi}
\chi(a)=-\omega_\fa(a,a_0).
\end{equation}
Besides, since $a_0$ is orthogonal with $x^*$, then
\[
\begin{array}{l}
0=\omega_\fg(x^*,[x^*,x^*]_\fg)=\omega_\fg(x^*,a_0+\nu x)=\nu
\end{array}
\]
Therefore, $\nu=0$.

The map $\psi$ is a bilinear on $\fa$ because $[\cdot , \cdot]_\fg$ is bilinear, so for convenience let us re-denote $\psi$ by $[\cdot , \cdot]_\fa$. We still have to show that $[\cdot , \cdot]_\fa$ satisfies the Jacobi identity.

Let us observe first that 
\[
\begin{array}{lcl}
0& = & [x^*,[x^*,x^*]_\fg]_\fg= [x^*, 2a_0+ \nu x]_\fg=2 ({\mathscr D}(a_0)+\chi(a_0)x),
\end{array}
\]
which implies that ${\mathscr D}(a_0)=0$ and $\chi(a_0)=0$ (which is certainly satisfied because of Eq. (\ref{chi})). 

For any $a \in\fa$, the Jacobi identity on $\fg$ reads
\[
(-1)^{p(a)}[x^*, [x^*,a]_\fg]_\fg=-(-1)^{p(a)}[a, [x^*,x^*]_\fg]_\fg + [x^*, [a, x^*]_\fg]_\fg .
\]
The LHS reads
\[
\begin{array}{lcl}
(-1)^{p(a)}[x^*, [x^*,a]_\fg]_\fg & = & (-1)^{p(a)}[x^*, {\mathscr D}(a)+\chi(a) x]_\fg =  (-1)^{p(a)} ( {\mathscr D}( {\mathscr D}(a))+\chi( {\mathscr D}(a)) x).
\end{array}
\]
The RHS reads
\[
\begin{array}{lcl}
-(-1)^{p(a)}[a, [x^*,x^*]_\fg]_\fg + [x^*, [a, x^*]_\fg]_\fg & = & -(-1)^{p(a)}[a, 2 a_0 + \nu x]_\fg - (-1)^{p(a)} [x^*, {\mathscr D}(a)+\nu x]_\fg \\[2mm]
& = & -(-1)^{p(a)}([a, 2 a_0 ]_\fa+ \phi(a, 2 a_0)x) \\[2mm]
& & - (-1)^{p(a)}  ({\mathscr D}({\mathscr D}(a))+\chi({\mathscr D}(a)) x). \\[2mm]
\end{array}
\]
It follows that ${\mathscr D}^2=\ad_{a_0}$ and $\chi({\mathscr D}(a))=-\phi(a,a_0)$ (which is true because ${\mathscr D}(a_0)=0$ and Eqs.  (\ref{phioBB}), (\ref{chi}).  )

For any $a,b \in\fa$, the Jacobi identity reads
\[
(-1)^{p(b)}[x^*, [a,b]_\fg]_\fg=-(-1)^{p(a)p(b)}[b, [x^*,a]_\fg]_\fg - (-1)^{p(a)} [a, [b, x^*]_\fg]_\fg .
\]
The LHS reads
\[
\begin{array}{lcl}
(-1)^{p(b)}[x^*, [a,b]_\fg]_\fg & =& (-1)^{p(b)} [x^*, [a,b]_\fa+ \phi(a,b)x]_\fg =(-1)^{p(b)}( {\mathscr D}([a,b]_\fa)+ \chi([a,b]_\fa)x).
\end{array}
\]
The RHS reads
\begin{align*}
&-(-1)^{p(a)p(b)}[b, [x^*,a]_\fg]_\fg - (-1)^{p(a)} [a, [b, x^*]_\fg]_\fg \\[2mm]
&= -(-1)^{p(a)p(b)}[b, {\mathscr D}(a)+\chi(a)x]_\fg +(-1)^{p(a)+p(b)} [a, {\mathscr D}(b)+\chi(b)x]_\fg\\[2mm]
& =  (-1)^{p(b)} [{\mathscr D}(a), b]_\fg + (-1)^{p(a)+p(b)} [a, {\mathscr D}(b)]_\fg\\[2mm]
& =  (-1)^{p(b)}([{\mathscr D}(a), b]_\fa +\phi({\mathscr D}(a),b)x)+(-1)^{p(a)+p(b)} ( [a, {\mathscr D}(b)]_\fa+ \phi(a, {\mathscr D}(b))x).\\[2mm]
\end{align*}
By comparaison, it follows that ${\mathscr D}([a,b]_\fa)=[{\mathscr D}(a),b]_\fa+(-1)^{p(a)}[a,{\mathscr D}(b)]_\fa$, namely ${\mathscr D}$ is a derivation, and 
\[
\begin{array}{lcl}
\chi([a,b]_\fa) & = & \phi({\mathscr D}(a),b)+ (-1)^{p(a)} \phi(a, {\mathscr D}(b))\\[2mm]
& = &  \omega_\fa( {\mathscr D}({\mathscr D}(a)), b) +(-1)^{p(a)+1} \omega_\fa({\mathscr D}(a), {\mathscr D}(b))\\[2mm]&&+(-1)^{p(a)}( \omega_\fa( {\mathscr D}(a), {\mathscr D}(b)) +(-1)^{p(a)} \omega_\fa(a, {\mathscr D}({\mathscr D}(b)))) \\[2mm]
&=&   \omega_\fa( ({\mathscr D} \circ {\mathscr D}- {\mathscr D}^* \circ {\mathscr D}^*)(a)  , b) .
\end{array}
\]
Therefore, the map 
\[
(a,b)\mapsto   \omega_\fa( {\mathscr D} \circ {\mathscr D}(a) - {\mathscr D}^* \circ {\mathscr D}^*(a)  , b)
\]
is a 2-coboundary.

For any $a,b$ and $c$ in $\fa$, we have
\[
\begin{array}{lcl}
0&=&(-1)^{p(a)p(c)}[a,[b,c]_\fg]_\fg+ \circlearrowleft (a,b,c)\\[2mm]
&=&(-1)^{p(a)p(c)}[a,[b,c]_\fa + \phi (b,c)x]_\fg + \circlearrowleft (a,b,c)\\[2mm]
&=&(-1)^{p(a)p(c)} [a,[b,c]_\fa]_\fa+ (-1)^{p(a)p(c)} \phi(a, [b,c]_\fa)x+ \circlearrowleft (a,b,c) \end{array}
\]
This implies that the bracket $[.,.]_\fa$ satisfied the Jacobi identity, and 
\begin{equation} \label{Dsympl}
(-1)^{p(a)p(c)} \phi(a, [b,c]_\fa) + (-1)^{p(c)p(b)} \phi(c, [a,b]_\fa)+ (-1)^{p(b)p(a)} \phi(b, [c,a]_\fa)=0.
\end{equation}
Before checking Eq. (\ref{Dsympl}) we will prove first that $\omega_\fa:={\omega_\fg}\vert_{\fa \times \fa}$ is closed.  Indeed, for any $a,b,c\in\fa$, we have
\[
\begin{array}{lcl}
(-1)^{p(a)p(c)}\omega_\fa(a,[b,c]_\fa)+\circlearrowleft (a,b,c)&=& (-1)^{p(a)p(c)}\omega_\fg(a,[b,c]_\fa)+\circlearrowleft (a,b,c) \\[2mm]
&=& (-1)^{p(a)p(c)}\omega_\fg(a,[b,c]_\fg- \phi(a,b)x)+\circlearrowleft (a,b,c) \\[2mm]
&=& (-1)^{p(a)p(c)}\omega_\fg(a,[b,c]_\fg)+\circlearrowleft (a,b,c) \\[2mm]
& = & 0,
\end{array}
\]
since $\omega_\fg$ is closed.  Now, let us rewrite Eq. (\ref{Dsympl}); explicitly:
\[
\begin{array}{l}
(-1)^{p(a)p(c)} \left ( \omega_\fa( {\mathscr D}(a), [b,c]_\fa) +(-1)^{p(a)}  \omega_\fa(a, [{\mathscr D}(b),c]_\fa)+(-1)^{p(a)+p(b)} \omega_\fa(a, [b,{\mathscr D}(c)]_\fa) \right )\\[2mm]
+(-1)^{p(c)p(b)} \left ( \omega_\fa( {\mathscr D}(c), [a,b]_\fa) + (-1)^{p(c)}  \omega_\fa(c, [{\mathscr D}(a),b]_\fa)+(-1)^{p(c)+p(a)}\omega_\fa(c, [a,{\mathscr D}(b)]_\fa) \right )\\[2mm]
+(-1)^{p(b)p(a)} \left ( \omega_\fa( {\mathscr D}(b), [c,a]_\fa) + (-1)^{p(b)}  \omega_\fa(b, [{\mathscr D}(c),a]_\fa)+ (-1)^{p(b)+p(c)}\omega_\fa(b, [c,{\mathscr D}(a)]_\fa) \right )=0.
\end{array}
\]
We see that all these terms cancel because $\omega_\fa$ is closed. The proof now is complete.\end{proof}
\section{Periplectic double extensions}
\subsection{${\mathscr D}_\ev$-extension}
Let $(\fa, \omega_\fa)$ be a~periplectic quasi-Frobenius Lie superalgebra. Let ${\mathscr D}\in \fder_\ev(\fa)$ be a derivation such that the following 2-cocycle
\begin{equation}\label{OmTriOE}
(a,b)\mapsto \Omega(a,b):= \omega_\fa \left ( ({\mathscr D} \circ {\mathscr D}+ 2 {\mathscr D}^* \circ {\mathscr D}+ {\mathscr D}^* \circ {\mathscr D}^* +\lambda  {\mathscr D} + \lambda  {\mathscr D} ^* )(a) , b \right )  
\end{equation}
is in $B^2(\fg;\Kee)$ \textup{(}i.e., a 2-coboundary\textup{)}. Let us write $\Omega= \delta (\chi)$ for some $\chi\in C^1(\fa; \Kee)$; namely, $\Omega(a,b)=\chi([a,b]_\fa)$. Since $\omega_\fa$ is non-degenerate, there exists $Z_\Omega\in \fa$ such that $\Omega(a,b)=\omega_\fa(Z_\Omega,[a,b]_\fa)$.

\sssbegin[(${\mathscr D}_\ev$-extension -- the case where $\omega$ is periplectic)]{Theorem}[${\mathscr D}_\ev$-extension -- the case where $\omega$ is periplectic ]\label{MainThOE} Let $(\fa, \omega_\fa)$ be a~periplectic quasi-Frobenius Lie superalgebra. Let ${\mathscr D}\in \fder_\ev(\fa)$ be a derivation satisfying condition \eqref{OmTriOE}. Then, there exists a~Lie superalgebra structure on $\fg:=\mathscr{K} \oplus \fa \oplus \mathscr{K} ^*$, where $\mathscr{K} :=\Span\{x\}$ for $x$ odd and $\mathscr{K} :=\Span\{e\}$ for $e$ even, defined as follows: \textup{(}for any $a, b\in \fa$\textup{)}
\[
\begin{array}{l}
[x ,e]_\fg=\lambda x, \; [a,b]_\fg  :=  [a,b]_\fa +  \left (\omega_\fa( {\mathscr D}(a), b) + \omega_\fa(a, {\mathscr D}(b))\right )x, \quad [e,a]_\fg :=  {\mathscr D}(a)+ \omega_\fa(Z_\Omega,a)x.
\end{array}
\]
There exists a~closed anti-symmetric periplectic form $\omega_\fg$ on $\fg$ defined as follows: 
\begin{eqnarray*}
{\omega_\fg}\vert_{\fa \times \fa}:= \omega_\fa, \quad \omega_\fg(\fa,\mathscr{K} ):=\omega_\fg(\fa,\mathscr{K} ^*):=0, \quad \omega_\fg(e,x):=1, \quad 
\omega_\fg(x,x):=\omega_\fg(e,e) :=0.
\end{eqnarray*}

\end{Theorem}
\begin{proof}

Similar to that of Theorem \ref{MainTh}\end{proof}

The converse of Theorem \ref{MainThOE} is given by the following theorem. 

\sssbegin[(Converse of Theorem \ref{MainThOE})]{Theorem}[Converse of Theorem \ref{MainThOE}]
\label{Rec3}
Let $(\fg,\omega_\fg)$ be a~periplectic quasi-Frobenius Lie superalgebra. Suppose that there exists  $0\not = x \in ([\fg, \fg])^\perp_\od$ such that $\mathscr{K}:=\Span\{x\}$ is an ideal. Then $(\fg,\omega_\fg)$ is obtained as an ${\mathscr D}_\ev$-extension from a~periplectic quasi-Frobenius Lie superalgebra $(\fa,\omega_\fa)$.
\end{Theorem}
\begin{proof}

Let us show that $\mathscr{K}^\perp$ is also an ideal in $(\fg,\omega_\fg)$. Indeed, let us show first that $(\mathscr{K}^\perp)_\od=\fg_\od$. This is true because $x$ is orthogonal to any odd element since $\omega_\fg$ is odd. Let $f\in \mathscr{K}^\perp$ and $g\in \fg$, then $[f,g]_\fg \in \mathscr{K}^\perp$ because
\[
\omega_\fg(x,[f,g]_\fg)=-(-1)^{p(g)p(f)} \omega_\fg(g, [x,f]_\fg) - \omega_\fg(f, [g,x]_\fg)=0-0=0.
\]
Since $\mathscr{K}$ is 1-dimensional and $\omega_\fg(x,x)=0$, it follows that $\mathscr{K}\subset \mathscr{K}^\perp$ and $\dim(\mathscr{K}^\perp)=\dim(\fg)-1$. Therefore, there exists $e \in \fg_\ev$ (since $(\mathscr{K}^\perp)_\od=\fg_\od$) such that
\[
\fg=\mathscr{K}^\perp\oplus \mathscr{K}^*, \quad \text{ where  $\mathscr{K}^*:=\Span\{x^*\}$}.
\]
We can normalize $e$ so that $\omega_\fg(e,x)=1$. 

Let us define $\fa:=(\mathscr{K} +\mathscr{K}^*)^\perp$. We then have a decomposition $\fg=\mathscr{K} \oplus \fa \oplus \mathscr{K}^*$. Now the proof follows as in Theorem \ref{Rec1}.\end{proof}
\subsection{${\mathscr D}_\od$-extension}
Let $(\fa, \omega_\fa)$ be a  symplectic Lie superalgebra. Let ${\mathscr D}\in \fder_\od(\fa)$ be a derivation such that the following conditions are satisfied: the map
\[
\begin{array}{l}
(a,b)\mapsto \Omega(a,b):= \omega_\fa( {\mathscr D} \circ {\mathscr D}(a)  - {\mathscr D}^* \circ {\mathscr D}^*(a)  , b) 
\end{array}
\]
is in $B^2(\fg;\Kee)$ \textup{(}i.e., a 2-coboundary\textup{)}. Let us write $\Omega= \delta (\chi)$ for some $\chi\in C^1(\fa; \Kee)$; namely, $\Omega(a,b)=\chi([a,b]_\fa)$. Since $\omega_\fa$ is non-degenerate, there exists $a_0\in \fa$ such that $\Omega(a,b)=\omega_\fa(a_0,[a,b]_\fa)$.

\sssbegin[(${\mathscr D}_\od$-extension -- the case where $\omega$ is periplectic)]{Theorem}[${\mathscr D}_\od$-extension -- the case where $\omega$ is periplectic]\label{MainThOO} Let $(\fa, \omega_\fa)$ be a~periplectic quasi-Frobenius Lie superalgebra. Let ${\mathscr D}\in \fder_\od(\fa)$ be a derivation and there exists $a_0\in \fa_\ev$ such that 
\[
{\mathscr D}^2=\ad_{a_0}, \quad  {\mathscr D}(a_0)={\mathscr D}^*(a_0)=0, \quad \text{and} \quad \Omega(a,b)=\omega_\fa(a_0, [a,b]_\fa).
\]
Then, there exists a Lie superalgebra structure on $\fg:=\mathscr{K} \oplus \fa \oplus \mathscr{K} ^*$, where $\mathscr{K} :=\Span\{x\}$ for $x$ even and $\mathscr{K}^* :=\Span\{e\}$ for $e$ odd, defined as follows: \textup{(}for any $a, b\in \fa$\textup{)}
\[
\begin{array}{ll}
[\mathscr{K} ,\fg]_\fg=0, & [a,b]_\fg  :=  [a,b]_\fa +  \left (\omega_\fa( {\mathscr D}(a), b) +(-1)^{p(a)}\omega_\fa(a, {\mathscr D}(b))\right )x, \\[2mm]
[e, e]_\fg=2 a_0, & [e,a]_\fg :=  {\mathscr D}(a) + \omega_\fa(a_0,a)x.
\end{array}
\]
There exists a~closed anti-symmetric periplectic form $\omega_\fg$ on $\fg$ defined as follows: 
\begin{eqnarray*}
{\omega_\fg}\vert_{\fa \times \fa}:= \omega_\fa, \quad \omega_\fg(\fa,\mathscr{K} ):=\omega_\fg(\fa,\mathscr{K} ^*):=0, \quad \omega_\fg(e,x):=1, \quad 
\omega_\fg(x,x):=\omega_\fg(e,e) :=0.
\end{eqnarray*}

\end{Theorem}
\begin{proof}

Similar to that of Thereom \ref{MainTh},\end{proof}

The converse of Theorem \ref{MainThOO} is given by the following Theorem.

\sssbegin[(Converse of Theorem \ref{MainThOO})]{Theorem}[Converse of Theorem \ref{MainThOO}]
\label{Rec4}
Let $(\fg,\omega_\fg)$ be a~periplectic quasi-Frobenius Lie superalgebra. Suppose that $\fz(\fg)_\ev\not=\{0\}.$ Then, $(\fg,\omega_\fg)$ is obtained as an ${\mathscr D}_\od$-extension of a~periplectic quasi-Frobenius Lie superalgebra $(\fa,\omega_\fa)$.

\end{Theorem}
\begin{proof}

Let $x$ be a non-zero element in $\fz(\fg)_\ev$. Since $x$ is central in $\fg$, it follows that the subspace $\mathscr{K}:=\Span\{x\}$ is an ideal in $(\fg, \omega_\fg)$. Moreover, $\mathscr{K}^\perp$ is also an ideal in $(\fg,\omega_\fg)$. Indeed, let us show first that $(\mathscr{K}^\perp)_\ev=\fg_\ev$. This is true because $x$ is orthogonal to any even  element since $\omega_\fg$ is odd. Let $f\in \mathscr{K}^\perp$ and $g\in \fg$, then $[f,g]_\fg \in \mathscr{K}^\perp$ because
\[
\omega_\fg(x,[f,g]_\fg)=-(-1)^{p(g)p(f)} \omega_\fg(g, [x,f]_\fg) - \omega_\fg(f, [g,x]_\fg)=0-0=0.
\]
Since  $\omega_\fg(x,x)=0$, it follows that $\mathscr{K}\subset \mathscr{K}^\perp$ and $\dim(\mathscr{K}^\perp)=\dim(\fg)-1$. Therefore, there exists $e \in \fg_\od$ (since $(\mathscr{K}^\perp)_\ev=\fg_\ev$) such that
\[
\fg=\mathscr{K}^\perp\oplus \mathscr{K}^*, \quad \text{ where  $\mathscr{K}^*:=\Span\{e\}$}.
\]
We can normalize $e$ so that $\omega_\fg(e,x)=1$. Besides, $\omega_\fg(x,x)=0$ since $\mathscr{K}=\mathscr{K}\cap \mathscr{K}^\perp$.

Let us define $\fa:=(\mathscr{K} +\mathscr{K}^*)^\perp$. We then have a decomposition $\fg=\mathscr{K} \oplus \fa \oplus \mathscr{K}^*$. 

Let us define a periplectic form on $\fa$ by setting: 
\[
\omega_\fa={\omega_\fg}\vert_{\fa\times \fa}.
\]  
The form $\omega_\fa$ is non-degenerate on $\fa$. Indeed, suppose there exists an $a\in\fa$ such that 
\[
\omega_\fa(a,f)=0\quad \text{ for any $f\in \fa$}.
\] 
But $a$ is also orthogonal to $x$ as well as to $e$. It follows that $\omega_\fg(a,f)=0$ for any $f\in\fg$. Hence, $a=0$, since $\omega_\fg$ is nondegenerate on $\fg$. 

Let us now show that vector space $\fa$ can be endowed with a Lie superalgebra structure, and there exists a~periplectic structure on $\fa$ for which $\fg$ is its double extension. Let $a,b \in \fa$. The bracket $[a,b]_\fg$ belongs to $\mathscr{K} \oplus \fa$ because $\fa \subset \mathscr{K}^\perp=\mathscr{K}\oplus \fa$ and the latter is an ideal. It follows that 
\[
[a,b]_\fg=\psi(a,b)+\phi(a,b)x, \quad [e,a]_\fg={\mathscr D}(a)+\chi(a)x,\quad [e,e]_\fg=2 a_0 + \nu x.
\]
where $\psi(a,b), {\mathscr D}(a) \in \fa, a_0\in \fa_\ev$, and $\phi(a,b), \chi (a) \in \Bbb K$. The proof  follows as in Theorem~\ref{MainThO}.\end{proof}
\section{The $T^*$-extensions  and $\Pi T^*$-extensions (Lagrangian extensions)}\label{TstarChap}
Following \cite{BC}, a {\it polarization} for a~quasi-Frobenius Lie superalgebra $(\fg, \omega)$ is a choice of a Lagrangian subalgebra $\fl$ of $(\fg, \omega)$. This terminology was motivated by the fact that  Lagrangian foliations in a manifold with a symplectic form is a polarization of the manifold, see \cite{BC}. A {\it strong polarization} of a~quasi-Frobenius Lie superalgebra $(\fg, \omega)$ is a pair $(\fa, N)$ consisting of a Lagrangian {\it ideal} $\fa\subset \fg$ and a complementary Lagrangian {\it subspace}  $N\subset \fg$. The quadruple $(\fg, \omega,\fa, N)$ is referred to as a {\it strongly polarized}  quasi-Frobenius Lie superalgebra.

An isomorphism of strongly polarized quasi-Frobenius Lie superalgebras $(\fg, \omega, \fa, N)\rightarrow (\fg', \omega', \fa', N')$ is an even isomorphism of quasi-Frobenius Lie superalgebras $(\fg, \omega) \rightarrow (\fg', \omega')$ wich maps the strong polarization $(\fa, N)$ to the strong polarization $(\fa',N')$.

Let $\fh$ be a Lie superalgebra and let $\nabla$ be a~flat torsion-free connection on $\fg$. We call $(\mathfrak{h}, \nabla)$ a~flat Lie superalgebra. 

We will introduce the notion of $T^*$- and $\Pi T^*$-extension of Lie superalgebras.  To our best knowledge,  the notion of $T^*$-extension for Lie algebras was introduced by Bordemann in \cite{Bor}. This notion was also studied in \cite{BC} in the context of quasi-Frobenius Lie algebras and called a {\it Lagrangian extension}.
\subsection{The construction of  $T^*$- and $\Pi T^*$-extensions (Lagrangian extensions)}\label{TstarSec}
Recall that since the connection $\nabla$ is flat, the map 
\[
\fh \rightarrow \End(\fh) \quad u\mapsto \nabla_u,
\]
defines a representation. We define a map $\rho: \fh \rightarrow \End(\fh^*)$, the dual representation, defined as follows:
\begin{equation}\label{eq5.1}
\rho: \fh \rightarrow \End(\fh^*) \quad u \mapsto \rho(u), \text{ where } \rho(u)\cdot \xi =- (-1)^{p(u) p(\xi)} \xi \circ \nabla_u\quad \text{ for all $u\in \fh$ and $\xi\in \fh^*$}.
\end{equation} 
Similarly, we have a representation 
\begin{equation}\label{eq5.1b}
\chi: \fh \rightarrow \End(\Pi(\fh^*)) \quad u\mapsto \chi(u):= (-1)^{p(u)}\, \Pi \circ \rho(u)\circ \Pi.
\end{equation} 

\sssbegin[(Dual representation)]{Lemma}[Dual representation]
The maps \eqref{eq5.1} and \eqref{eq5.1b} are indeed representations. 

\end{Lemma}
\begin{proof}

First, let us compute
\[
\rho(u)\circ \rho(v) (\xi)= \rho(u) (- (-1)^{p(\xi) p(v)} \xi \circ \nabla_v)= (-1)^{p(\xi) p(v)+p(u) (p(v)+p(\xi))} \xi \circ \nabla_v \circ \nabla_u
\]
Now, using the fact that the connection is flat we get (for all $u,v\in \fh$ and for all $\xi\in \fh^*$):
\begin{align*}
&\rho(u)\circ \rho(v) (\xi)- (-1)^{p(u)p(v)} \rho(v)\circ \rho(u) (\xi)- \rho([u,v])(\xi) \\
&= (-1)^{p(\xi) p(v)+p(u) (p(v)+p(\xi))} \xi \circ \nabla_v \circ \nabla_u -  (-1)^{p(\xi) p(u)+p(v) p(\xi)} \xi \circ \nabla_u \circ \nabla_v- (-1)^{p(\xi)(p(u)+p(v))} \xi \circ \nabla_{[u,v]}\\
& =(-1)^{p(\xi)(p(u)+p(v))} \xi \circ \left (\nabla_u \circ \nabla_v -(-1)^{p(u)p(v)} \nabla_v \circ \nabla_u- [u,v]  \right )=0.
\end{align*}
On the other hand, 
\begin{align*}
&\chi(u)\circ \chi(v) (\Pi(\xi))- (-1)^{p(u)p(v)} \chi(v)\circ \chi(u) (\Pi(\xi))- \chi([u,v])(\Pi(\xi)) \\
&=(-1)^{p(u)+p(v)}\Pi \circ (\rho(u)\circ \rho(v) (\xi)- (-1)^{p(u)p(v)} \rho(v)\circ \rho(u) (\xi)- \rho([u,v])(\xi) )=0.\qed
\end{align*}
\noqed
\end{proof}
By means of a 2-cocycle $\alpha \in Z^2(\fh, \fh^*)$ (resp. $\beta  \in Z^2(\fh, \Pi(\fh^*))$ we will construct a Lie superalgebra structure on $\fg:=\fh\oplus \fh^*$ (resp. $\fg:=\fh\oplus \Pi(\fh^*)$), called $T^*$-extension (resp.  $\Pi T^*$-extension).  Recall that the space $\fh^*$ and $\Pi(\fh^*)$ are $\fh$-modules by means of the representations  (\ref{eq5.1}) and (\ref{eq5.1b}).\\

\underline{On the space $\fg:=\fh\oplus \fh^*$:} The brackets are defined as follows: \textup{(}for all $u, v\in \fh$ and $\xi\in \fh^*$\textup{)}
\[
\begin{array}{l}
 [u,v]_\fg  :=  [u,v]_\fh +  \alpha(u,v), \quad [u,\xi]_\fg :=  \rho(u) \cdot \xi.
\end{array}
\]
We define an {\bf even} form $\omega$ on $\fg$ as follows: \textup{(}for all $u+\xi, v+\zeta \in \fh\oplus \fh^*$\textup{)}
\begin{equation}\label{omegaTstar}
\omega(u+\xi, v+ \zeta)= \xi(v) - (-1)^{p(\zeta) p(u)} \zeta(u).
\end{equation}
This construction will be referred to as the {\it $T^*$-extension}, and denoted by $(\fh\oplus \fh^*,\alpha, \nabla, \omega)$. 

\underline{On the space $\fk:=\fh\oplus \Pi(\fh^*)$:} The brackets are defined as follows: \textup{(}for any $u, v\in \fa$\textup{)}
\[
\begin{array}{l}
 [u,v]_\fg  :=  [u,v]_\fh +  \beta(u,v), \quad [u, \Pi(\xi)]_\fg :=  \chi (u) \cdot \Pi(\xi).
\end{array}
\]
We define an {\bf odd} form $\kappa$ on $\fg$ as follows:  \textup{(}for all $u+\Pi(\xi), v+\Pi(\zeta) \in \fh\oplus \fh^*$\textup{)}
\begin{equation}\label{omegaTstarodd}
\kappa(u+\Pi(\xi), v+ \Pi(\zeta))= \xi(v) - (-1)^{(p(\zeta)+1) p(u)} \zeta (u).
\end{equation}
This construction will be referred to as the {\it $\Pi T^*$-extension}, and denoted by $(\fh\oplus \Pi(\fh^*), \beta, \nabla, \kappa)$.  

\sssbegin[(Conditions on $\alpha$ and $\beta$)]{Proposition}[Conditions on $\alpha$ and $\beta$] \label{condalphabeta}\textup{(}i\textup{)} The form $\omega$ is closed on $\fg$ if and only if 
\begin{equation}\label{eq5.4}
(-1)^{p(u)p(w)}\alpha(u,v)(w)+ \circlearrowleft (u,v,w)=0 \text{ for all $u,v,w\in\fh$}.
\end{equation}
\textup{(}ii\textup{)} The form $\kappa$ is closed on $\fg$ if and only if 
\begin{equation}\label{eq5.4b}
(-1)^{p(u)p(w)}\Pi \circ \beta(u,v)(w)+ \circlearrowleft (u,v,w)=0 \text{ for all $u,v,w\in\fh$}.
\end{equation}

\end{Proposition}
\begin{proof} 

Let us prove Part (i). It is enough to check the cocycle condition of the form $\omega_\fg$ pertaining to the triple $(u,v+\zeta, w)$. First, we compute
\[
\begin{array}{l}
\omega ([u, v+\zeta]_\fg,w) = \omega_\fg ([u, v]_\fh+\alpha(u,v) +\rho(u)\circ \zeta,w)= \alpha(u,v)(w) + \rho(u)\circ \zeta(w) \\[2mm]
=\alpha(u,v)(w) -(-1)^{p(u)p(\zeta)} \zeta \circ \nabla_u(w).
\end{array}
\]
Similarly, 
\[
\begin{array}{lcl}
 \omega([w, u]_\fg, v+\zeta)&=&\alpha(w,u)(v) - (-1)^{p(\zeta)(p(w)+p(u))} \zeta ([w,u]_\fh)\\[2mm]
\omega ([v+\zeta, w]_\fg, u)&=&\alpha(v,w)(u) + \zeta \circ \nabla_w (u).
\end{array}
\]
Because $\omega_\fg$ is bilinear, we may assume that $p(v)=p(\zeta)$. 

In order to check the condition \[
(-1)^{p(u)p(w)}\omega_\fg ([u, v+\zeta]_\fg,w) + \circlearrowleft (u,v+\zeta,w)=0,
\] 
we see that the part in $\fh$ is equivalent to Eq. (\ref{eq5.4}). The part in $\fh^*$ contains the terms involving  $\zeta$. Let us collect them. Using the fact that the connection is flat, we have
\[
\begin{array}{l}
(-1)^{p(u)p(w)}\left ( -(-1)^{p(u)p(\zeta)} \zeta \circ \nabla_u(w)\right ) + (-1)^{p(w)p(v)}\left (  - (-1)^{p(\zeta)(p(w)+p(u))} \zeta ([w,u]_\fh) \right ) \\[2mm]
+(-1)^{p(v)p(u)}\left (  \zeta \circ \nabla_w (u) \right )=(-1)^{p(u)p(v)} \zeta \circ \left (\nabla_w (u) -(-1)^{p(u)p(w)} \nabla_u (w)  -[w,u]_\fh \right )=0.
\end{array}
\]
The proof of Part (ii) is similar. 
\end{proof}
As a result, we arrive at the following theorem.

\sssbegin [(Lagrangian or $T^*$- and $\Pi T^*$-extensions)]{Theorem} [Lagrangian or $T^*$- and $\Pi T^*$-extensions] \label{Tstar}
Let $(\fh, \nabla)$ be a flat Lie superalgebra. To every 2-cocycle $\alpha\in Z^2(\fh, \fh^*)$ \textup{(}resp. $\beta\in Z^2(\fh, \Pi(\fh^*))$\textup{)} which satisfies \eqref{eq5.4}  \textup{(}resp. \eqref{eq5.4b}\textup{)} one can canonically associate a strongly polarized quasi-Frobenius Lie superalgebra $\textup{(}\fg, \omega,\fh, \fh^*)$ \textup{(}resp. $\textup{(}\fk, \kappa,\fh, \Pi(\fh^*))$\textup{)}, where $\fg=\fh\oplus \fh^*$ \textup{(}resp. $\fk=\fh\oplus \Pi(\fh^*)$\textup{)} and $\omega$ \textup{(}resp. $\kappa$\textup{)} is given as in \eqref{omegaTstar} \textup{(}resp.~\eqref{omegaTstarodd}\textup{)}.
\end{Theorem}

We will now prove the converse of this theorem, namely every orthosymplectic or periplectic quasi-Frobenius Lie superalgebra $(\fg, \omega)$ that has a Lagrangian ideal $\fa$ can be obtained as a $T^*$-extension  or $\Pi T^*$-extension of a flat Lie superalgebra. 

Let $(\fg, \omega, \fa,N)$ be a strongly polarized quasi-Frobenius Lie superalgebra. Here the form $\omega$ is either orthosymplectic or periplectic. Consider the quotient space $\fh=\fg/\fa$ (recall that $\fa^\perp=\fa$). We will construct a flat connection $\nabla$ on $\fh$, superizing the construction of \cite{BC}. For every $u,v\in \fh$, we denote by $\tilde u$ and $\tilde v$ their lift to $\fg$. We then write 
\begin{equation}\label{eq5.5}
\omega_\fh(\nabla_uv, a)= - (-1)^{p(u)p(v)} \omega_\fg (\tilde v, [\tilde u, a]).
\end{equation}
The pair $(\fh,\nabla)$ is called {\it the quotient flat Lie superalgebra associated with the Lagrangian ideal} $\fa$ of $(\fg, \omega)$.

\sssbegin[(Converse of Theorem \ref{Tstar})]{Theorem} [Converse of Theorem \ref{Tstar}] \label{TstarConv} Let $(\fg, \omega, \fa,N)$ be a strongly polarized ortho-symplectic or periplectic quasi-Frobenius Lie superalgebra and $(\fh, \nabla)$ be its associated quotient flat Lie superalgebra. 

\textup{(}i\textup{)} If the form $\omega$ is orthosymplectic, then there exists $\alpha \in Z^2(\fh, \fh^*)$ 
satisfying Eq. \eqref{eq5.4}, such that $(\fg, \omega, \fa,N)$ is isomorphic to the $T^*$-extension of $(\fh,\nabla,\alpha)$.

\textup{(}ii\textup{)} If the form $\omega$ is periplectic, then there exists $\beta \in Z^2(\fh, \Pi(\fh^*))$ 
satisfying Eq. \eqref{eq5.4b}, such that $(\fg, \omega, \fa,N)$ is isomorphic to the $\Pi T^*$-extension of $(\fh,\nabla,\beta)$.

\end{Theorem}

Let's fix some notations before presenting the proof. We keep the notation $\omega_\fg$ when the form is orthosymplectic. The notation is changed to $\kappa_\fg$ when the form is periplectic. 

\begin{proof} 

First, observe that the representation 
\[
\ad_{\fh,\fa}: \fh \rightarrow \End(\fa, \fa) \quad u \mapsto [\tilde{u},\cdot ], \quad \text{ where $\tilde u$ is a lift of $u$ to $\fg$,}
\]
is well-defined, since $\fa$ is abelian.

Let $\pi_\fa :  \fg \rightarrow \fa$ and  $\pi_N :  \fg \rightarrow N$ be the projection maps induced by the strong polarization $\fg = \fa \oplus N$. For every $u, v \in  \fh$, let $\tilde u,  \tilde v \in \fg$ denote their lifts to $\fg$ with respect to the homomorphism $\fg \rightarrow  \fh$. We will always assume that this {\bf lifting preserves the parity}; namely,  $p(u)=p(\tilde u)$ for all $u\in \fh$. 

The expression 
\begin{equation}\label{alphapia}
\tilde \alpha (u, v) = \pi_a([\pi_N(\tilde u), \pi_N(\tilde v)]_\fg) \in \fa
\end{equation}
 is well-defined since the expressions $\pi_N(\tilde u)$ and $\pi_N(\tilde v)$ are unique. We will show that it defines an even 2-cocycle $\tilde \alpha \in  Z^2_\ev(\fh, \fa)$ for the representation $\ad_{\fh,\fa}$. We need the following lemma. 

\sssbegin[(Technical)]{Lemma}[Technical] \label{Lemtec}
\textup{(}i\textup{)} For all $u,v\in \fh$, we have \textup{(}where $\tilde x$ is a lift of $x$ to $\fg$\textup{)}:
\[ \pi_N(\widetilde{[u,v]_\fh})= \pi_N([\tilde u, \tilde v]_\fg) = \pi_N([\pi_N(\tilde u), \pi_N(\tilde v)]_\fg);\]

\textup{(}ii\textup{)} For all $u,v,w\in \fh$, we have \textup{(}where $\tilde x$ is the lift of $x$ to $\fg$\textup{)}:
\[u\cdot \tilde \alpha(v,w)+\tilde \alpha(u,[v,w])=  \pi_a([\pi_N(\tilde u), [\pi_N(\tilde v), \pi_N(\tilde w)]_\fg]_\fg).\]

\end{Lemma}
\begin{proof}

Let us just prove the second equality of part (i). Using the fact that $\fa$ is an ideal, we have
\[
 \pi_N([\tilde u, \tilde v]_\fg)=  \pi_N([\pi_N(\tilde u)+\pi_\fa(\tilde u) , \pi_N(\tilde v)+\pi_\fa(\tilde v)]_\fg) =  \pi_N([\pi_N(\tilde u), \pi_N(\tilde v)]_\fg).
\]
Let us now prove Part (ii). Using Part (i), we get
\begin{align*}
&u \cdot \tilde \alpha (v,w)+\tilde \alpha(u, [v,w]_\fh)=[\tilde u , \pi_\fa([\pi_N(\tilde v), \pi_N(\tilde w)]_\fg)]_\fg+\pi_\fa([\pi_N(\tilde u),  \pi_N(\widetilde{[v, w]_\fh)}]_\fg)\\
&=[\tilde u , \pi_\fa([\pi_N(\tilde v), \pi_N(\tilde w)]_\fg)]_\fg+\pi_\fa([\pi_N(\tilde u),  \pi_N([\pi_N(\tilde v), \pi_N(\tilde w)]_\fg)]_\fg)\\
&=\pi_\fa([\tilde u , \pi_\fa([\pi_N(\tilde v), \pi_N(\tilde w)]_\fg)]_\fg)+\pi_\fa([\pi_N(\tilde u),  \pi_N([\pi_N(\tilde v), \pi_N(\tilde w)]_\fg)]_\fg)\\
&=\pi_\fa([ \pi_N(\tilde u)+ \pi_\fa(\tilde u)  , \pi_\fa([\pi_N(\tilde v), \pi_N(\tilde w)]_\fg)]_\fg)+\pi_\fa([\pi_N(\tilde u),  \pi_N([\pi_N(\tilde v), \pi_N(\tilde w)]_\fg)]_\fg)\\
&=\pi_\fa([\pi_N(\tilde u),  [\pi_N(\tilde v), \pi_N(\tilde w)]_\fg]_\fg).\qed
\end{align*}
\noqed
\end{proof}

\sssbegin[($\tilde \alpha$ is a 2-cocycle)]{Lemma}[$\tilde \alpha$ is a 2-cocycle] We have $\tilde \alpha \in  Z^2_\ev(\fh, \fa)$.

\end{Lemma}
\begin{proof}

Using Lemma \ref{Lemtec}, and since the Jacobi identity holds on $\fg$, we have
\begin{align*}
& u \cdot \alpha(v,w)-(-1)^{p(u)p(v)} v \cdot \alpha(u,w)+ (-1)^{p(w)(p(u)+p(v))} w \cdot \alpha(u,v)- \alpha([u,v]_\fh, w)\\
& +(-1)^{p(v)p(w)}\alpha([u,w]_\fh, v)+\alpha(u, [v,w]_\fh) \\ 
&=\pi_a([\pi_N(\tilde u), [\pi_N(\tilde v), \pi_N(\tilde w)]_\fg]_\fg) - (-1)^{p(u)p(v)}\pi_a([\pi_N(\tilde v), [\pi_N(\tilde u), \pi_N(\tilde w)]_\fg]_\fg)\\
&+(-1)^{p(w)(p(u)+p(v))}\pi_a([\pi_N(\tilde w), [\pi_N(\tilde u), \pi_N(\tilde v)]_\fg]_\fg)=0.\qed
\end{align*}
\noqed
\end{proof}

Let $i_\omega: \fa \rightarrow \fh^*,  \; a \mapsto  \omega(a,\cdot)$ (resp. $i_\kappa: \fa \rightarrow \fh^*,  \; a \mapsto  \kappa(a,\cdot)$) be the identification of $\fa$ with $\fh^*$, which is
induced by $\omega$ (resp. $\kappa$). Namely, for every $u\in \fh$ and $a\in \fa$, we put $i_\omega(a)(u):= \omega(a, \tilde u)$ (resp. $i_\kappa(a)(u):= \kappa(a, \tilde u)$), where $\tilde u$ is a lift of $u$ to $\fg$ such that $p(u)=p(\tilde u)$. These two maps are well-defined because $\fa$ is Lagrangian. Besides, $p(i_\omega)=\ev$ (resp. $p(i_\kappa)=\od$). 

We need the following Lemma. 

\sssbegin[(Equivalence of $\ad_{h,a}$ and $ \rho$ (resp. $\chi$))]{Lemma}[Equivalence of $\ad_{h,a}$ and $ \rho$ (resp. $\chi$)] \label{actionequiv} With the definitions given in Eqs. \eqref{eq5.1} and \eqref{eq5.1b}, we have 
\begin{eqnarray} \label{eq5.8}
 \rho(u) \circ  i_\omega &  = & i_\omega \circ  \ad_{\fh,\fa}(u) \quad \text{ \textup{(}if the form $\omega$ is orthosymplectic\textup{)}}.\\
\label{eq5.8b}  \chi(u) \circ  \Pi\circ i_\kappa &  = & \Pi\circ  i_\kappa \circ  \ad_{\fh,\fa}(u) \quad \text{ \textup{(}if the form $\kappa$ is periplectic\textup{)}}.
\end{eqnarray}

\end{Lemma}
\begin{proof}

Eq. (\ref{eq5.5}) is actually equivalent to Eq. (\ref{eq5.8}). Let us only prove Eq. (\ref{eq5.8b}). Let us first check the equality 
\begin{equation} \label{rhodash}
\rho(u) \circ i_\kappa(a)=(-1)^{p(u)}i_\kappa(a) \circ \ad_{\fh, \fa}(u). 
\end{equation} 
Indeed,  for every $u,v\in \fh$ and for all $a\in \fa$, we have 
\begin{align*}
&\rho(u) \circ i_\kappa(a)(v)= -(-1)^{p(u)p(i_\kappa(a))} i_\kappa(a) \circ \nabla_u(v)=  -(-1)^{p(u)(p(a)+1)} i_\kappa(a, \nabla_u(v))\\
&=(-1)^{p(u)(p(a)+1)+p(u)p(a)} \kappa([\tilde u, a], \tilde v)=(-1)^{p(u)} i_\kappa  \circ \ad_{a, \fh} (u)(a)(v).
\end{align*}
Now, for every $u\in \fh$, and using Eq. (\ref{rhodash}), we have
\[
 \chi(u) \circ  \Pi\circ i_\kappa (a) = (-1)^{p(u)} \Pi \circ \rho(u) \circ i_\kappa(a)=  \Pi \circ i_\kappa \circ \ad_{\fh, \fa}(u)(a).\qed
\]
\noqed
\end{proof}

In particular, this shows that the representation $\rho$ of $\fh$ on $\fh^*$ (resp. $\chi$ of $\fh$ on $\Pi(\fh^*)$), which belongs to the flat connection $\nabla$ by (\ref{eq5.1}) (resp. by (\ref{eq5.1b})) is equivalent to $\ad_{\fh,\fa}.$ \\

Let $\pi_\fh:N \rightarrow \fh$ be the isomorphism of vector spaces induced by the projection map $\fg \rightarrow \fh.$ We have two cases.

\underline{The case where the form is orthosymplectic.} The isomorphisms $\pi_\fh$ and $i_\omega$ assemble to an iso-morphism
\[
\Phi:=\pi_\fh \oplus i_\omega:\fg = N \oplus \fa \rightarrow \fh \oplus \fh^*.
\]
Define 
\begin{equation}\label{eq5.9}
\alpha := i_\omega \circ \tilde \alpha \in Z^2_\ev(\fh, \fh^*).
\end{equation}
to be the {\it push-forward} of $\tilde \alpha$. Lemma \ref{actionequiv} shows that $\alpha \in Z^2_\ev(\fh, \fh^*)$. It is now easily verified that the map $\pi_\fh \oplus i_\omega: (\fg,\omega_\fg) \rightarrow (\fh\oplus \fh^*,\alpha,\nabla,\omega_{\fh\oplus \fh^*})$ defines an isomorphism of Lie superalgebras preserving the orthosymplectic forms. As a consequence, $\alpha$ satisfies (\ref{eq5.4}) by Proposition \ref{condalphabeta}.

In fact, (where $x=x_a+x_N$ and $y=y_a+y_N$):
\begin{align*}
&(\pi_\fh\oplus i_\omega)^* \omega_{\fh\oplus \fh^*}(x,y)=  \omega_{\fh\oplus \fh^*}(\pi_\fh(x_N)+i_\omega(x_a), \pi_\fh(y_N)+i_\omega(y_a))=\\[2mm]
 &i_\omega(x_a)(\pi_\fh(y_N))-(-1)^{Pty} i_\omega(y_a)(\pi_\fh(x_N))=\omega_\fh(x_a, \pi_\fh(y_N))-(-1)^{Pty}  \omega_\fh(y_a, \pi_\fh(x_N))\\[2mm]
&=\omega_\fg(x_a, \tilde  \pi_\fh(y_N))-(-1)^{Pty}  \omega_\fg(y_a,  \tilde \pi_\fh(x_N))= \omega_\fg(x_a, y_N)-(-1)^{Pty}  \omega_\fg(y_a,  x_N) \\[2mm]
&= \omega_\fg(x_a, y_N) +   \omega_\fg(x_N, y_a) = \omega_\fg(x_a, y_a+y_N) +   \omega_\fg(x_N, y_N+y_a) =\omega_\fg(x,y).
\end{align*}
Let us check that $\Phi$ preserves the bracket. Let us start with $x\in N$ and $y\in \fa$. Then, 
\[
\Phi([x,y])= i_\omega ([x,y]) = \omega_\fg([x,y], \cdot).
\]
On the other hand, 
\begin{align*}
&[\Phi(x), \Phi(y)]=[\pi_\fh(x) , i_\omega(y,\cdot)]=-(-1)^{p(x)p(y)} i_\omega(y, \cdot) \circ \nabla_{\pi_\fh(x)}=-(-1)^{p(x)p(y)} \omega_\fg(y,  \nabla_{\pi_\fh(x)} (\cdot))\\
&=-(-1)^{p(x)p(y)} \omega_\fg([y, \pi_\fh(x)], (\cdot))=-(-1)^{p(x)p(y)} \omega_\fg([y, x], (\cdot))= \omega_\fg([x, y], (\cdot)).
\end{align*}
Let us consider now $x\in N$ and $y\in N$. We have
\[
\Phi([x,y])= \Phi(\pi_N([x,y]) + \pi_\fa ([x,y]))= \pi_\fh(\pi_N([x,y]))+ i_\omega (\pi_\fa([x,y]))=\pi_\fh(\pi_N([x,y]))+ \alpha (x,y).
\]
On the other hand, 
\[
[\Phi(x), \Phi(y)]= [\pi_\fh(x), \pi_\fh(y)]+ \alpha(x,y)= \pi_\fh([x,y])+ \alpha(x,y)=\pi_\fh(\pi_N([x,y]))+ \alpha(x,y).
\]
Therefore, the equality $[\Phi(x), \Phi(y)]=\Phi([x,y])$ holds for every $x,y\in \fg$.\\

\underline{The case where the form is periplectic.} The isomorphisms $\pi_\fh$ and $\Pi \circ i_\kappa$ assemble to an iso-morphism
$\pi_\fh \oplus \Pi\circ i_\kappa:\fg = N \oplus \fa \rightarrow \fh \oplus \Pi(\fh^*).$ 
Define 
\begin{equation}\label{eq5.9b}
\beta := \Pi\circ i_\kappa \circ \tilde \alpha  \in Z^2_\ev(\fh, \Pi(\fh^*)).
\end{equation}
to be the push-forward of $\tilde \beta$. The map $\pi_\fh \oplus \Pi \circ i_\kappa: (\fg,\kappa_\fg) \rightarrow (\fh\oplus \Pi(\fh^*),\beta, \nabla,\kappa_{\fh\oplus \fh^*})$ defines an isomorphism of Lie superalgebras preserving periplectic forms. As a consequence, $\beta$  satisfies (\ref{eq5.4b}) by Prop. \ref{condalphabeta}.

In fact (where $x=x_a+x_N$ and $y=y_a+y_N$):
\begin{align*}
&(\pi_\fh\oplus \Pi\circ i_\kappa)^* \kappa_{\fh\oplus \Pi(\fh^*)}(x,y)=  \kappa_{\fh\oplus \fh^*}(\pi_\fh(x_N)+\Pi \circ i_\kappa(x_a), \pi_\fh(y_N)+\Pi \circ i_\kappa (y_a))\\
&= i_\kappa(x_a)(\pi_\fh(y_N))-(-1)^{Pty} i_\kappa(y_a)(\pi_\fh(x_N))=\kappa_\fh(x_a, \pi_\fh(y_N))-(-1)^{Pty}  \kappa_\fh(y_a, \pi_\fh(x_N))\\
&=\kappa_\fg(x_a, \tilde  \pi_\fh(y_N))-(-1)^{Pty}  \kappa_\fg(y_a,  \tilde \pi_\fh(x_N))= \kappa_\fg(x_a, y_N)-(-1)^{Pty}  \kappa_\fg(y_a,  x_N) \\
&= \kappa_\fg(x_a, y_N) +   \kappa_\fg(x_N, y_a) = \kappa_\fg(x_a, y_a+y_N) +   \kappa_\fg(x_N, y_N+y_a) =\kappa_\fg(x,y).
\end{align*}
We skip the proof that the isomorphism preserves the product.
\end{proof}
We call the triple $(\fh, \nabla, \alpha)$ (resp.  $(\fh, \nabla, \beta)$) constructed above the `even' (resp. `odd') extension triple associated to $(\fg, \omega, \fa, N)$, where $\omega$ is orthosymplectic (resp. periplectic). In general, triples $(\fh, \nabla, \alpha)$ (resp. $(\fh, \nabla, \beta)$), where $\alpha$ (resp. $\beta$) satisfies the orthosymplectic (resp. periplectic) extension cocycle  \eqref{eq5.4} (resp. \eqref{eq5.4b}), will be called {\it an `even' \textup{(}resp. `odd'\textup{)} flat Lie superalgebra with orthosymplectic \textup{(}resp. periplectic\textup{)} extension cocycle}.
\subsection{Functoriality of the correspondence}
Following \cite{BC}, we consider briefly the functorial properties of our constructions. 

\sssbegin[(Push forward of a connection)]{Lemma} [Push forward of a connection]\label{Funct1} Let $(\fg,\omega)$ be an~orthosympletic or periplectic quasi-Frobenius Lie superalgebra. Let $\fa$ be Lagrangian ideal of $\fg$ and $(\fh = \fg/\fa,\nabla)$ be the associated flat Lie superalgebra. Let $\Phi:(\fg,\omega) \rightarrow
(\fg',\omega')$ be an isomorphism of Lie superalgebras preserving the forms, and let $\Phi_\fh:\fh = \fg/\fa \rightarrow \fh' = \fg'/\fa'$ be the induced map on the quotient spaces, where $\fa'=\Phi(\fa)$. Then, $\nabla' = (\Phi_\fh)_*$ \textup{(}push-forward of $\nabla$\textup{)} is the associated flat connection on $\fh'$.

\end{Lemma}
\begin{proof} 

An immediate consequence of Eq.(\ref{eq5.5}). \end{proof}

Similarly, we can state the following lemma.
\sssbegin[(The preimage of $\alpha$ and $\beta$)]{Lemma}[The preimage of $\alpha$ and $\beta$] \label{Funct2} Let $\Phi:(\fg,\omega, \fa,N) \rightarrow (\fg', \omega', \fa', N')$ be an isomorphism of strongly polarized quasi-Frobenius Lie superalgebras. Then 
\[
\begin{array}{lcl}
\alpha_{(\fg' ,\omega', \fa', N')} &= & (\Phi_\fh)_* \alpha_{(\fg, \omega, \fa, N)} \quad \text{ if the form $\omega$ is orthosymplectic.}\\[2mm]
\beta_{(\fg' ,\kappa', \fa', N')} &= & (\Phi_\fh)_* \beta_{(\fg, \kappa, \fa, N)} \quad \text{ if the form $\kappa$ is periplectic.}
\end{array}
\]

\end{Lemma}
\begin{proof}

We will only give the prove in the case where the form is periplectic. Put $\beta = \beta_{(\fg,\kappa,\fa,N)}$ and $\beta' = \beta_{(\fg',\kappa',\fa',N')}$. Using repeatedly that $\Psi$ is an isomorphism of strongly polarized Lie superalgebras algebras and by equations (\ref{alphapia}) and (\ref{eq5.9b}), we have (for all $u,v,w\in \fh$):
\[
\begin{array}{lcl}
\beta'(\Phi_\fh u, \Phi_\fh v)(\Phi_\fh w) & = & \Pi\circ \kappa' (\pi_{\fa'} ([\pi_{N'}(\widetilde{\Phi_\fh u}), \pi_{N'}(\widetilde{\Phi_\fh v})]', \widetilde{\Phi_\fh w}) \\[2mm]
& = &\Pi \circ  \kappa'(\pi_{\fa'} ([\pi_{N'}(\Phi \tilde u),  \pi_{N'}(\Phi \tilde v)]', \Phi(\tilde  w)) \\[2mm]
& = & \Pi\circ \kappa'(\pi_{\fa'} (\Phi [\pi_N(\tilde u), \pi_N(\tilde v)], \Phi (\tilde w)) \\[2mm]
& = & \Pi \circ \kappa'(\Phi ( \pi_{\fa} ([ \pi_N(\tilde u), \pi_N(\tilde v)])), \Phi (\tilde w)) \\[2mm]
& =& \Pi \circ \kappa  (\pi_\fa([\pi_N(\tilde u), \pi_N(\tilde v)]), \tilde w)\\[2mm]
& = & \beta(u,v)(w).\qed
\end{array}
\]
\noqed\end{proof}
Let $(\fh, \nabla,\alpha)$ and $(\fh', \nabla',\alpha')$ be two flat Lie superalgebras with extension cocycles $\alpha$ and $\alpha'$ such that $p(\alpha)=p(\alpha')$.
An isomorphism of Lie superalgebras $\phi:\fh \rightarrow \fh'$ is an {\it isomorphism of flat Lie superalgebras
with extension cocycle} if $\phi^*\nabla'=\nabla$ and $\phi^* \alpha'=\alpha$.

By virtue of Lemma \ref{Funct1} and Lemma \ref{Funct2}, isomorphic strongly polarized Lie superalgebras give rise to isomorphic flat Lie superalgebras with extension
cocycle. Together with Theorem \ref{Tstar}, we therefore have the following corollary.

\sssbegin[(Isomorphism classes of strongly polarized Lie superalgebras vs isomorphism classes of flat Lie superalgebras with a cocycle)]{Corollary}[Isomorphism classes of strongly polarized Lie superalgebras vs isomorphism classes of flat Lie superalgebras with a cocycle] The correspondence which associates a strongly polarized Lie superalgebra $(\fg,\omega, \fa,N)$, where $\omega$ is orthosymplectic \textup{(}resp. periplectic\textup{)},  to its extension triple $(\fh,\nabla,\alpha)$ \textup{(}resp. $(\fh,\nabla,\beta)$\textup{)} induces a bijection between isomorphism classes of strongly polarized Lie superalgebras and isomorphism classes of `even' \textup{(}resp. `odd'\textup{)} flat Lie superalgebras with  orthosymplectic \textup{(}resp. periplectic\textup{)} extension cocycle.
\end{Corollary}
\subsection{Change of strong polarization }
Let $(\fg,\omega, \fa,N)$ and $(\fg,\omega, \fa,N')$ be two strong polariza-tions belonging to the same Lagrangian extension, where the form $\omega$ is either orthosymplectic or periplectic. We have showed in Theorem \ref{TstarConv} that each of them can be obtained as a $T^*$-extension or a $\Pi T^*$-extension by means of a 2-cocycle. We will show that these 2-cocycles differ by a coboundary. The following Lemma is a superization of a result in \cite{BC}.

\sssbegin[($\alpha$ and $\alpha'$ (resp. $\beta$ and $\beta'$) are cohomologous)]{Lemma}[$\alpha$ and $\alpha'$ (resp. $\beta$ and $\beta'$) are cohomologous]\label{classindN}
\textup{(}i\textup{)} If the form $\omega$ is orthosymplectic, then there exists $\sigma \in \Hom(\fh, \fh^*)$, satisfying
\begin{equation}\label{sigmaeven}
\sigma(u)(v) - (-1)^{p(u)p(v)}\sigma(v)(u) = 0 ,\quad \text{for all $u, v \in \fh$,} 
\end{equation}
such that $\alpha' =  \alpha + \partial_\rho \sigma$.

\textup{(}ii\textup{)} If the form $\omega$ is periplectic, the there exists $\mu \in \Hom(\fh, \Pi(\fh^*))$, satisfying
\begin{equation}\label{sigmaodd}
\mu(u)(v) - (-1)^{p(u)p(v)} \mu(v)(u) = 0 , \quad \text{for all $u, v \in \fh$,} 
\end{equation}
such that $\beta' = \beta + \partial_\chi \mu$.

\end{Lemma}
\begin{proof} 

We will denote the form by $\omega$ (resp. $\kappa$) if it is orthosymplectic (resp. periplectic) on $\fg$.

As $\fg=\fa\oplus N=\fa'\oplus N'$, let $\pi_\fa, \pi_N$, as well as $\pi_\fa, \pi_{N'}'$ be the corresponding projection 
on the summand of $\fg$. Let us put $\tau:= \pi_\fa-\pi'_\fa  \in \Hom(\fg, \fa)$. Clearly, $\tau$ is an~even map with $\fa \subseteq \text{Ker} (\tau)$. Moreover, since any element $x \in \fg$ can be uniquely expressed as $x=\pi_a(x)+\pi_N(x)=\pi_a'(x)+\pi_{N'}'(x)$, it follows that $\pi'_{N'}=\pi_N+\tau$. 

Since, all $N, N'$, and $\fa$ are Lagrangian, the homomorphism $\tau$ satisfies the condition (for all $n,m \in N$)
\begin{eqnarray}
\label{omegataueven} \omega(\tau(n),m) + \omega(n, \tau(m)) & = &  0  \quad \text{(if the form $\omega$ is orthosymplectic)}\\[2mm]
\label{omegatauodd}  \kappa(\tau(n),m) + \kappa(n, \tau(m)) & = &  0  \quad \text{(if the form $\kappa$ is periplectic).}
\end{eqnarray}
Let $u, v \in  \fh$  and $\tilde u,  \tilde v $ their respective lifts to $\fg$. We first compute the bracket, and using the fact that $\fa$ is abelian, we get
\begin{align*}
& [\pi'_{N'}(\pi_N(\tilde u)), \pi'_{N'}(\pi_N(\tilde v))] = [\pi_N(\tilde u)+ \tau (\pi_N(\tilde u)), \pi_N(\tilde v)+ \tau (\pi_N(\tilde v))] \\[2mm]
&=  [\pi_N(\tilde u), \pi_N(\tilde v)] + 
[\pi_N(\tilde u),  \tau (\pi_N(\tilde v))]  + [ \tau (\pi_N(\tilde u)), \pi_N(\tilde v)].
\end{align*}

Let us first prove Part (i). By Eq. (\ref{eq5.9}), we thus have
\begin{align*}
&\alpha_{(\fg,\omega,\fa,N')}(u, v)  =  \omega(\pi_\fa'([\pi'_{N'}(\tilde u), \pi'_{N'}(\tilde v)], \cdot ) =  \omega(\pi_\fa([\pi'_{N'}(\tilde u), \pi'_{N'}(\tilde v)] )-  \tau ( [\pi'_{N'}(\tilde u), \pi'_{N'}(\tilde v)]  ), \cdot )\\[2mm]
 & =  \omega(\pi_\fa( [\pi_N(\tilde u), \pi_N(\tilde v)] + 
[\pi_N(\tilde u),  \tau (\pi_N(\tilde v))]  + [ \tau (\pi_N(\tilde u)), \pi_N(\tilde v)]) -  \tau ( [\pi_{N}(\tilde u), \pi_{N}(\tilde v)]  ), \cdot )\\[2mm]
& =  \alpha_{(\fg,\omega,\fa,N)}(u, v) + \omega( 
[\pi_N(\tilde u),  \tau (\pi_N(\tilde v))]  + [ \tau (\pi_N(\tilde u)), \pi_N(\tilde v)]) -  \tau ( [\pi_{N}(\tilde u), \pi_{N}(\tilde v)]  ), \cdot ).
\end{align*}
Since $\fa$ is contained in $\text{Ker} (\tau)$, then $\tau$ defines an element $\bar \tau \in \Hom(\fh, \fa)$ given by $\bar \tau(u):=\tau(\tilde u)$ for all $u\in \fh $,  and therefore
$\sigma = i_\omega  \circ  \bar \tau \in \Hom(\fh, \fh^*)$. Using Eq. (\ref{eq5.9}), we deduce from the above that $\alpha' =
\alpha + \partial_\rho \sigma$. 
Now, using Eq.(\ref{omegataueven}), we get
\[
\begin{array}{lcl}
\sigma(u)(v)- (-1)^{p(u) p(v)} \sigma(v) (u) & = & i_\omega \circ \bar \tau(u)(v)- (-1)^{p(u) p(v)} i_\omega\circ \bar \tau(v)(u)\\[2mm]
& = & \omega( \bar \tau(u), \tilde v) -(-1)^{p(u) p(v)} \omega(\bar \tau(v), \tilde u)\\[2mm]
& =& \omega(\tau( \tilde u), \tilde v)+ \omega (\tilde u, \tau ( \tilde v))=0.
\end{array}
\]

For Part (ii) the proof is similar. In this case, choose $\mu =\Pi \circ  i_\kappa  \circ  \bar \tau \in \Hom(\fh, \Pi(\fh^*))$.
\end{proof}
\subsection{Equivalence classes of $T^*$-extensions and $\Pi T^*$-extensions}
Let $\fh$ be an arbitrary Lie superalgebra. An {\it orthosymplectic} (resp. {\it periplectic}) {\it $T^*$-extension } (resp. {\it $\Pi T^*$-extension}) $(\fg, \omega, \fa)$ of $\fh$ is an~orthosymplectic (resp. a~periplectic) quasi-Frobenius Lie superalgebra $(\fg,\omega)$, together with an extension of Lie superalgebras
\[
0 \rightarrow  \fa  \rightarrow  \fg  \rightarrow \fh  \rightarrow 0 ,
\]
such that the image of $\fa$ in $\fg$ is a Lagrangian ideal of $(\fg, \omega)$.

An {\it isomorphism} of $T^*$-extension or $\Pi T^*$-extension  of the Lie superalgebra $\fh$ is an isomorphism of Lie superalgebras $\Phi: (\fg, \omega) \rightarrow  (\fg', \omega')$ both equipped with orthosymplectic forms (or periplectic forms) such that the
diagram below is commutative:
\begin{equation}\label{diag} \begin{tikzcd}
 0 \arrow{r} & \fa \arrow{r} \arrow[swap]{d}{\Phi|_\fa} & \fg \arrow{r} \arrow{d}{\Phi} &  \fh \arrow{r} \arrow{d}{=} & 0 \\%
0 \arrow{r}&  \fa' \arrow{r} & \fg' \arrow{r} &\fh \arrow{r} & 0
\end{tikzcd}
\end{equation}

\sssbegin[(Isomorphic extensions give the same quotient Lie superalgebra)] {Proposition} [Isomorphic extensions give the same quotient Lie superalgebra] Isomorphic orthosymplectic $T^*$-extension or  periplectic $\Pi T^*$-extensions of $\fh$ give rise to the same associated quotient flat Lie superalgebra $(\fh,\nabla)$. 

\end{Proposition}
\begin{proof}

An immediate consequence of Lemma \ref{Funct1}.
\end{proof}

We now construct for any flat Lie superalgebra $(\fh,\nabla)$ a cohomology,
which describes all ortho-symplectic $T^*$-extensions or periplectic $\Pi T^*$-extensions of $\fh$ with associated flat Lie superalgebra $(\fh, \nabla)$. First, we define Lagrangian $1$- and $2$-cochains on $\fh$ as
\[ \begin{array}{lcl} 
C^1_L(\fh, \fh^*) & :=&  \{\phi \in C^1(\fh, \fh^*) \mid  \phi(u)(v) - (-1)^{p(u)p(v)}\phi(v)(u) = 0\quad  \text{for all } u, v \in \fh\},\\[2mm]
C^2_L(\fh, \fh^*) & := &\{\alpha \in  C^2(\fh, \fh^*) \mid  \alpha \text{ satisfies (\ref{eq5.4})} \}.
\end{array}
\]
\[ \begin{array}{lcl} 
C^1_L(\fh, \Pi(\fh^*)) &: =&  \{\phi \in C^1(\fh, \Pi(\fh^*)) \mid  \phi(u)(v) - (-1)^{p(u)p(v)}\phi(v)(u) = 0 \quad \text{for all } u, v \in \fh\},\\[2mm]
C^2_L(\fh, \Pi(\fh^*)) & := &\{\beta \in  C^2(\fh, \Pi(\fh^*)) \mid  \alpha \text{ satisfies (\ref{eq5.4b})} \}.
\end{array}
\]
Furthermore, let $\rho = \rho^\nabla$ (resp. $\chi = \chi^\nabla$) be the representation of $\fh$ on $\fh^*$ (resp. $\Pi(\fh^*)$) associated to $\nabla$, as
defined in (\ref{eq5.1}) (resp.  (\ref{eq5.1b})). Denote by $\partial_\rho = \partial^i_\rho$ (resp. $\partial_\chi = \partial^i_\chi$) the corresponding coboundary operators for cohomology with $\rho$-coefficients (resp. $\chi$-coefficients). 
\sssbegin[(The coboundary image of Lagrangian cochains)]{Lemma}[The coboundary image of Lagrangian cochains] \textup{(}i\textup{)} The coboundary operator $\partial_\rho : C^1(\fh, \fh^*) \rightarrow C^2(\fh, \fh^*)$ maps the subspace $C^1_L(\fh, \fh^*)$ to $C^2_L(\fh, \fh^*) \cap Z^2 (\fh, \fh^*)$.

\textup{(}ii\textup{)} The coboundary operator $\partial_\chi : C^1(\fh, \Pi(\fh^*)) \rightarrow C^2(\fh, \Pi(\fh^*))$ maps the subspace $C^1_L(\fh, \Pi(\fh^*))$ to $C^2_L(\fh, \Pi(\fh^*)) \cap Z^2(\fh, \Pi(\fh^*))$.
\end{Lemma}
\begin{proof}
Let us prove Part (i). Let $\phi \in C^1(\fh, \fh^*)$, and let $u,v,w\in \fh$. We have 
\[
\begin{array}{lcl}
(\partial_\rho \phi)(u,v) & = &(-1)^{p(u)p(\phi)}\rho(u)\cdot \phi(v)- (-1)^{p(v) (p(u) + p(\phi))} \rho(v)\cdot \phi(u) - \phi ([u,v])\\[2mm]
& = &  - (-1)^{p(u) p(v)} \phi(v)\circ \nabla_u  + \phi(u)\circ \nabla_v - \phi ([u,v]).
\end{array}
\]
Let us check that condition  (\ref{eq5.4}) is satisfied. Indeed, for $u,v,w\in \fh$ we have
\begin{align*}
& (-1)^{p(u)p(w)}(\partial_\rho \phi)(u,v)(w) +  \circlearrowleft (u,v,w)\\[2mm]
&=(-1)^{p(u)p(w)} ( - (-1)^{p(u) p(v)} \phi(v)\circ \nabla_u (w) + \phi(u)\circ \nabla_v(w) - \phi ([u,v])(w))
+  \circlearrowleft (u,v,w).
\end{align*}
Now, using the fact that the connection is flat, the equation above becomes
\begin{align*}
& (-1)^{p(u)p(w)}(\partial_\rho \phi)(u,v)(w) +  \circlearrowleft (u,v,w)\\[2mm]
&= (-1)^{p(u)p(v)} (\phi(v)([w,u])-(-1)^{p(v)(p(w)+p(u))}\phi ([w,u])(v))+  \circlearrowleft (u,v,w)=0
\end{align*}
This follows that $ \partial_\rho \phi \in C^2_L(\fh, \fh^*)$.
\end{proof}
Let $Z^2_L(\fh, \fh^*):=C^2_L(\fh, \fh^*) \cap Z^2(\fh, \fh^*)$ and $Z^2_L(\fh, \Pi(\fh^*)):=C^2_L(\fh, \Pi(\fh^*)) \cap Z^2(\fh, \Pi(\fh^*))$ denote the space of Lagrangian cocycles. We define the Lagrangian extension cohomology of  the flat Lie algebra $(\fh, \nabla)$ by 
\[
\begin{array}{rcl}
H^2_L(\fh, \fh^*) & := & \displaystyle \frac{Z^2_L(\fh, \fh^*)}{\partial_\rho C^1_L(\fh , \fh^*)},\\[4mm]
H^2_L(\fh, \Pi(\fh^*)) & := & \displaystyle \frac{Z^2_L(\fh, \Pi(\fh^*))}{\partial_\chi C^1_L(\fh , \Pi(\fh^*))}.
\end{array}
\]
\sssbegin[(A link between extensions and cohomology classes)]{Theorem}[A link between extensions and cohomology classes]
Every orthosymplectic $T^*$-extension or  periplectic $\Pi T^*$-extensions $(\fg, \omega, \fa)$ of the flat Lie algebra $(\fh, \nabla)$ gives rise to a characteristic extension class
\[
\begin{array}{ll}
[\alpha] \in H^2_{\ev,L}(\fh, \fh^*) & \text{\textup{(}if $\omega$ is orthosymplectic\textup{)}} ,\\[2mm]
[\beta] \in H^2_{\ev, L}(\fh, \Pi(\fh^*)) & \text{\textup{(}if $\omega$ is periplectic\textup{)}}.
\end{array}
\]
Two extensions of $(\fh, \nabla)$ are isomorphic if and only if they have the same extension class in either $H^2_L(\fh, \fh^*)$ or $H^2_L(\fh, \Pi(\fh^*))$.
\end{Theorem}
\begin{proof}
Let $(\fg, \omega, \fa)$ be a Lagrangian extension of $(\fh, \nabla)$, where $\omega$ is orthosymplectic (resp. periplectic). Choose a strong polarization $(\fg, \omega, \fa,N)$. By virtue of Theorem \ref{Tstar}, $(\fg, \omega, \fa,N)$ is a Lagrangian extension by means of $\alpha\in Z^2_{\ev, L}(\fh, \fh^*)$ (resp. $\beta\in Z^2_{\ev, L}(\fh, \Pi(\fh^*))$) if $p(\omega)=\ev$ (resp. $p(\omega)=\od$). Its cohomology class $[\alpha]$ (resp. $[\beta]$) is independent of $N$ by Lemma \ref{classindN}. 

Now suppose $\Phi: (\fg, \omega, \fa) \rightarrow (\fg', \omega',\fa')$ is an isomorphism of $\fh$ (see, diagram (\ref{diag})). Choose a complementary Lagrangian subspace $N$ for $(\fg, \omega, \fa)$. Then, $\Phi:(\fg, \omega, \fa, N)\rightarrow (\fg', \omega', \fa', N'=\Phi(N))$ is an isomorphism of strongly polarized orthosymplectic (resp.  periplectic) Lie superalgebras. By Lemma \ref{Funct2}, we have $\alpha_{\fg', \omega', \fa', N'}=(\Phi_\fh)_* \alpha_{\fg, \omega, \fa, N}$ (resp. $\beta_{\fg', \omega', \fa', N'}=(\Phi_\fh)_* \beta_{\fg, \omega, \fa, N}$), since by assumption $\Phi |_\fh=\id_\fh$. This shows that isomorphic extensions over $\fh$ have the same cohomology class. 

It remains to show that two extensions of $\fh$ with the same class are isomorphic. By Theorem \ref{Tstar}, it is enough to show that any two strongly polarized orthosymplectic (resp. periplectic) Lie superalgebras give rise to isomorphic extensions of $\fh$ if the corresponding classes $\alpha$ and $\alpha'$ if $p(\omega)=\ev$ (resp. $\beta$ and $\beta'$ if $p(\omega)=\od$) satisfy $[\alpha]=[\alpha']\in H^2_{\ev, L}(\fh, \fh^*)$ (resp. $[\beta]=[\beta']\in H^2_{\ev, L}(\fh, \Pi(\fh^*))$). That is $\alpha'=\alpha-\partial_\rho \sigma$ (resp. $\beta'=\beta-\partial_\chi \mu$).  The map $(u, \xi)\mapsto (u, \xi + \sigma(u))$ (resp. $(u, \xi)\mapsto (u, \xi + \mu(u))$) is an ismorphism. 
\end{proof}
\sssbegin[(A link between ordinary cohomology and Lagrangian cohomology)]{Remark}[A link between ordinary cohomology and Lagrangian cohomology]
The kernel of the map $H^2_L(\fg, \fg^*) \rightarrow H^2_\rho(\fg, \fg^*)$ (resp. $H^2_L(\fg, \Pi(\fg^*)) \rightarrow H^2_\chi(\fg, \Pi(\fg^*))$) is given by
\[
\kappa_L=\frac{B^2_\rho (\fg, \fg^*) \cap Z^2_L(\fg, \fg^*)}{B^2_L(\fg, \fg^*)}, \left ( \text{resp. } \kappa_L=\frac{B^2_\chi (\fg, \Pi(\fg^*)) \cap Z^2_L(\fg, \Pi(\fg^*))}{B^2_L(\fg, \Pi(\fg^*))}\right),
\]
where 
\[
B^2_\rho(\fg, \fg^*)=\{ \partial^1_\rho \lambda \mid \lambda \in \Hom(\fg, \fg^*)\} \text{ (resp. } B^2_\rho(\fg, \Pi(\fg^*))=\{ \partial^1_\rho \lambda \mid \lambda \in \Hom(\fg, \Pi(\fg^*))\})
\] 
is the set of ordinary two-coboundaries with $\rho$-coefficients (resp. $\chi$-coefficient) and 
\[
B^2_L(\fg, \fg^*)=\{ \partial^1_\rho \lambda \mid \lambda \in C^1_L(\fg, \fg^*)\} \text{ (resp. } B^2_L(\fg, \Pi(\fg^*))=\{ \partial^1_\rho \lambda \mid \lambda \in C^1_L(\fg, \Pi(\fg^*))\})
\]is the set of two-coboundaries for Lagrangian extension cohomology.
\end{Remark}
\section{Examples of double, $\Pi T^*$- and $T^*$-extensions}\label{Examples}
\ssbegin[(Examples of ${\mathscr D}_\od$-extensions)] {Proposition}[Examples of ${\mathscr D}_\od$-extensions] \label{centdame}  \textup{(}i\textup{)} 
The Lie superalgebra $C^1_{1/2}+A$ -- see Table \ref{tab2} -- cannot be obtained as a double extension even though it has a center. 

\textup{(}ii\textup{)} The Lie superalgebra $(2A_{1,1}+2A)^3_p$ \textup{(}for $p=1/2$\textup{)} -- see Table \ref{tab2} -- can be obtained as a double ${\mathscr D}_\od$-extension of the Abelian Lie superalgebra $\fa:=\mathbb R^{1|1}$ spanned by $\{e_2, e_3\}$ with \textup{(}for the notations, see Theorem \ref{MainThO}\textup{)}:
\[
\omega_\fa=e_2\wedge e_3, \quad x=e_1, \quad e=e_4, \quad a_0=\frac{1}{2}e_2, \quad  {\mathscr D}(e_3)=\frac{1}{2}e_2, \quad {\mathscr D}(e_2)=0, \quad  {\mathscr D}^* \equiv {\mathscr D}.
\]

\end{Proposition}
\begin{proof}

Let us prove Part (i). A direct computation shows that the center of $\fg:=C^1_{1/2}+A$ is $\Span\{e_4\}$. Observe that $\omega_\fg(e_4, e_4)\not =0$. Now, the dual to the center $x=e_4$ with respect to the form $\omega_\fg$ is $x^*=e_3+e_4$. If the Lie superalgebra is a double extension of a Lie sub-superalgebra $\fa$, then $\fa$ would be spanned by $e_1$ and $e_2$. Now the condition $[x^*, x^*]_\fg=2a_0$ implies that $a_0=\frac{1}{2} e_2$. On the other hand, the equation $ [x^*, a]_\fg={\mathscr D}(a)+\omega_\fa(a,a_0)x$ 
implies that ${\mathscr D}(e_1)=-\frac{1}{2}(e_3+e_4)\not \in \fa$ which is a contradiction. 

Let us prove Part (ii). For short, let us write $\fg:=(2A_{1,1}+2A)^3_p$, where $p=1/2$. A direct computation shows that ${\mathscr D}^2=\ad_{a_0}=0$. The map ${\mathscr D}$ is obviously a derivation since the Lie superalgebra $\fa$ is abelian.  Besides, the two conditions ${\mathscr D}(a_0)={\mathscr D}^*(a_0)=0$ are obviously satisfied, by construction of ${\mathscr D}$ and ${\mathscr D}^*$. 
 Now, the map 
\[
(a,b)\mapsto\omega_\fa( {\mathscr D} \circ {\mathscr D}(a) +(1+(-1)^{p(a)}) {\mathscr D}^*\circ {\mathscr D}(a)- {\mathscr D}^* \circ {\mathscr D}^*(a),b)
\]
is not only a coboundary but actually identically zero. Moreover, $\omega_\fa(a_0, [a,b]_\fa)=0$ for all $a,b\in \fa$ since $\fa$ is abelian. It remains to check the Lie brackets.  We have
\[
[e_3, e_3]_\fg=[e_3, e_3]_\fa+ ( \omega_\fa( \mathscr{D} (e_3),e_3)- \omega_\fa(e_3, \mathscr{D} (e_3)) )x= \left (\frac{1}{2}  \omega_\fa(e_2,e_3)- \frac{1}{2}\omega_\fa(e_3, e_2) \right) e_1=e_1,
\]
\[
[e_3, e_4]_\fg={\mathscr D}(e_3)+ \omega_\fa(a_0, e_3)x= \frac{1}{2}e_2+ \frac{1}{2} \omega_\fa(e_2, e_3) e_1=\frac{1}{2}(e_1+e_2).\qed
\]
\noqed
\end{proof}

\ssbegin[(Examples of $T^*$- and $\Pi T^*$-extensions)]{Proposition} [Examples of  $T^*$- and $\Pi T^*$-extensions] \textup{(}i\textup{)} The Lie superalgebra $D^5$ -- see Table \ref{tab1}  -- is a $\Pi T^*$-extension of the Lie superalgebra $\fh$ spanned by $e_1\mid e_4$ with the bracket $[e_1, e_4]=e_4$, and where $\Pi(\fh^*)$ is spanned by $e_2=\Pi(e_4^*) \mid e_3=\Pi(e_1^*)$. The $2$-form $\beta \equiv 0$ and the connection $\nabla$ on $\fh$ is given as follows:
\[
\nabla_{e_4}(e_1)=e_4, \quad \nabla_{e_1}(e_1)=e_1, \quad \nabla_{e_1}(e_4)=2 e_4, \quad \nabla_{e_4}(e_4)=0.
\]
\textup{(}ii\textup{)} The Lie superalgebra $C_1^1+A$ -- see Table \ref{tab2} for its structure -- is a $T^*$-extension of the abelian Lie superalgebra $\fh$ spanned by $e_1 \mid X:=e_3-2 e_4$ where the $2$-form $\alpha$ is given by
\[
\alpha=  e_1^*\otimes X \wedge X+X^*\otimes e_1 \wedge X.
\]
The dual space $\fh^*$ is spanned by $e_2=\frac{1}{2} e_1^* \mid e_3 =X^*$ and the connection $\nabla$ on $\fh$ is given by
\[
\nabla_{e_1}(e_1)=- e_1, \quad \nabla_{e_1}(X)=-X, \quad \nabla_X(e_1)=-X, \quad \nabla_X(X)=0.
\]
\textup{(}iii\textup{)} The Lie superalgebra $D^7_{-1,q}$ ($q\leq -1$) -- see Table \ref{tab1} -- can be obtained as a $\Pi T^*$-extension of the Lie superalgebra $\fh$ spanned by $e_1, e_2$ with the bracket $[e_1, e_2]=e_2$, and where the dual space $\Pi(\fh^*)$ is spanned by $e_3= \Pi(e_1^*)+ \Pi(e_2^*), e_4= \Pi(e_1^*)$. The $2$-form  $\beta \equiv0$ and the connection $\nabla$ on $\fh$ is given by 
\[
\nabla_{e_1}(e_2)= e_2, \quad \nabla_{e_1}(e_1)= -q e_1+(q+1)e_2, \quad \nabla_{e_2}(e_1)=\nabla_{e_2}(e_2)=0.
\]

\end{Proposition}
\begin{proof}

Let us just prove Part (ii) and skip Part (i) and (iii) as the computation is similar. For short, we put $\fk=\fh\oplus \Pi (\fh^*)$. Let us first show that the connection is torsion-free. Indeed,
\[
\nabla_{e_1}(X)-\nabla_X(e_1)-[e_1, X]=-X-(-X)-0=0.
\]
Let us check the flatness. Indeed, 
\[
\nabla_{e_1}\circ \nabla_X- \nabla_X\circ \nabla_{e_1}- \nabla_{[e_1, X]}=\nabla_{e_1}\circ \nabla_X- \nabla_X\circ \nabla_{e_1}.
\]
Evaluating this expression at $X$ gives zero. Let us evaluate it at $e_1$, we get
\[
\nabla_{e_1}\circ \nabla_X(e_1)- \nabla_X\circ \nabla_{e_1}(e_1)=\nabla_{e_1} (-X) - \nabla_X (-e_1)=X- X=0.
\]
Let us check the Lie bracket. We have
\[
[e_4, e_3]_\fk= [ \frac{1}{2} (e_3-X), e_3]_\fk=  - \frac{1}{2} [X, X^*]_\fk =  -\frac{1}{2} X^* \circ \nabla_X=\frac{1}{2} e_1^*=e_2, 
\]
\[
[e_1, e_3]_\fk=[e_1, X^*]_\fk=- X^* \circ \nabla_{e_1}= X^*=e_3, \quad 
[e_1, e_2]_\fk= [e_1, \frac{1}{2}e_1^*]_\fk=-\frac{1}{2} e_1^* \circ \nabla_{e_1}=\frac{1}{2}e_1^*=e_2,
\]
\begin{align*}
&[e_4, e_4]_\fk= [\frac{1}{2} (e_3-X), \frac{1}{2} (e_3-X)]_\fk=\frac{1}{4} \left ( [X^* -X, X^* -X]_\fk  \right )=\frac{1}{4} (- 2[X,X^*]_\fk + [X,X]_\fk)\\
&=  \frac{1}{4} ( -2 X^* \circ \nabla_X - 2 e_1^*) =   \frac{1}{4} ( 2 e_1^*  -2 e_1^*) =0.\qed
\end{align*}
\noqed
\end{proof}

\ssbegin[(Filiforms as double and $\Pi T^*$-extensions)]{Proposition}[Filiforms as double and $\Pi T^*$-extensions]
\textup{(}i\textup{)} The filiform Lie superalgebra $L^{n,n}$ can be obtained as a $\Pi T^*$-extension of the filiform Lie algebra $L^n$ where the connection $\nabla$ on $L^n$ is given by
\[
\nabla_{X_1}(X_i)=X_{i+1}\quad \text{ for $1\leq i \leq n-1$ \quad and \quad $\nabla_{X_j}(X_i)=0 $ \quad  for all $j>1$}.
\]
The odd part $\Pi({(L^n})^*)$ is given by 
\[
Y_i=(-1)^i \sum_{j=1}^{n-i+1} \Pi X_j^*\quad \text{ for $i=1,\ldots,n$.}
\]
\textup{(}ii\textup{)} For $n$ even and $m$ odd,  the filiform Lie superalgebra $L^{n,m}$ can be obtained as a ${\mathscr D}_\ev$-extension of the abelian Lie  superalgebra ${\mathbb R}^{n-2,1}$, followed by successive ${\mathscr D}_\od$-extensions.
\begin{proof}

For Part (i), let us first check that the connection is torsion-free. Indeed, 
\[
\nabla_{X_1}X_i- \nabla_{X_i}X_1-[X_1, X_i]=X_{i+1}-X_{i+1}=0.
\]
Now, the connection is flat since $\nabla_{X_i}\circ \nabla_{X_j}=0$ and $\nabla_{[X_i,X_j]}(X_k)=0$ for all $i,j,k=1,\ldots,n$. 
Let us now check the Lie brackets. Indeed, for all $i=1,\ldots, n-1$, we have 
\[
[X_1, Y_{i}]=[X_1,(-1)^i \sum_{j=1}^{n-i+1} \Pi X_j^* ]= -(-1)^i \sum_{j=1}^{n-i+1} \Pi X_j^* \circ \nabla_{X_1}=(-1)^{i+1} \sum_{j=1}^{n-i} \Pi X_j^* =Y_{i+1}.
\]
Also,
\[
[X_1, Y_{n}]=[X_1, (-1)^n \Pi(X_1^*)]=-(-1)^n \Pi(X_1^*)\circ \nabla_{X_1}=0.
\]

Let us prove Part (ii). First, let us observe that $L^{n,m}$ is a ${\mathscr D}_\od$-extension of $L^{n,m-2}$ as follows (for the notations, see Theorem \ref{MainThO}):
\[
a_0=0, \quad {\mathscr D}(X_1)=-(-1)^{\frac{m+1}{2}}Y_2, \quad  {\mathscr D}(a)=0 \text{ for all $a\not =X_1$}, \quad x^*= (-1)^{\frac{m+1}{2}} Y_1, \quad x=Y_m.
\]
Now, all the requirements are satisfied: ${\mathscr D}^2=0$  and the map $(a,b)\mapsto \omega({\mathscr D}\circ {\mathscr D} (a)-{\mathscr D}^*\circ {\mathscr D}^*(a),b)$ is not only a coboundary but is identically zero.

By induction, we deduce that $L^{n,m}$ is a successive ${\mathscr D}_\od$-extensions of the filiform Lie superalgebra $L^{n,1}$.

Let us now show that $L^{n,1}$ is a ${\mathscr D}_\ev$-extension of the Abelian Lie superalgebra ${\mathbb R}^{n-2|1}$. Denote by $X_1,\ldots, X_n|Y_1$ the generators of $L^{n,1}$. Let us choose (for the notations, see Theorem \ref{MainTh}): 
 \[
\lambda=Z_\Omega=0,\quad x=X_n, \quad x^*= X_1, \quad  {\mathscr D}(X_i)=X_{i+1} \text{ for $i=2,\ldots, n-2$}, \quad  {\mathscr D}(Y_1)= {\mathscr D}(X_{n-1})=0.
\]
Now, the map 
\[
(a,b)\mapsto\omega_\fa( {\mathscr D} \circ {\mathscr D}(a) +2 {\mathscr D}^*\circ {\mathscr D}(a)- {\mathscr D}^* \circ {\mathscr D}^*(a)+ \lambda({\mathscr D}+{\mathscr D^*})(a),b)
\]
is not only a coboundary, but it is identically zero. \end{proof}
\end{Proposition}
{\bf Acknowledgments.} We would like to thank S.  Benayadi for several stimulating and fruitful discussions. 



\begin{thebibliography}{999}

\bibitem[ABB]{ABB} Albuquerque H., Barreiro E. and Benayadi S., Quadratic Lie superalgebras with a reductive even part. J. Pure Appl. Algebra, {\bf 213} (2009), 724-731.

\bibitem[ABBQ]{ABBQ} Albuquerque H., Barreiro E. and Benayadi S., Odd quadratic Lie superalgebras. J. of Geometry and Physics, {\bf 60} (2010), 230-250.

\bibitem[Ba]{Ba} Backhouse N., A classification of four-dimensional Lie superalgebras, Journal of Mathematics Physics, {\bf 19}, 2400, (1987).

\bibitem[BaBe]{BaBe} Bajo I. and Benayadi S., Abelian para-Kähler structures on Lie algebras, Differential Geometry and its Applications, {\bf 29} (2011), 160--173.

\bibitem[BeB]{BeB} Benamor H. and Benayadi S., Double extension of quadratic Lie superalgebras. Comm. Algebra, {\bf 27}, No. 1 (1999), 67-88.

\bibitem[B]{B} Benayadi S., Quadratic Lie superalgebras with completely reductive action of the even part on the odd part. J. of Algebra, {\bf 223} (2000), 344-366.

\bibitem[BeBou]{BeBou} Benayadi S. and Bouarroudj S., Double extensions of Lie superalgebras in characteristic $2$ with non-degenerate invariant supersymmetric bilinear forms, J. of Algebra, {\bf 510} (2018), 141--179; arXiv:1707.00970.

\bibitem[BB]{BB} Benayadi S. and Bouarroudj S., Manin triples and non-degenerate anti-symmetric bilinear forms on Lie superalgebras in characteristic $2$, (2021); https://doi.org/10.1016/j.jalgebra.2022.09.019; \texttt{arXiv:2110.05141 } 

\bibitem[BBH]{BBH} Benayadi S., Bouarroudj S., Hajli M.,  Double extensions of restricted Lie (super)algebras,  Arnold. Math. J.  {\bf 6} (2020),  231 -- 269; \texttt{arXiv:1810.03086}

\bibitem[BC]{BC} Baues O. and Cort\'es V., Symplectic Lie groups I-III, Symplectic Reduction, Lagrangian extensions, and existence of Lagrangian normal subgroups. Ast\'etrique {\bf 379}, Soci\'et\'e math\'ematique de France, 2016.

\bibitem[BYC]{BYC} Bon-Yao-Chu, Symplectic homogeneous spaces, Trans. Am. Math. Soc. {\bf 197}
(1974), 145--159.

\bibitem[Bor]{Bor} Bordemann M., Nondegenerate invariant bilinear forms on nonassociative algebras. Acta Math. Univ. Comenianae, Vol. LXVI, 2 (1997), 151--201.

\bibitem[BGKN]{BGKN} Bordemann M, Gomez J.R., Khakimdjanov, Yu. and Navarro R.M., Some deformations of nilpotent Lie superalgebras. J. Geom. Phys. 57 (2007), no. 5, 1391--1403.





\bibitem[BGL]{BGL} Bouarroudj S., Grozman P., Leites D., Classification of finite-dimensional modular Lie superalgebras with indecomposable Cartan matrix. Symmetry, Integrability and Geometry: Methods and Applications (SIGMA), 5 (2009), 060, 63 pages; arXiv:math.RT/0710.5149

\bibitem[BGLLS]{BGLLS}
Bouarroudj S., Grozman P., Lebedev A., Leites D. and Shchepochkina I., Simple vectorial Lie algebras in characteristic
$2$ and their superizations. Symmetry, Integrability and Geometry: Methods and Applications (SIGMA) {\bf 16} (2020), 089, 101 pages; \texttt{arXiv:1510.07255}

\bibitem[BN]{BN} Bouarroudj S, Navarro R. M.,  Cohomologically rigid solvable Lie superalgebras with model filiform and model nilpotent nilradical, Communications in Algebra (2020); https://doi.org/10.1080/00927872.2021.1936541


\bibitem[E]{E} Elashvili A. G., Frobenius Lie algebras, Functional Analysis and Its Applications, (1982), 16(4), 326--328.

\bibitem[F]{F} Fuchs D. B., {\it Cohomology of Infinite-Dimensional Lie Algebras}, Contemp. Soviet.
Math., (Consultants Bureau, New York, 1986).

\bibitem[O]{O} Ooms I. O., On Frobenius Lie algebras, Communications in Algebra,  (1980),  8(1),  13--52.


\bibitem[H]{H} Hegazi A., Classification of Nilpotent Lie superalgebras of dimension five. I, International Journal of Theoretical Physics, Vol. {\bf 38}, No. 6, (1999),  1735--1740.


\bibitem[K]{K} Kac V., Lie superalgebras, Adv. Math. {\bf 26} (1977), 8--96.



\bibitem[L]{L} Leites D., Lie superalgebra cohomology. Funct. Analysis and Its Applications, (1975), 9:4, 340--341.

\bibitem[LSoS]{LSoS}
Leites D. (ed.) \textit{Seminar on supersymmetry v. $1$. Algebra and
Calculus: Main chapters}, (J.~Bernstein, D.~Leites, V.~Molotkov,
V.~Shander), MCCME, Moscow, 2012, 410 pp (in Russian; a~version in
English is in preparation but available for perusal).


\bibitem[MR1]{MR1} Medina A. and Revoy P., Alg\`ebres de Lie et produit scalaire invariant. Ann. Scient. \'Ec. Norm. Sup., 4 s\'erie, {\bf 18} (1985), 553--561.

\bibitem[MR2]{MR2} Medina A. and Revoy P.,  Groupes de Lie \`a structure symplectique invariante, in P. Dazord and A. Weinstein (Eds.)``Symplectic Geometry, Groupoids and Integrable Systems, S\'eminaire Sud-Rhodanien de G\'eom\'etrie'', Mathematical Sciences Research Institute Publications, pp. 247--266, Springer-Verlag, Berlin New York, 1991.

\bibitem[MD]{MD} Medina A. and Dardi\'e J.-M., Double extension symplectique d'un groupe de Lie symplectique, Advances in Mathematics, {\bf 117} (1996), 208--227.

\bibitem[V]{V} Vergne M., Cohomologie des alg\`ebres de Lie nilpotentes. Application \`a l’\'etude de la vari\'et\'e des alg\`ebres de Lie nilpotentes, Bull. Sot. Math. France {\bf 98} (1970), 81--116.






\end{thebibliography}
\end{document}